\documentclass[12pt]{article}
\usepackage{a4, amssymb, latexsym, amsmath, exscale}

\newsavebox{\ttt}
\sbox{\ttt}{}
\newcommand{\startsection}[1]
    {\section[#1]{#1}
    \sbox{\ttt}{\thesection\ \ \textsc{#1}}
    \thispagestyle{plain}
}

\newenvironment{evlist}[2]{
\begin{list}{}{
\setlength{\topsep}{0.5ex plus0.2ex minus0.1ex} 
\setlength{\leftmargin}{#1}
\setlength{\itemsep}{#2 plus0.2ex}
\setlength{\parsep}{0ex plus0.2ex} }}
{\end{list}}

\newcommand{\makespace}[2]{\rule[-#1mm]{0mm}{#2mm}}
\newcommand{\Nat}{\mathbb{N}}
\newcommand{\Int}{\mathbb{Z}}
\newcommand{\id}{\mathrm{id}}
\newcommand{\definition}[1]{\textit{#1}}
\newcommand{\proof}{{\textit{Proof}\enspace}}
\newcommand{\Self}[1]{\mathrm{T}_{#1}}

\newcommand{\eop}{\ \vbox{\hrule
                       \hbox{\vrule
                             \hskip 6pt
                             \vrule height 6pt width 0pt
                             \vrule}%
                       \hrule}%
                     \vspace{\medskipamount}
                }

\newtheorem{lemma}{Lemma}[section]
\newtheorem{proposition}{Proposition}[section]
\newtheorem{theorem}{Theorem}[section]
\newcommand{\prop}{\mathsf{P}}

\newcommand{\ssets}[1]{\#_{\mathcal{P}}(#1)}

\newcommand\mythicklines{}
\newcommand\sput[3]{\put(#1,#2){$\scriptstyle{#3}$}}

\newlength{\graphicthick}
\newlength{\graphicmid}
\newlength{\graphicthin}

\begin{document}

\begin{titlepage}

\begin{center}
\phantom{q}
\vspace{60pt}

{\huge Finite Sets and Counting}

\bigskip
{\LARGE Chris Preston}

\bigskip
{\large June 2010}
\end{center}

\bigskip
\bigskip
\bigskip
\bigskip
\begin{quote}

We start by presenting a theory of finite sets using the approach which
is essentially that taken by Whitehead and Russell in 
\textit{Principia Mathematica}, and which does not involve the natural numbers 
(or any other infinite set).
This theory is then applied to prove results about structures which, like the natural
numbers, satisfy the principle of mathematical induction, but do not necessarily satisfy
the remaining Peano axioms.

\end{quote}

\end{titlepage}

\thispagestyle{empty}

\addtocontents{toc}{\vskip 20pt}
\addtolength{\parskip}{-5pt}
\tableofcontents
\addtolength{\parskip}{5pt}

\startsection{Introduction}

\label{intro}

The notion of what it means for a set to be finite is usually presented within the framework of the natural
numbers $\Nat = \{0,1,\ldots\,\}$ (and for our purposes it is convenient to regard $0$ as being a natural number, 
so $0$ can be used to count the number of elements in the empty set). The standard approach is to define a set 
$A$ to be finite if there exists $n \in \Nat$ and a bijective mapping $h : [n] \to A$, where $[0] = \varnothing$ 
and $[n] = \{0,1,\ldots,n-1\}$ for each $n \in \Nat$ with $n \ne 0$. If this approach appears to be rather 
straightforward then it is probably because the various facts about $\Nat$ which are needed to define the set 
$\{0,1,\ldots,m\}$ have been taken for granted.

We are going to define finite sets within a framework which does not involve the natural numbers. (The definition 
will of course meet the obvious requirement of being equivalent to the usual one.) One reason for doing this is 
that it seems somewhat strange to have to define finite sets in terms of the infinite set $\Nat$, and in fact
we don't want to assume in general anything about the existence of such sets. The definition of being finite to 
be used here is essentially that employed by Whitehead  and Russell in \textit{Principia Mathematica} \cite{wr}. 
There are similar definitions which also appeared in the first decades of the previous century and excellent 
treatments of this topic can be found in Levy \cite{levy} and Suppes \cite{suppes}.

Let us start by discussing a simple situation involving counting. This is done very informally, taking for 
granted all the standard properties of the natural numbers. The aim is to make a plausible transition from the 
old to the new definition of being finite.

Suppose there is a box containing table-tennis balls, and we want to know how many balls there are in the box.
Let us count them by taking them out of the box one at a time and marking them, say with a red felt marker, 
marking the first ball taken out with a $1$, the second with a $2$ and so on. Eventually we take the last ball 
out  and mark it, and the number we write on this last ball tells us how many balls there were in the box. 

Now repeat the above procedure: Put the balls back in the box, mix them up and then take them out again one 
by one, this time marking them with a green felt marker. Then the last ball taken out will be marked with the 
same number (now in green) that appears on the last ball that was taken out when the red marker was used. 

Each ball now has two numbers on it, one written in red and one in green, and in general there will be no 
discernible relationship between the red and green numbers to be found on any ball. However, this is not 
important; what is important is that the number written on the last ball in red is the same as the number written 
on the last ball in green. This is one of the reasons that counting makes sense.

Consider the following simple example:

\begin{center}
\begin{tabular}{|r||c|c|c|c|c|c|c|c|c|}\hline
 \makespace{2}{7}     red number: &  1 & 2 & 3 & 4 & 5 & 6 & 7 & 8 & 9 \\ \hline
 \makespace{2}{7}   green number: &  8 & 9 & 3 & 2 & 4 & 6 & 5 & 1 & 7 \\ \hline
\end{tabular}
\end{center}

The number written on the last ball in red is $9$, which is the same as the number written on the last ball in 
green. Note that the last ball taken out when using the red marker was the $7$th ball taken out using the green 
marker, and the last ball taken out when using the green marker was the $2$nd ball taken out using the red 
marker. For $9$ balls there are   $9! = 362880$ different ways of taking the balls one at a time out of the box. 
However, they all give the same result as far as determining the number of balls is concerned: The last ball taken 
out and counted will always be assigned the number $9$.

The above discussion illustrates a \definition{basic counting principle} which says that if a set of objects 
can be counted then it doesn't matter in which order the counting is carried out. Most people's experience of 
counting things leads them at some point to accept the principle as something which is true about the world.
Thereafter they are usually not even aware of the role which it plays.

Note that the basic counting principle only applies to sets of objects which can be counted in the sense that
the process of removing the objects from the set one at a time will eventually terminate. Of course, it is not
really clear what `eventually terminates' means, but intuitively it characterises the finite sets. Let us call a set 
with this rather vague property \definition{$\Nat$-finite} in order to distinguish it from the new definition of 
being finite to be given below.

We next reformulate the basic counting principle in a form in which the natural numbers do not occur. Consider 
the set $S$ whose elements are the balls in the box, and suppose the box contains $N$ balls. In order to refer 
to the balls we assume that the numbers made with the red felt marker are still there, and we refer to the ball
with the number $n$ written on it in red by $b_n$.

Let $f : S \to S$ be an injective mapping, so $f(s_1) \ne f(s_2)$ whenever $s_1 \ne s_2$. We can write this 
mapping as a table:

\begin{center}
\begin{tabular}{|c||c|c|c|c|c|c|c|}\hline
\makespace{2}{7} $b_n$ &  $b_1$ & $b_2$ & $b_3$ & $b_4$ & $\cdots$ & $b_{N-1}$ & $b_N$ \\ \hline
\makespace{2}{7} $f(b_n)$ &  $f(b_1)$ & $f(b_2)$ & $f(b_3)$ & $f(b_4)$ & $\cdots$ & $f(b_{N-1})$ & $f(b_N)$\\ \hline
\end{tabular}
\end{center}

An example of such a mapping (with $9$ balls) is the following:

\begin{center}
\begin{tabular}{|c||c|c|c|c|c|c|c|c|c|}\hline
\makespace{2}{7}    $b_n$ &  $b_1$ & $b_2$ & $b_3$ & $b_4$ & $b_5$ & $b_6$ & $b_7$ & $b_8$ & $b_9$ \\ \hline
\makespace{2}{7}   $f(b_n)$ &  $b_8$ & $b_9$ & $b_3$ & $b_2$ & $b_4$ & $b_6$ & $b_5$ & $b_1$ & $b_7$ \\ \hline
\end{tabular}
\end{center}

In this example the ball marked with the red $1$ is mapped to the ball marked with the red $8$, the ball marked 
with the red $2$ is mapped to the ball marked with the red $9$, and so on.

Now do the following: Find the ball $f(b_1)$ and write $1$ on it in blue, then find the ball $f(b_2)$ and write 
$2$ on it in blue. Repeat this until finally finding the ball $f(b_N)$ and writing $N$ on it in blue. Then for 
each number between $1$ and $N$ there is exactly one ball with this number written on it in blue. (No ball can 
have more than one blue number written on it because the mapping $f$ is injective.)

For the above example this means that the ball with the red $8$ also has the blue $1$ on it, the ball with the 
red $9$ also has the blue $2$ on it, and so on:

\begin{center}
\begin{tabular}{|r||c|c|c|c|c|c|c|c|c|}\hline
 \makespace{2}{7}  red number: &  8 & 9 & 3 & 2 & 4 & 6 & 5 & 1 & 7 \\ \hline
 \makespace{2}{7} blue number: &  1 & 2 & 3 & 4 & 5 & 6 & 7 & 8 & 9 \\ \hline
\end{tabular}
\end{center}

Now count the balls in the box using the the blue numbers: First take the ball marked with the blue $1$ out of 
the box, then the ball marked with the blue $2$ and so on. Eventually the ball marked with the blue $N$ will be 
taken out, and then there can be no more balls left in the box: If this were not the case then we could continue 
the counting by taking the remaining balls out in any order and would conclude that there were more than $N$ 
balls originally in the box, which by the basic counting principle is not possible. But if there are no more 
balls left in the box after the ball marked with the blue $N$ is taken out then every ball is marked with a blue 
number. It follows that the mapping $f$ is surjective, since a ball marked with a blue $n$ is the image under $f$
of the ball marked with a red $n$. Thus $f$, being injective and surjective, is bijective.

Let us say that a set $A$ is \definition{tame} if every injective mapping $f : A \to A$ is bijective. We have 
thus seen that the basic counting principle implies that any $\Nat$-finite set is tame.

The converse also holds: Suppose that $A$ is an $\Nat$-finite set which is tame and consider two ways of counting 
out the elements of $A$. The first time we count we mark the elements in red and suppose that the last element 
is marked with the number $M$; the second time we mark the elements in green and suppose that the last element 
is marked with the number $N$. Without loss of generality we can assume that $M \le N$ and so we can define a 
mapping $f : A \to A$ as follows: For each number $m$ between $1$ and $M$ let the image under $f$ of the element 
marked with $m$ in red be the element marked with $m$ in green. Then $f$ is injective and thus bijective, since 
$A$ is tame, which implies that $M = N$. This shows that the basic counting principle is valid for $A$.

We have now established that for $\Nat$-finite sets the basic counting principle is equivalent to being tame, and 
note that being tame has nothing to do with the natural numbers. Let us next consider an important property of 
tame sets which will be used to help motivate the new definition of being finite.

\begin{lemma}\label{lemma_intro_11}
Let $A$ be tame and $b$ be some element not in $A$. Then $A' = A \cup \{b\}$ is also tame.
\end{lemma}

\proof
Consider an injective mapping $f : A' \to A'$; then there are two cases:

1.\enskip 
$f(A) \subset A$. Then the restriction $f_{|A} : A \to A$ of $f$ to $A$ is injective
and hence bijective, since $A$ is tame. If $f(b) \in A$ then $f(a) = f_{|A}(a) = f(b)$ for 
some $a \in A$, since $f_{|A}$ is surjective, which contradicts the fact that $f$ is injective. Thus $f(b) = b$, and 
it follows that $f$ is bijective.

2.\enskip 
$f(A) \not\subset A$. In this case there exists an element $c \in A$ with $f(c) = b$ and, since $f$ is injective, 
we must have $f(a) \in A$ for all $a \in A \setminus \{c\}$ and $f(b) \in A$. This means there is an injective 
mapping $g : A \to A$ defined by letting
\[ g(a) = \left\{ \begin{array}{cl}
                  f(a) &\ \mbox{if}\ a \in A \setminus \{c\}\;,\\
                  f(b)   &\  \mbox{if}\ a = c\;\\
\end{array} \right. \]
and then $g$ is bijective, since $A$ is tame. Therefore $f$ is again bijective.

Thus in both cases $f$ is bijective, and this shows that $A'$ is tame. 
\eop

Note that as a very special case the empty set $\varnothing$ is tame, since there is only one mapping 
$\epsilon : \varnothing \to \varnothing$ which is trivially bijective. For the same reason each singleton set 
$\{a\}$ is tame. This also follows from Lemma~\ref{lemma_intro_11}, since $\varnothing$ is tame and 
$\{a\} = \varnothing \cup \{a\}$.

Let $A$ be an $\Nat$-finite set; then we can apply Lemma~\ref{lemma_intro_11} to show that $A$ is tame without 
appealing to the basic counting principle: Start with the set $A$ and the empty set $\varnothing$; next remove 
an element $a_1$ from $A$ and add it to $\varnothing$ to give the set $\varnothing \cup \{a_1\} = \{a_1\}$;
then remove an element $a_2$ from $A \setminus \{a_1\}$ and add it to $\{a_1\}$ to give the set 
$\{a_1\} \cup \{a_2\} = \{a_1,a_2\}$. Repeat the procedure, at each stage removing an element from what is left 
from $A$ and adding it to the set consisting of the elements of $A$ which have already been removed. This 
process will eventually come to a halt, because at some point all the elements will have been removed from $A$. 
The other set will then consist of all these elements and so will be equal to $A$. As noted above $\varnothing$ 
is tame; thus by Lemma~\ref{lemma_intro_11} $\{a_1\} = \varnothing \cup \{a_1\}$ is tame and so again 
by Lemma~\ref{lemma_intro_11} $\{a_1,a_2\} = \{a_1\} \cup \{a_2\}$ is tame. Repeating this until the original set is 
empty, Lemma~\ref{lemma_intro_11} shows that at each stage the set consisting of those elements which have been
removed from $A$ will be tame, and therefore $A$ itself must have the property.

The set of all subsets of a set $A$ (the \definition{power set of $A$}) will be denoted as usual by
$\mathcal{P}(A)$. The new definition of being finite will involve what is known as an 
\definition{inductive system in $\mathcal{P}(A)$} (or just an \definition{inductive system} if $A$ can be 
determined from the context): This is any subset $\mathcal{K}$ of $\mathcal{P}(A)$ such that 
$\varnothing \in \mathcal{K}$ and for which $B \cup \{a\} \in \mathcal{K}$ for all $B \in \mathcal{K}$, 
$b \in A \setminus B$. In particular, $\mathcal{P}(A)$ itself is always an inductive system. Another example is 
given by Lemma~\ref{lemma_intro_11} which shows that
$\mathcal{K} = \{ B \in \mathcal{P}(A) : \mbox{$B$ is tame} \}$ is an inductive system.

Now consider an $\Nat$-finite set $A$ and let $\mathcal{K}$ be an inductive system. Let $B \in \mathcal{P}(A)$; 
then $B$ is $\Nat$-finite, and the argument used to show that an $\Nat$-finite set is tame shows that 
$B \in \mathcal{K}$. Thus $\mathcal{K} = \mathcal{P}(A)$. In other words, if $A$ is $\Nat$-finite then 
$\mathcal{P}(A)$ is the only inductive system.

On the other hand, for any set $A$ we have that 
$\mathcal{K} = \{ B \in \mathcal{P}(A) : \mbox{$B$ is $\Nat$-finite}\}$ is an inductive system: The empty set 
is trivially $\Nat$-finite, and if $B \in \mathcal{K}$ and $a \in A \setminus B$ then 
$B \cup \{a\} \in \mathcal{K}$. (First count the elements in $B$ and then count $a$.) In particular, if 
$\mathcal{P}(A)$ is the only inductive system then $\mathcal{K} = \mathcal{P}(A)$ and hence $A \in \mathcal{K}$,
i.e., $A$ is $\Nat$-finite.

The above (rather informal) arguments imply that a set $A$ is $\Nat$-finite if and only if $\mathcal{P}(A)$ is 
the only inductive system, which leads to the new definition: A set $A$ will be called \definition{finite} if
$\mathcal{P}(A)$ itself is the only inductive system. As mentioned above, this definition is essentially that 
employed by Whitehead  and Russell \cite{wr}. Note that it has nothing to do with the natural numbers.

\begin{proposition}\label{prop_intro_11}
Every finite set is tame. Thus if $A$ is finite then every injective mapping $f : A \to A$ is bijective.
\end{proposition}

\proof 
Let $A$ be a finite set; as noted above $\mathcal{K} = \{ B \subset A : \mbox{$B$ is tame} \}$ is an inductive 
system and hence $\mathcal{K} = \mathcal{P}(A)$, since $A$ is finite. In particular $A \in \mathcal{K}$, which 
shows that $A$ is tame.
\eop

\begin{lemma}\label{lemma_intro_21}
Let $A$ be finite and $b$ be some element not in $A$. Then $A' = A \cup \{b\}$ is also finite.
\end{lemma}

\proof 
Let $\mathcal{K}'$ be an inductive system in $\mathcal{P}(A')$; then 
$\mathcal{K} = \{ B \in \mathcal{K}' : B \subset A \}$ is clearly an inductive system in $\mathcal{P}(A)$ and 
hence $\mathcal{K} = \mathcal{P}(A)$. Thus $\mathcal{P}(A) \subset \mathcal{K}'$. Moreover, 
$B \cup \{a\} \in \mathcal{K}'$ for all $B \in \mathcal{P}(A)$, since $\mathcal{K}'$ is an inductive system in
$\mathcal{P}(A')$, and therefore $\mathcal{K}' = \mathcal{P}(A')$. It follows that $A'$ is finite.
\eop

Lemma~\ref{lemma_intro_21} shows that sets defined by explicitly giving each of their elements are finite. For 
example, consider the set $A = \{a,b,c,d\}$, where the elements $a,\,b,\,c,\,d$ are all different. Then $\varnothing$ 
is finite implies $\{a\} = \varnothing \cup \{a\}$ is finite implies $\{a,b\} = \{a\} \cup \{b\}$ is finite implies 
$\{a,b,c\} = \{a,b\} \cup \{c\}$ is finite implies $A = \{a,b,c,d\} = \{a,b,c\} \cup \{d\}$ is finite.

In Section~\ref{fs} we establish the properties of finite sets, using the above definition of being finite. There is 
nothing surprising here: Every subset of a finite set is finite, if $A$ and $B$ are finite sets then so are 
$A \cup B$, $A \times B$ and $B^A$ (the set of all mappings from $A$ to $B$) and the power set of a finite set is 
finite. If $A$ is a finite set and $f : A \to A$  a mapping then  $f$ is injective if and only if it is surjective 
(and thus if and only if it is bijective). If $A$ is finite and there exists an injective mapping $f : B \to A$ or a 
surjective mapping $f : A \to B$ then $B$ is also finite. If $A$ and $B$ are finite sets then either there exists an 
injective mapping $f : A \to B$ or an injective mapping $g : B \to A$, and if there exists both an injective mapping 
$f : A \to B$ and an injective mapping $g : B \to A$ then $A \approx B$, where we write $A \approx B$ if there exists
a bijective mapping $h : A \to B$. If $B$ is a subset of a finite set $A$ with $B \approx A$ then $B = A$. The proofs
are mostly very straightforward and they all follow the same pattern (and so tend to become somewhat monotonous).

Having developed a theory of finite sets we turn in Section~\ref{cfs} to the question of what it means to count the 
elements in such a set. Of course, counting is usually associated with the natural numbers, and let us start by 
giving an informal discussion of this case. Let $\textsf{s} : \Nat \to \Nat$ be the successor operation (with 
$\textsf{s}(0) = 1$, $\textsf{s}(1) = 2$ and so on). For each finite set $A$ there is an element $\#(A)$ of $\Nat$ 
which tells us how many elements $A$ contains, and which is usually referred to as the \definition{cardinality of $A$}.
The assignment $A \mapsto \#(A)$ has the properties that $\#(\varnothing) = 0$ and 
$\#(A \cup \{a\}) = \textsf{s}(\#(A))$ whenever $A$ is a finite set and $a$ is an element not in $A$. Moreover, it 
is uniquely determined by these requirements. The importance of the assignment $\#$ is that if $A$ and $B$ are 
finite sets then $\#(B) = \#(A)$ if and only if $B \approx A$.

Note that the existence of $\#$ implies that the basic counting principle introduced above is valid. For example, if 
$A = \{a,b,c\}$ then 
\begin{eqnarray*}
\#(A) &=& \#(\{a,b\} \cup \{c\}) = \textsf{s}(\#(\{a,b\})) 
= \textsf{s}(\#(\{a\}\cup \{b\}))  = \textsf{s}(\textsf{s}(\#(\{a\})))\\  
&=& \textsf{s}(\textsf{s}(\#(\varnothing \cup \{a\})))
= \textsf{s}(\textsf{s}(\textsf{s}(\#(\varnothing))))  
= \textsf{s}(\textsf{s}(\textsf{s}(0))) = 3
\end{eqnarray*}
which corresponds to counting the elements of $A$ by first removing $c$, then $b$ and finally $a$. However, the 
result $3$ does not depend on the order in which the elements are removed.

The triple $(\Nat,\textsf{s},0)$ is a special case of what we call a \definition{counting system}, which is defined 
to be any triple $(X,f,x_0)$ consisting of a set $X$, a mapping $f : X \to X$ and an element $x_0 \in X$. Consider a 
counting system $(X,f,x_0)$ and for each finite set $A$ let $\#(A)$ be an element of $X$. Then the assignment 
$A \mapsto \#(A)$ will be called an \definition{iterator for $(X,f,x_0)$} if $\#(\varnothing) = x_0$ and 
$\#(A \cup \{a\}) = f(\#(A))$ whenever $A$ is a finite set and $a \notin A$. (This means, for example, that 
$\#(\{a\}) = f(x_0)$ for each element $a$, $\#(\{a,b\}) = f(f(x_0))$ for distinct elements $a$ and $b$, and so on.)
In particular, the operation which assigns to each finite set $A$ its cardinality $\#(A)$ defines an iterator for 
$(\Nat,\textsf{s},0)$. In fact, we will see that there exists a unique iterator $\#$ for each counting system 
$(X,f,x_0)$ and it has the property that $\#(B) = \#(A)$ whenever $A$ and $B$ are finite sets with $B \approx A$. 
(These statements will be established in Theorem~\ref{theorem_cfs_11} and Proposition~\ref{prop_cfs_21} respectively.)

The iterator $\#$ for $(\Nat,\textsf{s},0)$ (with $\#(A)$ the cardinality of $A$) has the additional crucial 
property that $B \approx A$ whenever $A$ and $B$ are finite sets with $\#(B) = \#(A)$. If an iterator $\#$ for 
a counting system $(X,f,x_0)$ has this property then we will say that $\#$ is \definition{complete}. In 
general $\#$ will not be complete, as can be seen by looking at the trivial example in which $X$ consist of 
the single element $x_0$, and so $f$ is the mapping with $f(x_0) = x_0$. Then $\#(A) = x_0$ for each finite 
set $A$ and here $\#(B) = \#(A)$ for all finite sets $A$ and $B$. The completeness of the iterator for 
$(\Nat,\textsf{s},0)$ must therefore depend on some special properties of the natural numbers. These come about 
because the counting system $(\Nat,\textsf{s},0)$ is assumed to satisfy what are usually called the Peano axioms 
(although, as was acknowledged by Peano, they were introduced earlier by Dedekind in \cite{dedekind}).

One of these axioms is the \definition{principle of mathematical induction}, which holds for a counting system 
$(X,f,x_0)$ if whenever a proposition $\prop$ is given about the elements of $X$ (meaning for each $x \in X$ 
we have  a proposition $\prop(x)$) such that
\begin{evlist}{30pt}{6pt}
\item[($\diamond$)]  $\prop(x_0)$ holds,

\item[($\star$)]  
$\prop(f(x))$ holds for every $x \in X$ for which $\prop(x)$ holds,
\end{evlist}
then $\prop(x)$ holds for all $x \in X$. However, it turns out to be more convenient to work with the following 
property, which Lemma~\ref{lemma_intro_31} below shows is equivalent  to the principle of mathematical induction
holding: A counting system $(X,f,x_0)$ is said to be \definition{minimal} if the only $f$-invariant subset of $X$ 
containing $x_0$ is $X$ itself, where a subset $Y \subset X$ is \definition{$f$-invariant} if $f(Y) \subset Y$.

\begin{lemma}\label{lemma_intro_31}
The principle of mathematical induction holds for a counting system if and only if it is minimal.
\end{lemma}

\proof
Let $(X,f,x_0)$ be a counting system and let $\prop$ be a proposition about the elements of $X$ satisfying 
($\diamond$) and ($\star$). Then the subset $X' = \{ x \in X : \mbox{$\prop(x)$ holds} \}$ of $X$ is 
$f$-invariant and contains $x_0$. Therefore if $(X,f,x_0)$ is minimal then $X' = X$, and which means the 
principle of mathematical induction holds for $(X,f,x_0)$. Suppose conversely that $(X,f,x_0)$ is not minimal; 
then there exists an $f$-invariant subset $X'$ of $X$ containing $x_0$ with $X' \ne X$. For each $x \in X$ let 
$\prop(x)$ be the proposition that $x \in X'$; then ($\diamond$) and ($\star$) are satisfied by $\prop$, but 
$\prop(x)$ does not hold for $x \in X \setminus X'$ and so the principle of mathematical induction does not 
hold for $(X,f,x_0)$. \eop

The other two Peano axioms, when stated in terms of a counting system $(X,f,x_0)$, require that the mapping $f$ 
should be injective and that $f(x) \ne x_0$ for all $x \in X$ (i.e., that $x_0 \notin f(X)$), and a counting 
system satisfying these two conditions will be called \definition{standard}. The Peano axioms thus require that 
$(\Nat,\textsf{s},0)$ should be a minimal standard counting system, and such a counting system will be called 
a \definition{Dedekind system}.

It is the second property which implies completeness: Theorem~\ref{theorem_cfs_21} states that if $(X,f,x_0)$ 
is standard then the unique iterator for $(X,f,x_0)$ is complete. This confirms that the iterator for 
$(\Nat,\textsf{s},0)$ is complete.

Although the requirement that the counting system be minimal is not involved here, it will be needed in 
Theorem~\ref{theorem_cfs_31}, which states that if $(X,f,x_0)$ is a Dedekind system then for each counting system 
$(Y,g,y_0)$ there exists a unique mapping $h : X \to Y$ with $h(x_0) = y_0$ such that $h \circ f = g \circ h$. This 
result, which is known as  the \definition{recursion theorem} (at least 
when applied with $(X,f,x_0) = (\Nat,\textsf{s},0)$) is of fundamental importance, since it provides the justification 
for making recursive or inductive definitions.

Theorem~\ref{theorem_cfs_51} (a result of Lawvere \cite{lawvere}) will show that $(X,f,x_0)$ being a Dedekind system 
is necessary for the statement in the recursion theorem to hold, and so it is worth noting that the existence of a 
Dedekind system depends on the existence of a non-finite set. More precisely, a set $Y$ is said to be 
\definition{Dedekind-infinite} if there exists a mapping $g : Y \to Y$ which is injective but not surjective, and 
thus by Proposition~\ref{prop_intro_11} a Dedekind-infinite set cannot be finite.

\begin{proposition}\label{prop_intros_21}
The existence of a Dedekind system is equivalent to that of a Dedekind-infinite set.
\end{proposition}

\proof 
If $(X,f,x_0)$ is a standard counting system then the set $X$ is Dedekind-infinite, since $f : X \to X$ is an 
injective mapping which is not surjective. For the converse consider  a Dedekind-infinite set $X$, and so there 
exists an injective mapping $f : X \to X$ which is not surjective. Choose an element $x_0 \notin f(X)$, which 
gives us a counting system $(X,f,x_0)$. Now let $X_0$ be the least $f$-invariant subset of $X$ containing $x_0$
and $f_0$ be the restriction of $f$ to $X_0$, considered as a mapping  $X_0 \to X_0$. Then it is easy to see 
that the counting system $(X_0,f_0,x_0)$ is minimal. Moreover, $x_0 \notin f(X) \supset f(X_0) = f_0(X_0)$ and 
$f_0$, as the restriction of an injective mapping, is itself injective. Hence $(X_0,f_0,x_0)$ is also standard
i.e., $(X_0,f_0,x_0)$ is a Dedekind system.
\eop

As an application of the recursion theorem it is shown in Theorem~\ref{theorem_cfs_41} that if $(X,f,x_0)$ is a 
Dedekind system then there exists a unique mapping $[\,\cdot\,]$ from $X$ to the set of finite subsets of $X$
with $[x_0] = \varnothing$ and such that $[f(x)] = [x] \cup \{x\}$ for all $x \in X$. Moreover, $A \approx [\#(A)]$ 
holds for each finite set $A$, and in particular a set $A$ is finite if and only if $A \approx [x]$ for some 
$x \in X$. 

For the Dedekind system $(\Nat,\textsf{s},0)$ we have $[0] = \varnothing$ and $[n] = \{0,1,\ldots,n-1\}$ for each 
$n \in \Nat$ with $n \ne 0$. This shows that the definition of a finite set being employed here is equivalent to
the usual one.

Section~\ref{mcs} gives a more detailed account of minimal counting systems, and in particular of those which are 
not Dedekind systems. It is shown in Theorem~\ref{theorem_mcs_11} that a minimal counting system $(X,f,x_0)$ is 
standard (and thus a Dedekind system) if and only if the iterator is complete and that this is the case if and only 
if the set $X$ is not finite. The rest of the section is taken up with an analysis of 
minimal counting systems $(X,f,x_0)$ for which $X$ is finite. There are two cases.
In the first case $f$ is bijective (and so $x_0 \in f(X)$):

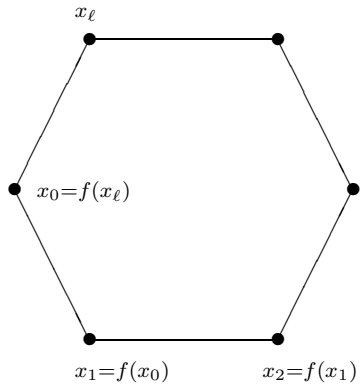
\setlength{\graphicthick}{0.1mm}
\setlength{\graphicmid}{0.1mm}
\setlength{\graphicthin}{0.1mm}

\begin{center}
\setlength{\unitlength}{1.0mm}
\begin{picture}(140,50)

\linethickness{\graphicthick}

\linethickness{\graphicthin}

\linethickness{\graphicmid}
\mythicklines

\put(49,24){$\bullet$}
\sput{53}{24}{x_0 = f(x_\ell)}

\sput{58}{0}{x_1 = f(x_0)}
\put(59,4){$\bullet$}

\sput{83}{0}{x_2 = f(x_1)}
\put(84,4){$\bullet$}

\put(59,44){$\bullet$}
\sput{58}{48}{x_\ell}

\put(84,44){$\bullet$}
\put(94,24){$\bullet$}

\put(50,25){\line(1,2){10}}
\put(50,25){\line(1,-2){10}}
\put(60,45){\line(1,0){25}}
\put(60,5){\line(1,0){25}}
\put(85,45){\line(1,-2){10}}
\put(85,5){\line(1,2){10}}

\end{picture}

\end{center}

\bigskip
In the second $f$ is not injective and $x_0 \notin f(X)$:

\bigskip

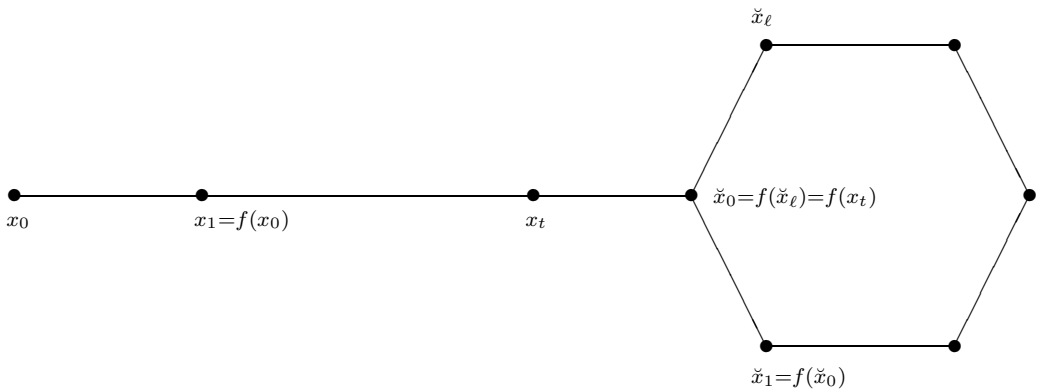
\begin{center}
\setlength{\unitlength}{1.0mm}
\begin{picture}(140,45)

\linethickness{\graphicthick}

\linethickness{\graphicthin}

\linethickness{\graphicmid}
\mythicklines

\put(5,25){\line(1,0){90}}

\put(95,25){\line(1,2){10}}
\put(95,25){\line(1,-2){10}}
\put(105,45){\line(1,0){25}}
\put(105,5){\line(1,0){25}}
\put(130,45){\line(1,-2){10}}
\put(130,5){\line(1,2){10}}

\sput{4}{21}{x_0}
\put(4,24){$\bullet$}

\sput{29}{21}{x_1 = f(x_0)}
\put(29,24){$\bullet$}

\sput{73}{21}{x_t}
\put(73,24){$\bullet$}

\sput{98}{24}{\breve{x}_0 = f(\breve{x}_\ell) = f(x_t)}
\put(94,24){$\bullet$}

\sput{103}{0}{\breve{x}_1 = f(\breve{x}_0)}
\put(104,4){$\bullet$}

\sput{103}{48}{\breve{x}_\ell}
\put(104,44){$\bullet$}

\put(129,4){$\bullet$}
\put(129,44){$\bullet$}
\put(139,24){$\bullet$}

\end{picture}

\end{center}

In Section~\ref{am} we show how an addition and a multiplication can be defined for any minimal counting system
$(X,f,x_0)$. These operations can be specified by the rules (a0), (a1), (m0) and (m1) below, which are usually 
employed when defining the operations on $\Nat$ via the Peano axioms.

Theorem~\ref{theorem_am_11} deals with the addition and states that there exists a unique binary operation $\oplus$ 
on $X$ such that
\[ 
\#(A) \oplus \#(B) = \#(A \cup B) 
\]
whenever $A$ and $B$ are disjoint finite sets, where $\#$ is the iterator for $(X,f,x_0)$. This operation $\oplus$ 
is both associative and commutative, $x \oplus x_0 = x$ for all $x \in X$ and for all $x_1,\,x_2 \in X$ there is an 
$x \in X$ such that either $x_1 = x_2 \oplus x$ or $x_2 = x_1 \oplus x$. Moreover, $\oplus$ is the unique binary 
operation $\oplus$ on $X$ 
such that
\begin{evlist}{32pt}{6pt}
\item[(a0)]  $x \oplus x_0 = x$ for all $x \in X$. 

\item[(a1)]  $x \oplus f(x') = f(x \oplus x')$ for all $x,\,x' \in X$.
\end{evlist}

Theorem~\ref{theorem_am_21} treats the multiplication and states that there exists a unique binary operation 
$\otimes$ on $X$ such that
\[ 
\#(A) \otimes \#(B) = \#(A\times B) 
\]
for all finite sets $A$ and $B$. This operation $\otimes$ is both associative and commutative, 
$x \otimes x_0 = x_0$ and $x \otimes f(x_0) = x$ for all $x \in X$ (and so $f(x_0)$ is a multiplicative unit) and 
the distributive law holds for $\oplus$ and $\otimes$: For all $x,\,x_1,\,x_2 \in X$
\[     
x \otimes (x_1 \oplus x_2) = (x \otimes x_1) \oplus (x \otimes x_2)\;. 
\]
Moreover, $\otimes$ is the unique binary operation on $X$ such that 
\begin{evlist}{32pt}{6pt}
\item[(m0)]  $x \otimes x_0 = x_0$ for all $x \in X$. 

\item[(m1)]  $x \otimes f(x') = x \oplus (x \otimes x')$ for all $x,\,x' \in X$.
\end{evlist}

\bigskip
Finally, Section~\ref{mam} presents alternative proofs for Theorems \ref{theorem_am_11} and \ref{theorem_am_21}.


\startsection{Finite sets}

\label{fs}

Let us first repeat the definition to be used here of a set being finite.

The set of all subsets of a set $A$ (the \definition{power set of $A$}) will be denoted as usual by $\mathcal{P}(A)$.
A subset $\mathcal{K}$ of $\mathcal{P}(A)$ is called an \definition{inductive system in $\mathcal{P}(A)$} (or just 
an \definition{inductive system} if $A$ can be determined from the context) if $\varnothing \in \mathcal{K}$ and 
$B \cup \{a\} \in \mathcal{K}$ for all $B \in \mathcal{K}$, $b \in A \setminus B$. In particular, $\mathcal{P}(A)$ 
itself is always an inductive system. A set $A$ will be called \definition{finite} if $\mathcal{P}(A)$ itself is the 
only inductive system. As was already mentioned, this definition is essentially that employed by Whitehead and 
Russell \cite{wr}. 

To keep the section self-contained let us also repeat Lemma~\ref{lemma_intro_21}.

\begin{lemma}\label{lemma_fs_21}
Let $A$ be finite and $b$ be some element not in $A$. Then $A' = A \cup \{b\}$ is also finite.
\end{lemma}

\proof 
Let $\mathcal{K}'$ be an inductive system in $\mathcal{P}(A')$; then 
$\mathcal{K} = \{ B \in \mathcal{K}' : B \subset A \}$ is clearly an inductive system in $\mathcal{P}(A)$ and 
hence $\mathcal{K} = \mathcal{P}(A)$. Thus $\mathcal{P}(A) \subset \mathcal{K}'$. Moreover, 
$B \cup \{a\} \in \mathcal{K}'$ for all $B \in \mathcal{P}(A)$, since $\mathcal{K}'$ is an inductive system in
$\mathcal{P}(A')$, and therefore $\mathcal{K}' = \mathcal{P}(A')$. It follows that $A'$ is finite.
\eop

Most proofs about finite sets take the following form: To show that every finite set $A$ has a particular 
property we consider the subset $\mathcal{K}$ of $\mathcal{P}(A)$ consisting of those subsets of $A$ which have 
the property. We then show that $\mathcal{K}$ is an inductive system and conclude that 
$\mathcal{K} = \mathcal{P}(A)$. In particular, it then follows that $A \in \mathcal{K}$, which shows that $A$ 
has the property. This means that what we usually need is not that $\mathcal{P}(A)$ is the only inductive system, 
but the apparently somewhat weaker statement that every inductive system in $\mathcal{P}(A)$ contains $A$.
However, as the next result shows, this statement is actually equivalent to $A$ being finite.

\begin{lemma}\label{lemma_fs_31}
The set $A$ is finite if and only if every inductive system in $\mathcal{P}(A)$ contains $A$.
\end{lemma}

\proof 
If $A$ is finite then $\mathcal{P}(A)$ is the only inductive system, and $\mathcal{P}(A)$ contains $A$.
Conversely, suppose that every inductive system in $\mathcal{P}(A)$ contains $A$, and consider the system of 
subsets $\mathcal{K} = \{ B \in \mathcal{P}(A) : \mbox{$B$ is finite} \}$. Then $\varnothing \in \mathcal{K}$ and 
if $B \in \mathcal{K}$ and $a \in A \setminus B$ then by Lemma~\ref{lemma_fs_21} $B \cup \{a\} \in \mathcal{K}$. 
Thus $\mathcal{K}$ is an inductive system and hence $A \in \mathcal{K}$, i.e., $A$ is finite.
\eop

We now establish the usual properties of finite sets. The proofs are mostly very straightforward and they all 
follow the same pattern (and so tend to become somewhat monotonous).

\begin{proposition}\label{prop_fs_21}
Every subset of a finite set is finite.
\end{proposition}

\proof 
Let $A$ be finite and put $\mathcal{K} = \{ B \in \mathcal{P}(A) : \mbox{$B$ is finite} \}$. Then 
$\varnothing \in \mathcal{K}$ and if $B \in \mathcal{K}$ and $a \in A \setminus B$ then by 
Lemma~\ref{lemma_fs_21} $B \cup \{a\} \in \mathcal{K}$. Thus $\mathcal{K}$ is an inductive system and hence 
$\mathcal{K} = \mathcal{P}(A)$, i.e., every subset of $A$ is finite. 
\eop

\begin{proposition}\label{prop_fs_31}
If $A$ and $B$ are finite sets then so is $A \cup B$.
\end{proposition}

\proof
Consider the system of subsets $\mathcal{K} = \{ C \in \mathcal{P}(A) : \mbox{$C \cup B$ is finite} \}$. Then 
$\varnothing \in \mathcal{K}$, since by assumption $\varnothing \cup B = B$ is finite and if $C \in \mathcal{K}$
(i.e., $C \cup B$ is finite) and $a \in A \setminus C$ then by Lemma~\ref{lemma_fs_21} 
$(C \cup \{a\}) \cup B = (C \cup B) \cup \{a\} \in \mathcal{K}$. Thus $\mathcal{K}$ is an inductive system in 
$\mathcal{P}(A)$ and therefore $\mathcal{K} = \mathcal{P}(A)$. In particular, $A \in \mathcal{K}$, i.e.,
$A \cup B$ is finite. 
\eop

\begin{proposition}\label{prop_fs_41}
Let $A$ and $B$ be sets with $A$ finite.

(1)\enskip 
If there exists an injective mapping $f : B \to A$ then $B$ is also finite.

(2)\enskip 
If there exists a surjective mapping $f : A \to B$ then $B$ is again finite.
\end{proposition}

\proof
(1)\enskip
Let $\mathcal{K}$ consist of those subsets $C$ of $A$ such that if $D$ is any set for which there exists an 
injective mapping $f : D \to C$ then $D$ is finite. Then $\varnothing \in \mathcal{K}$, since there can only 
exist a mapping $f : D \to \varnothing$ if $D = \varnothing$ and the empty set $\varnothing$ is finite. Let 
$C \in \mathcal{K}$ and $a \in A \setminus C$. Consider a set $D$ for which there exists an injective mapping 
$f : D \to C \cup \{a\}$. There are two cases:

($\alpha$)\enskip $f(d) \in A$ for all $d \in D$. Here we can consider $f$ as a mapping from $D$ to $C$ and as 
such it is still injective. Thus $D$ is finite since  $C \in \mathcal{K}$.

($\beta$)\enskip 
There exists an element $b \in D$ with $f(b) = a$. Put $D' = D \setminus \{b\}$. Now since $f$ is injective it 
follows that $f(d) \ne a$ for all $d  \in D'$, and thus we can define a mapping $g : D' \to C$ by letting 
$g(d) = f(d)$ for all $d \in D'$. Then $g : D' \to C$ is also injective (since $f : D \to C \cup \{a\}$ is) and 
therefore $D'$ is finite since $C \in \mathcal{K}$ holds. Hence by Lemma~\ref{lemma_fs_21} $D = D' \cup \{b\}$ 
is finite.

This shows that $C \cup \{a\} \in \mathcal{K}$ and therefore $\mathcal{K}$ is an inductive system in 
$\mathcal{P}(A)$. Thus $\mathcal{K} = \mathcal{P}(A)$ and in particular $A \in \mathcal{K}$, which means that if 
there exists an injective mapping $f : B \to A$ then $B$ is also finite.

(2)\enskip
Let $\mathcal{K}$ consist of those subsets $C$ of $A$ such that if $D$ is any set for which there exists a 
surjective mapping $f : C \to D$ then $D$ is finite. Then $\varnothing \in \mathcal{K}$, since there can only 
exist a surjective mapping $f : \varnothing \to D$ if $D = \varnothing$ and the empty set $\varnothing$ is finite.
Let $C \in \mathcal{K}$ and $a \in A \setminus C$. Consider a set $D$ for which there exists a surjective mapping 
$f : C \cup \{a\} \to D$. There are again two cases:

($\alpha$)\enskip The restriction $f_{|C} : C \to D$ of $f$ to $C$ is still surjective. Then $D$ is finite since 
$C \in \mathcal{K}$.

($\beta$)\enskip The restriction $f_{|C}$ is not surjective. Put $b = f(a)$ and $D' = D \setminus \{b\}$. Then 
$f(c) \ne b$ for all $c \in C$ (since $f_{|C}$ is not surjective) and therefore we can define a mapping 
$g : C \to D'$ by letting $g(c) = f(c)$ for all $c \in C$. But $f : C \cup \{a\} \to D$ is surjective and hence 
$g : C \to D'$ is also surjective. Thus $D'$ is finite since $C \in \mathcal{K}$ holds, and so by 
Lemma~\ref{lemma_fs_21} $D = D' \cup \{b\}$ is finite.

This shows that $C \cup \{a\} \in \mathcal{K}$ and therefore $\mathcal{K}$ is an inductive system in 
$\mathcal{P}(A)$. Thus $\mathcal{K} = \mathcal{P}(A)$ and in particular $A \in \mathcal{K}$, which means that
if there exists a surjective  mapping $f : A \to B$ then $B$ is also finite. 
\eop

\begin{proposition}\label{prop_fs_51}
If $A$ is a finite set then so is the power set $\mathcal{P}(A)$.
\end{proposition}

\proof 
Consider the system of subsets $\mathcal{K} = \{ B \in \mathcal{P}(A) : \mbox{$\mathcal{P}(B)$ is finite}\}$.
Then by Lemma~\ref{lemma_fs_21} $\varnothing \in \mathcal{K}$, since 
$\mathcal{P}(\varnothing) = \{\varnothing\} = \varnothing \cup \{\varnothing\}$. Let $B \in \mathcal{K}$ and 
$a \in A \setminus B$. Then $\mathcal{P}(B \cup \{a\}) = \mathcal{P}(B) \cup \mathcal{P}_a(B)$, where 
$\mathcal{P}_a(B) = \{ C \cup \{a\} : C \in \mathcal{P}(B) \}$ and the mapping $C \mapsto C \cup \{a\}$ from
$\mathcal{P}(B)$ to $\mathcal{P}_a(B)$ is surjective. It follows from Proposition~\ref{prop_fs_41}~(2) that 
$\mathcal{P}_a(B)$ is finite and thus by Proposition~\ref{prop_fs_31} $\mathcal{P}(B \cup \{a\})$ is finite, 
i.e., $B \cup \{a\} \in \mathcal{K}$. This shows that $\mathcal{K}$ is an inductive system in $\mathcal{P}(A)$. 
Hence $\mathcal{K} = \mathcal{P}(A)$ and in particular $A \in \mathcal{K}$, i.e., the power set $\mathcal{P}(A)$ 
is finite. 
\eop

\begin{proposition}\label{prop_fs_61}
If $A$ and $B$ are finite sets then so their product $A \times B$.
\end{proposition}

\proof
Put $\mathcal{K} = \{ C \in \mathcal{P}(A) : \mbox{$C \times B$ is finite}\}$. Then $\varnothing \in \mathcal{K}$, 
since $\varnothing \times B = \varnothing$. Let $C \in \mathcal{K}$ and $a \in A \setminus C$. Then 
$(C \cup \{a\}) \times B = (C \times B) \cup (\{a\} \times B)$ and by Proposition~\ref{prop_fs_41}~(2) 
$\{a\} \times B$ is finite since the mapping $f : B \to \{a\} \times B$ with $f(b) = (a,b)$ for all $b \in B$ is 
surjective. Thus by Proposition~\ref{prop_fs_21} $(C \cup \{a\}) \times B$ is finite, i.e., 
$C \cup \{a\} \in \mathcal{K}$. This shows that $\mathcal{K}$ is an inductive system in $\mathcal{P}(A)$. Hence 
$\mathcal{K} = \mathcal{P}(A)$ and in particular $A \in \mathcal{K}$, i.e., $A \times B$ is finite. 
\eop

\begin{proposition}\label{prop_fs_71}
If $A$ and $B$ are finite sets then so is $B^A$, the set of all mappings from $A$ to $B$.
\end{proposition}

\proof
Since $B^A$ is a subset of $\mathcal{P}(A \times B)$ and by Propositions~\ref{prop_fs_51} and \ref{prop_fs_61} 
the set $\mathcal{P}(A \times B)$ is finite it follows from Proposition~\ref{prop_fs_21} that $B^A$ is finite.
\eop

\begin{proposition}\label{prop_fs_81}
Let $A$ be a finite set and $f : A \to A$ be a mapping. Then  $f$ is injective if and only if it is 
surjective (and thus if and only if it is bijective).
\end{proposition}

\proof 
We first show that an injective mapping is bijective. Let $\mathcal{K}$ consist of those subsets $B$ of $A$ having 
the property that every injective mapping $f : B \to B$ is bijective. Then $\varnothing \in \mathcal{K}$, since the 
only mapping $f : \varnothing \to \varnothing$ is bijective. Let $B \in \mathcal{K}$ and $a \in A \setminus B$; 
consider an injective mapping $f : B \cup \{a\} \to B \cup \{a\}$. There are two cases:

($\alpha$)\enskip 
$f(B) \subset B$. Then the restriction $f_{|B} : B \to B$ of $f$ to $B$ is injective
and hence bijective, since $B \in \mathcal{K}$. If $f(a) \in B$ then $f(b) = f_{|B}(b) = f(a)$ for 
some $b \in B$, since $f_{|B}$ is surjective, which contradicts the fact that $f$ is injective. Thus $f(a) = a$, and 
it follows that $f$ is bijective.

($\beta$)\enskip 
$f(B) \not\subset B$. In this case there exists $b \in B$ with $f(b) = a$ and, since $f$ is injective, we must have 
$f(c) \in B$ for all $c \in B \setminus \{b\}$ and $f(a) \in B$. This means there is an injective mapping 
$g : B \to B$ defined by letting
\[ g(c) = \left\{ \begin{array}{cl}
                  f(c) &\ \mbox{if}\ c \in B \setminus \{b\}\;,\\
                  f(a)   &\  \mbox{if}\ c = b\;\\
\end{array} \right. \]
and then $g$ is bijective, since $B \in \mathcal{K}$. Therefore $f$ is again bijective.
 
This shows that $B \cup \{a\} \in \mathcal{K}$ and thus that $\mathcal{K}$ is an inductive system. Hence 
$\mathcal{K} = \mathcal{P}(A)$ and in particular $A \in \mathcal{K}$, i.e., every injective mapping 
$f : A \to A$ is bijective.

We now show  a surjective mapping is bijective, and here let $\mathcal{K}$ consist of those subsets $B$ of $A$ 
having the property that every surjective mapping $f : B \to B$ is bijective. Then $\varnothing \in \mathcal{K}$, 
again since the only mapping $f : \varnothing \to \varnothing$ is bijective. Let $B \in \mathcal{K}$ and 
$a \in A \setminus B$; consider a surjective mapping $f : B \cup \{a\} \to B \cup \{a\}$. Let 
$D = \{ b \in B : f(b) = a \}$; there are three cases:

($\alpha$)\enskip 
$D = \varnothing$. Then $f(a) = a$, since $f$ is surjective, thus the restriction $f_{|B} : B \to B$ of $f$ to 
$B$ is surjective and hence bijective (since $B \in \mathcal{K}$ holds), and this means $f$ is bijective. 

($\beta$)\enskip 
$D \ne \varnothing$ and $f(a) \in B$. Here we can define a surjective mapping $g : B \to B$ by letting
\[ 
   g(c) = \left\{ \begin{array}{cl}
                  f(c) &\ \mbox{if}\ c \in B \setminus D\;,\\
                  f(a)   &\  \mbox{if}\ c \in D\;.\\
\end{array} \right. 
\] 
Thus $g$ is bijective (since $B \in \mathcal{K}$), which implies that $D = \{b\}$ for some $b \in C$ and in 
particular $f$ is also injective. 

($\gamma$)\enskip 
$D \ne \varnothing$ and $f(a) = a$. This is not possible since then $f(B \setminus D) = B$ and so, choosing any 
$b \in D$, the mapping $h : B \to B$ with
\[ 
   h(c) = \left\{ \begin{array}{cl}
                  f(c) &\ \mbox{if}\ c \in B \setminus D\;,\\
                  b   &\  \mbox{if}\ c \in D\\
\end{array} \right. 
\] 
would be surjective but not injective (since there also exists $c \in B \setminus D$ with $f(c) = b$).

This shows that $B \cup \{a\} \in \mathcal{K}$ and thus that $\mathcal{K}$ is an inductive system. Hence 
$\mathcal{K} = \mathcal{P}(A)$ and in particular $A \in \mathcal{K}$, i.e., every surjective mapping 
$f : A \to A$ is bijective.
\eop

If $A$ and $B$ are sets then we write $A \approx B$ if there exists a bijective mapping $f : A \to B$.
Proposition~\ref{prop_fs_41} implies that if $A \approx B$ then $A$ is finite if and only if $B$ is. It is clear 
that $\approx$ defines an equivalence relation on the class of all finite sets.

\begin{proposition}\label{prop_fs_91}
Let $A$ and $B$ be finite sets. Then either there exists an injective mapping $f : A \to B$ or an injective 
mapping $g : B \to A$. Moreover, if there exists both an injective mapping $f : A \to B$ and an injective mapping
$g : B \to A$ then $A \approx B$.
\end{proposition}

\proof
Let $\mathcal{K}$ consist of those subsets $C$ of $A$ for which there either there exists an injective mapping 
$f : C \to B$ or an injective mapping $g : B \to C$. Then $\varnothing \in \mathcal{K}$, since the only mapping 
$f : \varnothing \to B$ is injective. Let $C \in \mathcal{K}$ and let $a \in A \setminus C$. There are two cases:

($\alpha$)\enskip 
There exists an injective mapping $g : B \to C$.  Then $g$ is still injective when considered as a mapping from 
$B$ to $C \cup \{a\}$.

($\beta$)\enskip 
There exists an injective mapping $f : C \to B$. If $f$ is not surjective then it can be extended to an injective 
mapping $f' : C \cup \{a\} \to B$ (with $f'(a)$ chosen to be any element in $B \setminus f(C)$). On the other 
hand, if $f$ is surjective (and hence a bijection) then the inverse mapping $f^{-1} : B \to C$ is injective and 
so is still injective when considered as a mapping from $B$ to $C \cup \{a\}$. 

This shows that $B \cup \{a\} \in \mathcal{K}$ and thus that $\mathcal{K}$ is an inductive system. Hence 
$\mathcal{K} = \mathcal{P}(A)$ and in particular $A \in \mathcal{K}$, i.e., there either exists an injective 
mapping $f : A \to B$ or an injective mapping $g : B \to A$.

Suppose there exists both an injective mapping $f : A \to B$ and an injective mapping $g : B \to A$. Then 
$f \circ g : B \to B$ is an injective mapping, which by Proposition~\ref{prop_fs_81} is bijective. In particular 
$f$ is surjective and therefore bijective, i.e., $A \approx B$. 
\eop

\begin{lemma}\label{lemma_fs_41}
Let $A$ and $B$ be finite sets. Then there exists either a subset $B'$ of $A$ with $B' \approx B$ or a subset $A'$ 
of $B$ with $A' \approx A$.
\end{lemma}

\proof 
By Proposition~\ref{prop_fs_91} there either exists an injective mapping $f : A \to B$ or an injective mapping 
$g : B \to A$. Suppose that the former is the case and put $A' = f(A)$; then $A' \subset B$ with $A' \approx A$ 
(since $f$ as a mapping from $A$ to $A'$ is a bijection). If there exists an injective mapping $g : B \to A$ then 
in the same way there exists a subset $B'$ of $A$ with $B' \approx B$.
\eop

\begin{proposition}\label{prop_fs_101}
If $B$ is a subset of a finite set $A$ with $B \approx A$ then $B = A$.
\end{proposition}

\proof 
There exists a bijective mapping $f : A \to B$ and then the restriction $f_{|B} : B \to B$ of $f$ to $B$ is 
injective; thus by Proposition~\ref{prop_fs_81} $f_{|B}$ is bijective. But this is only possible if $B = A$, 
since if $a \in A \setminus B$ then $f(a) \notin f_{|B}(B)$.
\eop

If $A$ is a set and $\mathcal{S}$ a non-empty subset of $\mathcal{P}(A)$ then $B \in \mathcal{S}$ is said to be 
\definition{minimal} if $B' \notin \mathcal{S}$ for each proper subset $B'$ of $B$. The statement in the 
following result is Tarski's definition \cite{tarski} of a set being finite.

\begin{proposition}\label{prop_fs_111}
Let $A$ be a set; then each non-empty subset of $\mathcal{P}(A)$ contains a minimal element if and only if $A$ 
is finite.
\end{proposition}

\proof 
Let $A$ be a finite set and let $\mathcal{K}$ consist of those subsets $B$ of $A$ such that each non-empty subset
of $\mathcal{P}(B)$ contains a minimal element. Then $\varnothing \in \mathcal{K}$, since the only non-empty 
subset of $\mathcal{P}(\varnothing)$ is $\{\varnothing\}$ and then $\varnothing$ is the required minimal element.
Let $B \in \mathcal{K}$ and $a \in A \setminus B$, and let $\mathcal{S}$ be a non-empty subset of 
$\mathcal{P}(B \cup \{a\})$. Put $\mathcal{S}' = \mathcal{S} \cap \mathcal{P}(B)$; there are two cases:

($\alpha$)\enskip 
$\mathcal{S}' \ne \varnothing$. Here $\mathcal{S}'$ is a non-empty subset of $\mathcal{P}(B)$ and thus contains 
a minimal element $C$ which is then a minimal element of $\mathcal{S}$, since each set in 
$\mathcal{S} \setminus \mathcal{S}'$ contains $a$ and so is not a proper subset of $C$.
 
($\beta$)\enskip 
$\mathcal{S}' = \varnothing$ (and so each set in $\mathcal{S}$ contains $a$). Put 
$\mathcal{S}'' = \{ C \subset B : C \cup \{a\} \in \mathcal{S} \}$; then $\mathcal{S}''$ is a non-empty subset 
of $\mathcal{P}(B)$ and thus contains a minimal element $C$. It follows that $C' = C \cup \{a\}$ is a minimal 
element of $\mathcal{S}$: A proper subset of $C'$ has either the form $D$ with $D \subset C$, in which case 
$D \notin \mathcal{S}$ (since each set in $\mathcal{S}$ contains $a$), or has the form $D \cup \{a\}$ with $D$ a 
proper subset of $C$ and here $D \cup \{a\} \notin \mathcal{S}$, since $D \notin \mathcal{S}''$.

This shows that $B \cup \{a\} \in \mathcal{K}$ and thus that $\mathcal{K}$ is an inductive system. Hence 
$\mathcal{K} = \mathcal{P}(A)$ and in particular $A \in \mathcal{K}$, i.e., non-empty subset of $\mathcal{P}(X)$ 
contains a minimal element.

Conversely, suppose $A$ is not finite and let 
$\mathcal{S} = \{ B \in \mathcal{P}(A) : \mbox{$B$ is not finite} \}$; then $\mathcal{S}$ is non-empty since it 
contains $A$. However $\mathcal{S}$ cannot contain a minimal element: If $B$ is a minimal element of 
$\mathcal{S}$ then $B \ne \varnothing$, since $\varnothing$ is finite. Choose $b \in B$; then $B \setminus \{b\}$
is a proper subset of $B$ and thus $B \setminus \{b\} \notin \mathcal{S}$, i.e., $B \setminus \{b\}$ is finite. 
But then by Lemma~\ref{lemma_fs_21} $B = (B \setminus \{b\}) \cup \{b\}$ would be finite. 
\eop

For what we consider in later sections it is useful to employ another technique for establishing statements about 
finite sets. This involves the following \definition{induction principle for finite sets} which first appeared in 
a 1909 paper of Zermelo \cite{zermelo}:

\newpage

\begin{theorem}\label{theorem_fs_11}
Let $\prop$ be a proposition about finite sets (meaning that for each finite set $A$ we have some proposition 
$\prop(A)$). Suppose that
\begin{evlist}{16pt}{5pt}
\item[$\mathrm{(\diamond)}$]  $\prop(\varnothing)$ holds.

\item[$\mathrm{(\star)}$]  
If $A$ is a finite set for which $\prop(A)$ holds then $\prop(A \cup \{a\})$ holds for each 
\phantom{xix}element $a \notin A$.
\end{evlist}

Then $\prop$ is a property of finite sets, i.e., $\prop(A)$ holds for every finite set $A$.
\end{theorem}

\proof
Let $A$ be a finite set and put $\mathcal{K} =  \{ B \in \mathcal{P}(A) : \mbox{$\prop(B)$ holds} \}$. In 
particular ($\diamond$) implies $\varnothing \in \mathcal{K}$. Consider $B \in \mathcal{K}$ (and so $\prop(B)$ 
holds) and let $a \in A \setminus B$; then $\prop(B \cup \{a\})$ holds by ($\star$), and hence 
$B \cup \{a\} \in \mathcal{K}$. This shows that $\mathcal{K}$ is an inductive system in $\mathcal{P}(A)$ and so 
$\mathcal{K} = \mathcal{P}(A)$. In particular $A \in \mathcal{K}$, i.e., $\prop(A)$ holds, and since $A$ is 
arbitrary $\prop(A)$ holds for every finite set $A$. 
\eop

All the properties of finite sets presented above could have been deduced from Theorem~\ref{theorem_fs_11} 
(together with Lemma~\ref{lemma_fs_21}, which states that $A \cup \{a\}$ is finite for each finite set $A$ and 
each element $a \notin A$, and the fact that the empty set $\varnothing$ is finite). For example, consider 
Proposition~\ref{prop_fs_31} which states that if $A$ and $B$ are finite sets then so is $A \cup B$. To establish 
this using the induction principle for finite sets regard $B$ as being fixed and for each finite 
set $A$ let $\prop(A)$ be the proposition that $A \cup B$ is finite. Then:

($\diamond$)\enskip $\prop(\varnothing)$ holds because $\varnothing \cup B = B$.

($\star$)\enskip Let $A$ be a finite set for which $\prop(A)$ holds and let $a \notin A$. Now $A \cup B$ is 
finite, since $\prop(A)$ holds, and thus by Lemma~\ref{lemma_fs_21} $(A \cup B) \cup \{a\}$ is also finite 
(since this holds immediately if $a \in A \cup B$). But $(A \cup \{a\}) \cup B = (A \cup B) \cup \{a\}$; i.e.,
$(A \cup \{a\}) \cup B$ is finite. This shows that $\prop(A \cup \{a\})$ holds.

Therefore by Theorem~\ref{theorem_fs_11} $\prop(A)$ holds for each finite set $A$, and thus for all finite sets 
$A,\,B$ the set $A \cup B$ is finite.

The reader is left to check that the proofs of all the other results about finite sets can be obtained in
this manner.


\startsection{Counting systems}

\label{cfs}

Having introduced a theory of finite sets we now turn to the question of what it means to count the elements in such 
a set. Let us first recall some of the definitions made in Section~\ref{intro}.

A triple $(X,f,x_0)$ consisting of a set $X$, a mapping $f : X \to X$ and an element $x_0 \in X$ will be called a 
\definition{counting system}. A counting system $(X,f,x_0)$ is said to be \definition{minimal} if the only 
$f$-invariant subset of $X$ containing $x_0$ is $X$ itself, where a subset $Y \subset X$ is 
\definition{$f$-invariant} if $f(Y) \subset Y$. Moreover, it will be called \definition{standard} if the mapping 
$f$ is injective and $f(x) \ne x_0$ for all $x \in X$ (i.e., $x_0 \notin f(X)$). The Peano axioms thus require that 
$(\Nat,\textsf{s},0)$ should be a minimal standard counting system, and such a counting system will be called a 
\definition{Dedekind system}.

Consider a counting system $(X,f,x_0)$ and for each finite set $A$ let $\#(A)$ be an element of $X$. Then the 
assignment $A \mapsto \#(A)$ will be called an \definition{iterator for $(X,f,x_0)$} if $\#(\varnothing) = x_0$ and 
$\#(A \cup \{a\}) = f(\#(A))$ whenever $A$ is a finite set and $a \notin A$. (This means, for example, that 
$\#(\{a\}) = f(x_0)$ for each element $a$, $\#(\{a,b\}) = f(f(x_0))$ for distinct elements $a$ and $b$, and so on.)
Theorem~\ref{theorem_cfs_11} states that there exists a unique iterator $\#$ for each counting system $(X,f,x_0)$ 
and by Proposition~\ref{prop_cfs_21} $\#(B) = \#(A)$ whenever $A$ and $B$. 

If the iterator $\#$ for $(X,f,x_0)$ has the additional property that $B \approx A$ whenever $A$ and $B$ are finite 
sets with $\#(B) = \#(A)$ then we say that it is \definition{complete}. In general $\#$ will not have this property, 
as can be seen by looking at the trivial example in which $X$ consist of the single element $x_0$, and so $f$ is the 
mapping with $f(x_0) = x_0$. Then $\#(A) = x_0$ for each finite set $A$ and here $\#(B) = \#(A)$ for all finite sets 
$A$ and $B$. However, Theorem~\ref{theorem_cfs_21} states that if $(X,f,x_0)$ is standard then the unique iterator 
is complete. 

In this section most of the statements about finite sets will be established with the help of the 
\definition{induction principle for finite sets} (Theorem~\ref{theorem_fs_11}). The counting system $(X,f,x_0)$ is 
considered to be fixed in what follows.

\begin{theorem}\label{theorem_cfs_11}
There exists a unique iterator $\#$ for $(X,f,x_0)$.
\end{theorem}

\proof 
We first consider a version of the theorem restricted to the subsets of a finite set. Let $A$ be a finite set; 
then a mapping $\#_A : \mathcal{P}(A) \to X$ will be called an \definition{$A$-iterator} if 
$\#_A(\varnothing) = x_0$ and $\#_A(B \cup \{a\}) = f(\#_A(B))$ for each $B \subset A$ and each 
$a \in A \setminus B$.

\begin{lemma}\label{lemma_cfs_21}
For each finite set $A$ there exists a unique $A$-iterator.
\end{lemma}

\proof
For each finite set $A$ let $\prop(A)$ be the proposition that there exists a unique $A$-iterator.

($\diamond$)\enskip 
$\prop(\varnothing)$ holds, since the mapping $\#_\varnothing : \mathcal{P}(\varnothing) \to X$ with 
$\#_\varnothing(\varnothing) = x_0$ is clearly the unique $\varnothing$-iterator.

($\star$)\enskip Let $A$ be a finite set for which $\prop(A)$ holds, and let $a \notin A$; put 
$A' = A \cup \{a\}$. By assumption there exists a unique $A$-iterator $\#_A$. Now $\mathcal{P}(A')$ is 
the disjoint union of the sets $\mathcal{P}(A)$ and $\{ B \cup \{a\} : B \subset A\}$ and so we can define 
a mapping $\#_{A'} : \mathcal{P}(A') \to X$ by $\#_{A'}(B) = \#_A(B)$ and $\#_{A'}(B\cup \{a\}) = f(\#_A(B))$ 
for each $B \subset A$. Then $\#_{A'}(\varnothing) = \#_A(\varnothing) = x_0$, and for all $B \subset A$ and 
all $b \in A \setminus B$ 
\begin{eqnarray*}
\#_{A'}(B \cup \{b\}) &=& \#_A(B \cup \{b\}) = f(\#_A(B)) = f(\#_{A'}(B)) \;,\\
 \#_{A'}(B \cup \{a\}) &=& f(\#_A(B)) = f(\#_{A'}(B)) \;,\\
 \#_{A'}(B \cup \{a\} \cup \{b\}) &=& f(\#_A(B \cup \{b\})) = f(f(\#_A(B))) =  f(\#_{A'}(B \cup \{a\})) \;,
\end{eqnarray*}
i.e., $\#_{A'}(B' \cup \{b\}) = f(\#_{A'}(B')$ for all $B' \subset A'$ and $b \in A' \setminus B'$ and this 
means that $\#_{A'}$ is an $A$-iterator in $(X,f,x_0)$. Let $\#'_{A'}$ be an arbitrary $A'$-iterator. Then in 
particular $\#'_{A'}(B \cup \{b\}) = f(\#'_{A'}(B)$ for all $B \subset A$ and all $b \in A \setminus B$, and 
from  the uniqueness of the $A$-iterator $\#_A$ it follows that $\#'_{A'}(B) = \#_A(B)$ and thus also
$\#'_{A'}(B \cup \{a\}) = f(\#'_{A'}(B)) = f(\#_{A'}(B)) = \#_{A'}(B \cup \{a\})$ for all $B \subset A$, i.e., 
$\#'_{A'} = \#_{A'}$. This shows that $\prop(A \cup \{a\})$ holds.

Therefore by the induction principle for finite sets $\prop(A)$ holds for every finite set $A$, i.e., for each 
finite set $A$ there exists a unique $\#_A$-iterator. \eop

\begin{lemma}\label{lemma_cfs_31}
Let $A$ and $B$ be finite sets with $B \subset A$; then the unique $B$-iterator $\#_B$ is the restriction of 
the unique $A$-iterator $\#_A$ to $\mathcal{P}(B)$, i.e., $\#_B(C) = \#_A(C)$ for all $C \in \mathcal{P}(B)$.
\end{lemma}

\proof 
This follows immediately from the uniqueness of $\#_B$.
\eop

Now for each finite set $A$ let $\#_A$ be the unique $A$-iterator and put $\#(A) = \#_A(A)$. In particular 
$\#(\varnothing) = \#_\varnothing(\varnothing) = x_0$. For each finite set $A$ and each element $a \notin A$ 
it follows from Lemma~\ref{lemma_cfs_31} that
\[ 
\#(A\cup \{a\}) = \#_{A\cup \{a\}}(A\cup \{a\}) = f(\#_{A\cup \{a\}}(A) = f(\#_A(A)) = f(\#(A))  
\]
and hence $\#$ is an iterator for $(X,f,x_0)$.

Finally, consider an arbitrary iterator $\#'$ for $(X,f,x_0)$ and for each finite set $A$ let $\prop(A)$ be 
the proposition that $\#'(A) = \#(A)$.

($\diamond$)\enskip 
$\prop(\varnothing)$ holds since $\#'(\varnothing) = x_0 = \#(\varnothing)$.

($\star$)\enskip Let $A$ be a finite set for which $\prop(A)$ holds (and so $\#'(A) = \#(A)$) and let
$a \notin A$. Then $\#'(A \cup \{a\}) = f(\#'(A)) = f(\#(A)) = \#(A \cup \{a\})$ and so $\prop(A \cup \{a\})$ 
holds.

Hence by the induction principle for finite sets $\prop(A)$ holds for every finite set $A$, which means  
$\#'(A) = \#(A)$ for each finite set $A$, i.e., the iterator $A \to \#(A)$ for $(X,f,x_0)$ is unique. This 
completes the proof of Theorem~\ref{theorem_cfs_11}.
\eop

In what follows let $\#$ be the unique iterator for $(X,f,x_0)$.

\begin{proposition}\label{prop_cfs_21}
If $A$ and $B$ are finite sets with $B \approx A$ then $\#(B) = \#(A)$.
\end{proposition}

\proof 
For each finite set $A$ let $\prop(A)$ be the proposition that $\#(B) = \#(A)$ whenever $B$ is a finite set 
with $B \approx A$.

($\diamond$)\enskip 
$\prop(\varnothing)$ holds since $B \approx \varnothing$ if and only if $B = \varnothing$.

($\star$)\enskip Let $A$ be a finite set for which $\prop(A)$ holds and let $a \notin A$. Consider a finite set 
$B$ with $B \approx A \cup \{a\}$, put $b = h(a)$ and let $B' = B \setminus \{b\}$; then $B' \approx A$, thus
$\#(B') = \#(A)$, since $\prop(A)$ holds, and it follows that
\[
\#(B) = \#(B' \cup \{b\}) = f(\#(B')) = f(\#(A)) = \#(A \cup \{a\})\;.
\]
This shows that $\prop(A \cup \{a\})$ holds.

Therefore by the induction principle for finite sets $\prop(A)$ holds for every finite set $A$, and hence 
$\#(B) = \#(A)$ whenever $B$ is a finite set with $B \approx A$.
\eop

\begin{theorem}\label{theorem_cfs_21}
If $(X,f,x_0)$ is standard then the iterator $\#$ is complete.
\end{theorem}

\proof 
For each finite set $A$ let $\prop(A)$ be the proposition that $B \approx A$ whenever $B$ is a finite set with 
$\#(B) = \#(A)$.

($\diamond$)\enskip 
Let $B$ be a finite set with $B \ne \varnothing$, let $b \in B$ and put $B' = B \setminus \{b\}$. Then 
$\#(B) = \#(B' \cup \{b\}) = f(\#(B'))$, and so $\#(B) \ne x_0$, since $x_0 \notin f(X)$. Thus 
$\#(B) \ne \#(\varnothing)$, which shows that $\prop(\varnothing)$ holds, since $\varnothing \approx B$ if and 
only if $B = \varnothing$.

($\star$)\enskip Let $A$ be a finite set for which $\prop(A)$ holds and let $a \notin A$. Consider a finite set 
$B$ with $\#(B) = \#(A \cup \{a\})$; then $\#(B) = f(\#(A)) \in f(X)$, hence $\#(B) \ne x_0$ and so 
$B \ne \varnothing$. Let $b \in B$ and put $B' = B \setminus \{b\}$; then
\[ 
f(\#(B')) = \#(B' \cup \{b\}) = \#(B) = f(\#(A))\;,
\] 
and thus $\#(B') = \#(A)$, since $f$ is injective, and it follows that $B' \approx A$, since $\prop(A)$ holds.
But $B = B' \cup \{b\}$ with $b \notin B'$, $a \notin A$ and $B' \approx A$, and therefore 
$B = B' \cup \{b\} \approx A \cup \{a\}$. This shows that $\prop(A \cup \{a\})$ holds.

Thus by the induction principle for finite sets $\prop(A)$ holds for every finite set $A$, which means that if 
$A$ and $B$ are finite sets with $\#(B) = \#(A)$ then $B \approx A$.
\eop

The next lemma indicates the advantage of having a minimal counting system for obtaining information using the 
iterator $\#$. Note that an arbitrary intersection of $f$-invariant subsets of $X$ is again $f$-invariant and 
$X$ is itself an $f$-invariant subset containing $x_0$. There is thus a least $f$-invariant subset of $X$ 
containing $x_0$ (namely the intersection of all such subsets). If $(X,f,x_0)$ is minimal then the least 
$f$-invariant subset containing $x_0$  is of course $X$ itself.

\begin{lemma}\label{lemma_cfs_41}
Let $X_0$ be the least $f$-invariant subset of $X$ containing $x_0$. Then
\[ 
X_0 = \{ x \in X : \mbox{$x = \#(A)$ for some finite set $A$} \}\;.
\]
In particular, if $(X,f,x_0)$ is minimal then for each $x \in X$ there exists a finite set $A$ such that 
$x = \#(A)$.
\end{lemma}

\proof 
Put $X'_0 = \{ x \in X : \mbox{$x = \#(A)$ for some finite set $A$} \}$. Now let $X'$ be an $f$-invariant subset 
of $X$ containing $x_0$ and for each finite set $A$ let $\prop(A)$ be the proposition that $\#(A) \in X'$.

($\diamond$)\enskip 
$\prop(\varnothing)$ holds since $\#(\varnothing) = x_0 \in X'$.

($\star$)\enskip Let $A$ be a finite set for which $\prop(A)$ holds (and so $\#(A) \in X'$) and let $a \notin A$. 
Then $\#(A \cup \{a\}) = f(\#(A))\in X'$, since $X'$ is $f$-invariant, and hence $\prop(A \cup \{a\})$ holds.

Therefore by the induction principle for finite sets $\prop(A)$ holds for every finite set $A$, which means that 
$\#(A) \in X'$ for each finite set $A$, i.e., $X'_0 \subset X'$. In particular, $X'_0 \subset X_0$.

It remains to show that $X'_0$ is itself an $f$-invariant subset of $X$ containing $x_0$, and clearly 
$x_0 \in X'_0$ since $x_0 = \#(\varnothing)$. Thus let $x \in X'_0$, and so there exists a finite set $A$ with 
$x = \#(A)$. Let $a$ be an element not in $A$;
it then follows that $\#(A \cup \{a\}) = f(\#(A)) = f(x)$, which implies that $f(x) \in X'_0$. Hence 
$X'_0$ is $f$-invariant. 
\eop

\textit{Remark:}\enskip
In the above proof, and also below, we use the fact that for each set $A$ there exists an element $a$ not in $A$. 
In fact there must exist an element in $\mathcal{P}(A) \setminus A$. If this were  not the case then 
$\mathcal{P}(A) \subset A$, and we could define a surjective mapping $f : A \to \mathcal{P}(A)$ by letting 
$f(x) = x$ if $x \in \mathcal{P}(A)$ and $f(x) = \varnothing$ otherwise. But by Cantor's diagonal argument (which 
states that a mapping $f : X \to \mathcal{P}(X)$ cannot be surjective) this is not possible.

Here is the recursion theorem (which first appeared in Dedekind \cite{dedekind}).

\begin{theorem}\label{theorem_cfs_31}
If $(X,f,x_0)$ is a Dedekind system then  for each counting system $(Y,g,y_0)$ there exists a unique mapping
$h : X \to Y$ with $h(x_0) = y_0$ such that $h \circ f = g \circ h$.
\end{theorem}

\proof 
Here we have the unique iterator $\#$ for $(X,f,x_0)$ as well as the unique iterator $\#'$ for $(Y,g,y_0)$.
Let $x \in X$; by Lemma~\ref{lemma_cfs_41} there exists a finite set $A$ with $x = \#(A)$ (since $(X,f,x_0)$ is 
minimal). Moreover, if $B$ is another finite set with $x = \#(B)$ then $\#(B) = \#(A)$ and hence by 
Theorem~\ref{theorem_cfs_21} $B \approx A$ (since $(X,f,x_0)$ is standard). Therefore by 
Proposition~\ref{prop_cfs_21} $\#'(B) = \#'(A)$. This implies  there exists a unique mapping $h : X \to Y$ such 
that $h(\#(A)) = \#'(A)$ for each finite set $A$. In particular, 
$h(x_0) = h(\#(\varnothing)) = \#'(\varnothing) = y_0$. Let $x \in X$; as above there exists a finite set $A$ 
with $x = \#(A)$, and there exists an element $a$ not contained in $A$. Hence 
\begin{eqnarray*}
  h(f(x)) = h(f(\#(A))) &=& h(\#(A \cup \{a\}))\\ 
 &=& \#'(A \cup \{a\}) = g(\#'(A)) = g(h(\#(A))) = g(h(x))
\end{eqnarray*}
and this shows that $h \circ f = g \circ h$. The proof of the uniqueness only uses the fact that $(X,f,x_0)$ is 
minimal: Let $h' : X \to Y$ be a further mapping with $h'(x_0) = y_0$ and $h' \circ f = g \circ h'$ and let
$X' = \{ x \in X : h(x) = h'(x) \}$. Then $x_0 \in X'$, since $h(x_0) = y_0 = h'(x_0)$, and if $x \in X'$ then
$h'(f(x)) = g(h'(x)) = g(h(x)) = h(f(x))$, i.e., $f(x) \in X'$. Thus $X'$ is an $f$-invariant subset of $X$ 
containing $x_0$ and so $X' = X$, i.e., $h' = h$.
\eop

The recursion theorem will now be applied to show that the definition of a finite set we are working with is 
equivalent to the usual one. The usual definition of $A$ being finite is that there exists $n \in \Nat$ and a 
bijective mapping $h : [n] \to A$, where $[0] = \varnothing$ and $[n] = \{0,1,\ldots,n-1\}$ for 
$n \in \Nat \setminus \{0\}$. Moreover, if there is a bijective mapping $h : [n] \to A$, then $n$ is the 
cardinality of $A$, i.e., $n = \#(A)$, and so $A \approx [\#(A)]$ for each finite set $A$. The problem here is 
to assign a meaning to the expression $\{0,1,\ldots,n-1\}$, and one way to do this is to make use of the fact 
that $\{0,1,\ldots,n\} = \{0,1,\ldots,n-1\} \cup \{n\}$ and hence that $[\textsf{s}(n)] = [n] \cup \{n\}$ for 
all $n \in \Nat$. The corresponding approach works with any Dedekind system.

If $X$ is a set then put $\mathcal{F}(X) = \{ A \subset X : \mbox{$A$ is finite} \}$.

\begin{theorem} \label{theorem_cfs_41}
Let $(X,f,x_0)$ be a Dedekind system. Then there exists a unique mapping $[\,\cdot\,] : X \to \mathcal{F}(X)$ 
with $[x_0] = \varnothing$ such that $[f(x)] = [x] \cup \{x\}$ for all $x \in X$. Moreover, $x \notin [x]$ for 
each $x \in X$ and $A \approx [\#(A)]$ holds for each finite set $A$. In particular, a set $A$ is finite if and 
only if $A \approx [x]$ for some $x \in X$.
\end{theorem}

\proof
If $A \in \mathcal{F}(X)$ then by Propositions \ref{prop_fs_31} and \ref{prop_fs_41}~(2)
$\{x_0\} \cup f(A) \in \mathcal{F}(X)$ and so there is a mapping $F : \mathcal{F}(X) \to \mathcal{F}(X)$ given 
by $F(A) = \{x_0\} \cup f(A)$ for all $A \in \mathcal{F}(X)$. (Note that here $F(A)$ is the result of applying 
the mapping $F$ to the argument $A$ while $f(A) = \{ f(a) : a \in A \}$.) Consider the counting system 
$(\mathcal{F}(X),F,\varnothing)$. By Theorem~\ref{theorem_cfs_31} there exists a unique mapping
$[\,\cdot\,] : X \to \mathcal{F}(X)$ with $[x_0] = \varnothing$ and such that 
$[f(x)] = F([x]) = \{x_0\} \cup f([x])$ for all $x \in X$. Let $X_0 = \{ x \in X : [f(x)] = [x] \cup \{x\} \}$. 
Then 
\[ 
[f(x_0)] = \{x_0\} \cup f([x_0]) = \{x_0\} \cup f(\varnothing) = \{x_0\}
 = \varnothing \cup \{x_0\} = [x_0] \cup \{x_0\}
\]
and so $x_0 \in X_0$. Let $x \in X_0$; then
\begin{eqnarray*}
 [f(f(x))] =  \{x_0\} \cup f([f(x)]) &=& \{x_0\} \cup f([x] \cup \{x\})\\
    &=& \{x_0\} \cup f([x]) \cup \{f(x)\} =   [f(x)] \cup \{f(x)\} 
\end{eqnarray*}
and so $f(x) \in X_0$. Thus $X_0$ is an $f$-invariant subset of $X$ containing $x_0$ and hence $X_0 = X$, i.e., 
$[f(x)] = [x] \cup \{x\}$ for all $x \in X$. The uniqueness of the mapping $[\,\cdot\,] : X \to \mathcal{F}(X)$ 
(satisfying $[x_0] = \varnothing$ and $[f(x)] = [x] \cup \{x\}$ for all $x \in X$) follows immediately from the 
fact that $(X,f,x_0)$ is minimal. 

Next consider the set $X_1 = \{ x\in X : x \notin [x] \}$; in particular $x_0 \in X_1$, since 
$x_0 \notin \varnothing = [x_0]$. If $x \in X_1$ then $f(x) \notin f([x])$ (since $f$ is injective) and 
$f(x) \notin \{x_0\}$ (since $x_0 \notin f(X)$) and so $f(x) \notin \{x_0\} \cup f([x]) = [f(x)]$. Hence 
$f(x) \in X_1$. Thus $X_1$ is an $f$-invariant subset of $X$ containing $x_0$ and therefore $X_1 = X$, i.e., 
$x \notin [x]$ for all $x \in X$.

Finally, for each finite set $A$ let $\prop(A)$ be the proposition that $A \approx [\#(A)]$.

($\diamond$)\enskip 
$\prop(\varnothing)$ holds since $[\#(\varnothing)] = [x_0] = \varnothing$.

($\star$)\enskip Let $A$ be a finite set for which $\prop(A)$ holds (so $A \approx [\#(A)]$) and let 
$a \notin A$. Then $A \approx f(A)$, since $f$ is injective, and $x_0 \notin f(A)$, and it therefore follows 
that $\{x_0\} \cup f(A) \approx A \cup \{a\}$. Moreover,
\[
[\#(A \cup \{a\})] = [f(\#(A))] = \{x_0\} \cup f([\#(A)]) = \{x_0\} \cup f(A)
\]
and thus $A \cup \{a\} \approx [\#(A \cup \{a\})]$, i.e., $\prop(A \cup \{a\})$ holds.

Hence by the induction principle for finite sets $\prop(A)$ holds for every finite set $A$, which means  
$A \approx [\#(A)]$ for each finite set $A$. This completes the proof of Theorem~\ref{theorem_cfs_41}.
\eop

We now consider the converse of the recursion theorem (Theorem~\ref{theorem_cfs_31}), and for this it is useful 
to to be more explicit about the structure preserving mappings between counting systems. If $(X,f,x_0)$ and 
$(Y,g,y_0)$ are counting systems then a mapping $\pi : X \to Y$ is said to be a \definition{morphism} from 
$(X,f,x_0)$ to $(Y,g,y_0)$ if $\pi(x_0) = y_0$ and $g\circ \pi = \pi \circ f$. This will also be expressed by 
saying that $\pi : (X,f,x_0) \to (Y,g,y_0)$ is a morphism.

\begin{lemma}\label{lemma_cfs_51}
(1)\enskip
For each counting system $(X,f,x_0)$ the identity mapping $\id_X$ is a morphism from $(X,f,x_0)$ to $(X,f,x_0)$. 

(2)\enskip
If $\pi : (X,f,x_0) \to (Y,g,y_0)$ and $\sigma : (Y,g,y_0) \to (Z,h,z_0)$ are morphisms then $\sigma \circ \pi$ 
is a morphism from $(X,f,x_0)$ to $(Z,h,z_0)$. 
\end{lemma}

\proof 
(1)\enskip
This is clear, since $\id_X(x_0) = x_0$ and $f \circ \id_X = f = \id_X \circ f$.

(2)\enskip
This follows since $(\sigma \circ \pi)(x_0) = \sigma(\pi(x_0)) = \sigma(y_0) = z_0$ and
\[
  h\circ (\sigma \circ \pi) =  (h\circ \sigma) \circ \pi =  (\sigma \circ g) \circ \pi 
=  \sigma \circ (g \circ \pi) =  \sigma \circ (\pi \circ f) = (\sigma \circ \pi) \circ f\;.\ \eop 
\]

If $\pi : (X,f,x_0) \to (Y,g,y_0)$ is a morphism then clearly $\pi \circ \id_X = \pi = \id_Y \circ \pi$, and 
if $\pi,\,\sigma$ and $\tau$ are morphisms for which the compositions are defined then
$(\tau \circ \sigma) \circ \pi = \tau \circ (\sigma \circ \pi)$. This means that counting systems are the 
objects of a concrete category, whose morphisms are those defined above.

A counting system $(X,f,x_0)$ is said to be \definition{initial} if for each counting system $(Y,g,y_0)$ there 
is a unique morphism from $(X,f,x_0)$ to $(Y,g,y_0)$. Theorem~\ref{theorem_cfs_31} thus states that a Dedekind 
system is initial. The following result of Lawvere \cite{lawvere} shows that the converse of the recursion 
theorem holds.

\begin{theorem}\label{theorem_cfs_51}
An initial counting system $(X,f,x_0)$ is a Dedekind system.
\end{theorem}

\proof
We first show that an initial counting system is minimal, and then that it is standard.

\begin{lemma}\label{lemma_cfs_61}
An initial counting system $(X,f,x_0)$ is minimal.
\end{lemma}

\proof 
Let $X_0$ be the least $f$-invariant subset of $X$ containing $x_0$ and $f_0$ be the restriction of $f$ to $X_0$, 
considered as a mapping  $X_0 \to X_0$. Then it is easy to see that the counting system $(X_0,f_0,x_0)$ is minimal
and clearly the inclusion mapping $\mathrm{inc} : X_0 \to X$ defines a morphism from $(X_0,f_0,x_0)$ to 
$(X,f,x_0)$. Let $\sigma : (X,f,x_0) \to (X_0,f_0,x_0)$ be the unique morphism; then
$\mathrm{inc} \circ \sigma = \id_X$, since by Lemma~\ref{lemma_cfs_51} $\mathrm{inc} \circ \sigma$ and $\id_X$ 
are both morphisms from $(X,f,x_0)$ to $(X,f,x_0)$ (and there is only one such morphism, since $(X,f,x_0)$ is 
initial). In particular, $\mathrm{inc}$ is surjective, which implies that $X_0 = X$, i.e., $(X,f,x_0)$ is minimal.
\eop

\begin{lemma}\label{lemma_cfs_71}
An initial counting system $(X,f,x_0)$ is standard.
\end{lemma}

\proof
Let $\diamond$ be an element not contained in $X$, put $X_\diamond = X \cup \{\diamond\}$ and define 
$f_\diamond : X_\diamond \to X_\diamond$ by letting $f_\diamond(x) = f(x)$ for $x \in X$ and 
$f_\diamond(\diamond) = x_0$; thus $(X_\diamond,f_\diamond,\diamond)$ is a counting system. Since $(X,f,x_0)$ is 
initial there exists a unique morphism $\pi : (X,f,x_0) \to (X_\diamond,f_\diamond,\diamond)$. Consider the set 
$X' = \{ x \in X : f_\diamond(\pi(x)) = x \}$; then $x_0 \in X'$, since
$f_\diamond(\pi(x_0)) = f_\diamond(\diamond) = x_0$ and if $x \in X'$ then $f_\diamond(\pi(x)) = x$ and so
\[  f_\diamond(\pi(f(x))) = f_\diamond(f_\diamond(\pi(x))) = f_\diamond(x) = f(x)\;, \]
i.e., $f(x) \in X'$. Thus $X'$ is an $f$-invariant subset of $X$ containing $x_0$ and hence $X' = X$, since by 
Lemma~\ref{lemma_cfs_61} $(X,f,x_0)$ is minimal. Thus $\pi(f(x)) = f_\diamond(\pi(x)) = x$ for all $x \in X$, which 
implies that $f$ is injective. Moreover, $x_0 \notin f(X)$, since
$\pi(f(x)) = f_\diamond(\pi(x)) \ne \diamond = \pi(x_0)$ for all $x \in X$. Hence $(X,f,x_0)$ is standard.
\eop

This completes the proof of Theorem~\ref{theorem_cfs_51}.
\eop


\startsection{Minimal counting systems}

\label{mcs}

In this section we give a more detailed analysis of minimal counting systems, and in particular of those which
are not Dedekind systems. In the following let $(X,f,x_0)$ be a minimal counting system and let $\#$ be the unique 
iterator for $(X,f,x_0)$. For each finite set $A$ define a subset $\ssets{A}$ of $X$ by 
\[ \ssets{A} = \{ x \in X : \mbox{$x = \#(C)$ for some $C \subset A$} \}\;.\]
Thus $\ssets{\varnothing} = \{x_0\}$, $\ssets{\{a\}} = \{x_0,f(x_0)\}$, $\ssets{\{a,b\}} = \{x_0,f(x_0),f(f(x_0))\}$ 
for distinct elements $a$ and $b$, and so on. If $A \subset B$ then clearly $\ssets{A} \subset \ssets{B}$.
A finite set $A$ will be called \definition{$\#$-regular} if $\#(B) \ne \#(A)$ for each proper subset $B$ of $A$.

The first result shows in particular that for minimal counting systems the converse of Theorem~\ref{theorem_cfs_21} 
holds.

\begin{theorem}\label{theorem_mcs_11}
The  following statements are equivalent:
\begin{evlist}{8pt}{6pt}
\item[(1)]  $(X,f,x_0)$ is standard (and thus a Dedekind system).

\item[(2)]  The iterator $\#$ is complete.

\item[(3)]  Each finite set is $\#$-regular.

\item[(4)]  $\ssets{A} \ne X$ for each finite set $A$.

\item[(5)]  The set $X$ is not finite.
\end{evlist}
\end{theorem}

For the proof we need two lemmas:

\begin{lemma}\label{lemma_mcs_11}
(1)\enskip 
For each finite set $A$ the set $\ssets{A}$ is finite.

(2)\enskip
If $A$ and $B$ are finite sets with $A \approx B$ then $\ssets{A} = \ssets{B}$.
\end{lemma}

\proof 
(1)\enskip
By Proposition~\ref{prop_fs_51} the set $\mathcal{P}(A)$ is finite and $\#$ maps $\mathcal{P}(A)$ onto $\ssets{A}$. 
Thus by Proposition~\ref{prop_fs_41}~(2) $\ssets{A}$ is finite.

(2)\enskip
If $h : A \to B$ is a bijection then the mapping $h_* : \mathcal{P}(A) \to \mathcal{P}(B)$ given by $h_*(C) = h(C)$ 
for each $C \in \mathcal{P}(A)$ is also a bijection with $h_*(C) \approx C$ and hence by 
Proposition~\ref{prop_cfs_21} with $\#(h_*(C)) = \#(C)$ for each $C \in \mathcal{P}(A)$. It follows that 
$\ssets{A} = \ssets{B}$.
\eop

\begin{lemma}\label{lemma_mcs_21}
Let $A$ be a finite set and $a \notin A$. Then $A \cup \{a\}$ is $\#$-regular if and only if $\ssets{A} \ne X$.
\end{lemma}

\proof
If $\ssets{A} = X$ then $\#(A \cup \{a\}) \in \ssets{A}$ and so $\#(A \cup \{a\}) = \#(B)$ for some $B \subset A$. 
Hence $A \cup \{a\}$ is not $\#$-regular, since $B$ is a proper subset of $A \cup \{a\}$. Suppose conversely that 
$A \cup \{a\}$ is not $\#$-regular and so there exists a proper subset $B'$ of $A \cup \{a\}$ with 
$\#(B') = \#(A \cup \{a\})$. There then exists $B \subset A$ with $B \approx B'$ and by Proposition~\ref{prop_cfs_21}
$\#(B) = \#(B') = \#(A \cup \{a\})$. Consider the set 
$X_0 = \{ x \in X : \mbox{$x = \#(C)$ for some $C \subset A$} \}$ and so in particular 
$x_0 = \#(\varnothing) \in X_0$. Let $x = \#(C) \in X_0$ with $C \subset A$; if $C \ne A$ then 
$f(x) = f(\#(C)) = \#(C \cup \{c\})$, with $c$ any element in $A \setminus C$, and thus $f(x) \in X_0$. On the other 
hand, if $C = A$ then $f(x) = f(\#(A)) = \#(A \cup \{a\}) = \#(B)$ and again $f(x) \in X_0$, since $B \subset A$. 
Thus $X_0$ is an $f$-invariant subset of $X$ containing $x_0$ and therefore $X_0 = X$, since $(X,f,x_0)$ is minimal. 
But this implies that $\ssets{A} = X$.
\eop

\textit{Proof of Theorem~\ref{theorem_mcs_11}:}\enskip
(1) $\Rightarrow$ (5):\enskip
This follows from Proposition~\ref{prop_fs_41}~(1).

(5) $\Rightarrow$ (4):\enskip
This follows from Lemma~\ref{lemma_mcs_11}~(1).

(4) $\Leftrightarrow$ (3):\enskip
This follows directly from Lemma~\ref{lemma_mcs_21}, since the empty set is clearly $\#$-regular.

(3) $\Leftrightarrow$ (2):\enskip
Suppose $\#$ is complete, let $A$ be a finite set and $B$ be a subset of $A$ with $\#(B) = \#(A)$. Then $B \approx A$ 
and hence by Proposition~\ref{prop_fs_101} $B = A$. Thus $A$ is $\#$-regular. Now suppose conversely that each 
finite set is $\#$-regular and let $A,\, B$ be finite sets with $\#(A) = \#(B)$. Then by Lemma~\ref{lemma_fs_41}, 
and without loss of generality, there exists $C \subset A$ with $C \approx B$ and so by 
Proposition~\ref{prop_cfs_21} $\#(C) = \#(B) = \#(A)$. Thus $C = A$, since $A$ is $\#$-regular, which implies that 
$A \approx B$. This shows that $\#$ is complete.

(3) $\Rightarrow$ (1):\enskip
We assume that  $(X,f,x_0)$ is not standard and show there exists a finite set which is not $\#$-regular. Suppose 
first that $x_0 = f(x)$ for some $x \in X$. By Lemma~\ref{lemma_cfs_41} there exists a finite set $A$ with 
$x = \#(A)$; let $a$ be some element not in $A$. Then $\#(A \cup \{a\}) = f(\#(A)) = f(x) = x_0 = \#(\varnothing)$ 
and $\varnothing$ is a proper subset of $A \cup \{a\}$; hence $A \cup \{a\}$ is not $\#$-regular. Suppose now that 
$f$ is not injective and so there exist $x,\,x' \in X$ with $x \ne x'$ and $f(x) = f(x')$. By 
Lemma~\ref{lemma_cfs_41} there exist finite sets $A$ and $B$ with $x = \#(A)$ and $x' = \#(B)$ and by 
Lemma~\ref{lemma_fs_41} and Proposition~\ref{prop_cfs_21} we can assume that $B \subset A$. Thus $B$ is a proper 
subset of $A$, since $\#(A) = x \ne x' = \#(B)$. Let $a \notin A$; then $B \cup \{a\}$ is a proper subset of 
$A \cup \{a\}$ with $\#(B \cup \{A\}) = f(\#(B)) = f(x') = f(x) = f(\#(A)) = \#(A \cup \{a\})$, and hence again 
$A \cup \{a\}$ is not $\#$-regular.
\eop

If the minimal counting system $(X,f,x_0)$ is not a Dedekind system then by Theorem~\ref{theorem_mcs_11} $X$ is 
finite, and we next look at the structure of such counting systems.

A counting system $(Y,g,y_0)$ will be called \definition{$z$-minimal}, where $z \in Y$, if the only subset $Y'$ of 
$Y$ which contains $y_0$ and is such that $g(Y' \setminus \{z\}) \subset Y'$ is $Y$ itself.

\begin{lemma}\label{lemma_mcs_31}
If $(Y,f,x_0)$ is $z$-minimal then it is also minimal. Moreover, if $A$ is any finite set with $\#(A) = z$ (and such 
a set exists by Lemma~\ref{lemma_cfs_41}) then $\ssets{A} = Y$, and in particular $Y$ is  finite.
\end{lemma}

\proof 
If $Y'$ is a $g$-invariant subset of $Y$ then $g(Y' \setminus \{z\}) \subset g(Y') \subset Y'$. Thus if also
$y_0 \in Y'$ then $Y' = Y$, which shows that $(Y,g,y_0)$ is minimal. Now by Lemma~\ref{lemma_cfs_41} there exists a 
finite set $A$ with $z = \#(A)$, and consider the subset 
$Y' = \{ y \in Y : \mbox{$y = \#(B)$ for some $B \subset A$} \}$; in particular $y_0 = \#(\varnothing) \in Y'$. Let 
$y \in  Y' \setminus \{z\}$; then $y = \#(B)$ for some proper subset $B$ of $A$. Let $a \in A \setminus B$; then 
$B \cup \{a\} \subset A$ and so $g(y) = g(\#(B)) = \#(B \cup \{a\}) \in Y'$. Therefore $Y' = Y$, since $(Y,g,y_0)$ 
is $z$-minimal, which implies that $\ssets{A} = Y$ and in particular that $Y$ is finite.
\eop

\begin{theorem}\label{theorem_mcs_21}
Let $X$ be finite. Then:

(1)\enskip
There exists a unique element $z \in X$ such that $(X,f,x_0)$ is $z$-minimal.

(2)\enskip
$f$ maps $X \setminus \{z\}$ bijectively onto $X \setminus \{x_0\}$.

(3)\enskip
If $A$ is any $\#$-regular finite set $A$ with $\ssets{A} = X$ then $z = \#(A)$. Moreover, a finite set $A$ is 
$\#$-regular and satisfies $\ssets{A} = X$ if and only if $A \approx X \setminus \{c\}$ for some $c \in X$ (i.e., if 
and only if $A$ has one less element than $X$). 

\end{theorem}

\proof 
(1)\enskip
We first show there exists a $\#$-regular finite set $A$ with $\ssets{A} = X$. By Theorem~\ref{theorem_mcs_11} 
there exists a finite set $A'$ such that $\ssets{A'} = X$. Therefore the set 
$\mathcal{S} = \{ C \in \mathcal{P}(A) : \ssets{C} = X \}$ is non-empty and so by Proposition~\ref{prop_fs_111} it 
contains a minimal element $A$, which is in fact $\#$-regular: This holds trivially if $A = \varnothing$, and if $A$ 
is non-empty and $B = A \setminus \{a\}$ with $a \in A$ then $\ssets{B} \ne X$, since $B$ is a proper subset of $A$, 
and hence by Lemma~\ref{lemma_mcs_21} $A = B \cup \{a\}$ is $\#$-regular.

Put $z = \#(A)$; we next show that $(X,f,x_0)$ is $z$-minimal, so let $X' \subset X$ with $x_0 \in X'$ and 
$f(X' \setminus \{z\}) \subset X'$, and consider $\mathcal{K} = \{ B \in \mathcal{P}(A) : \#(B) \in X' \}$. In 
particular, $\varnothing \in \mathcal{K}$, since $\#(\varnothing) = x_0 \in X'$. Now let $B \in \mathcal{K}$ and 
$b \in A \setminus B$; then $B$ is a proper subset of $A$, which means $\#(B) \ne \#(A) = z$, since $A$ is 
$\#$-regular. Hence $\#(B) \in X' \setminus \{z\}$ and it follows that $\#(B \cup \{b\}) = f(\#(B)) \in X'$. Thus 
$\mathcal{K}$ is an inductive system and so $\mathcal{K} = \mathcal{P}(A)$, which implies that 
$X = \ssets{A} \subset X'$, i.e., $X' = X$. Therefore $(X,f,x_0)$ is $z$-minimal.

Suppose $(X,f,x_0)$ is also $z'$-minimal for some $z' \in X$. Then, since $\ssets{A} = X$, there exists $B \subset A$
with $z' = \#(B)$ and by Lemma~\ref{lemma_mcs_31} $\ssets{B} = X$. In particular $z \in \ssets{B}$ and hence 
$\#(A) = z = \#(C)$ for some $C \subset B$. But this is only possible if $C = B = A$, since $A$ is $\#$-regular and 
therefore $z' = \#(B) = \#(A) = z$, which shows there exists a unique element $z \in X$ such that $(X,f,x_0)$ is 
$z$-minimal.

For the other parts we need the following facts:

\begin{lemma}\label{lemma_mcs_41}
Let $A$ be $\#$-regular. Then:

(1)\enskip
A finite set $A'$ with $\#(A') = \#(A)$ is $\#$-regular if and only if $A' \approx A$.

(2)\enskip
$B$ is $\#$-regular for each $B \subset A$.

(3)\enskip
$\#(B) = \#(C)$ holds for subsets $B$ and $C$ of $A$ if and only if $B \approx C$.

(4)\enskip
$\ssets{A} \approx A \cup \{a\}$ for any element $a \notin A$.
\end{lemma}

\proof 
(1)\enskip
Let $A'$ be a finite set with $A' \approx A$ and $h : A \to A'$ be bijective; then, as in the proof of
Lemma~\ref{lemma_mcs_11}~(2), the induced mapping $h_* : \mathcal{P}(A) \to \mathcal{P}(A')$ is also a bijection 
with $\#(h_*(C)) = \#(C)$ for each $C \in \mathcal{P}(A)$. It follows that $A'$ a $\#$-regular. Conversely, let $A$ 
and $A'$ be $\#$-regular with $\#(A') = \#(A)$; by Lemma~\ref{lemma_fs_41} (and without loss of generality) we can 
assume there exists $B' \subset A'$ with $B' \approx A$ and hence by Proposition~\ref{prop_cfs_21} with 
$\#(B') = \#(A) = \#(A')$. Thus $B' = A'$, since $A'$ is $\#$-regular, which shows that $A \approx A'$.

(2)\enskip
Let $B$ be a non-empty subset of $A$ and let $b \in B$; put $A' = A \setminus \{b\}$ and $B' = B \setminus \{b\}$.
Then by Lemma~\ref{lemma_mcs_21} $\ssets{A'} \ne X$, since $A' \cup \{b\} = A$ is $\#$-regular, and 
$\ssets{B'} \subset \ssets{A'}$, since $B' \subset A'$. Thus $\ssets{B'} \ne X$ and so by Lemma~\ref{lemma_mcs_21} 
$B = B' \cup \{b\}$ is $\#$-regular. Since the empty set is also $\#$-regular it follows that every subset of $A$ 
is $\#$-regular.

(3)\enskip
This follows immediately from (1) and (2).

(4)\enskip
Let $a \notin A$ and put $\mathcal{K} = \{ B \in \mathcal{P}(A) : \ssets{B} \approx B \cup \{a\}\}$; in particular
$\varnothing \in \mathcal{K}$, since $\ssets{\varnothing} = \{x_0\} \approx \{a\} = \varnothing \cup \{a\}$. Let 
$B \in \mathcal{K}$ (and so $\ssets{B} \approx B \cup \{a\}$) and let $b \in A \setminus B$; then by (2) 
$B \cup \{b\}$ is $\#$-regular and hence $\#(B \cup \{b\}) \notin \ssets{B}$, since every subset of $B$ is a proper 
subset of $B \cup \{b\}$. Moreover
\begin{eqnarray*}
\ssets{B \cup \{b\}} &=& \{ x \in X : \mbox{$x = \#(C)$ for some $C \subset B \cup \{b\}$} \}\\
    &=& \{ x \in X : \mbox{$x = \#(C)$ for some $C \subset B$} \} \cup \{\#(B \cup \{b\})\}\\ 
    &=& \ssets{B} \cup \{\#(B \cup \{b\})\}
\end{eqnarray*}
since if $C$ is a proper subset of $B \cup \{b\}$ then $C \approx C'$ for some $C' \subset B$, and so by
Proposition~\ref{prop_cfs_21} $\#(C') = \#(C)$. Therefore
\[ \ssets{B \cup \{b\}} = \ssets{B} \cup \{\#(B \cup \{b\}) \approx (B \cup \{a\}) \cup \{b\} 
= (B \cup \{b\}) \cup \{a\}\;,\]
i.e., $B \cup \{b\} \in \mathcal{K}$, which shows that $\mathcal{K}$ is an 
inductive system and hence that $\mathcal{K} = \mathcal{P}(A)$. In particular $A \in \mathcal{K}$, i.e., 
$\ssets{A} \approx A \cup \{a\}$. 
\eop

We now continue with the proof of Theorem~\ref{theorem_mcs_21}.

\medskip

(2)\enskip
Let $A$ be $\#$-regular with $\ssets{A}$, and thus  $z = \#(A)$. Let $x \in X \setminus \{z\}$, and so $x = \#(B)$ for 
some proper subset $B$ of $A$. If $a \in A \setminus B$ then by Lemma~\ref{lemma_mcs_41}~(3)
$f(x) = f(\#(B)) = \#(B \cup \{a\}) \ne x_0$, since $B \cup \{a\} \not\approx \varnothing$, 
and it thus follows that $f(X \setminus \{z\}) \subset X \setminus \{x_0\}$. On the 
other hand, if $x \in X \setminus \{x_0\}$ then by Lemma~\ref{lemma_mcs_41}~(3)
$x = \#(B)$ for some non-empty $B \subset A$. Let $b \in B$
and put $B' = B \setminus \{b\}$; then $x' = \#(B') \in X \setminus \{z\}$, since $B'$ is a proper subset of 
the $\#$-regular set $A$, and $f(x') = f(\#(B')) = \#(B' \cup \{b\}) = \#(B) = x$. Hence 
$f(X \setminus \{z\}) = X \setminus \{x_0\}$. Now let $x_1,\,x_2 \in X \setminus \{z\}$ with 
$f(x_1) = f(x_2)$; there then exist proper subsets $B_1,\,B_2$ of $A$ with $x_1 = \#(B_1)$ and $x_2 = \#(B_2)$ and 
by Lemma~\ref{lemma_fs_41} and Proposition~\ref{prop_cfs_21}, and without loss of generality, we can assume 
$B_2 \subset B_1$. Let $a \in A \setminus B_1$; then 
\[
\#(B_1 \cup \{a\}) = f(\#(B_1)) = f(x_1) = f(x_2) = f(\#(B_2)) = \#(B_2 \cup \{a\})
\]
and so by Lemma~\ref{lemma_mcs_41}~(3) $B_1 \cup \{a\} \approx B_2 \cup \{a\}$, which implies that 
$B_1 \approx B_2$, and then by Proposition~\ref{prop_cfs_21} $x_1 = \#(B_1) = \#(B_2) = x_2$. This shows that the 
restriction of $f$ to $X \setminus \{z\}$ is injective. 

(3)\enskip
The proof of (1) shows there exists a $\#$-regular finite set $A'$ with $\ssets{A'} = X$, and that $z = \#(A)$ for any 
$\#$-regular finite set $A$ with $\ssets{A} = X$. Moreover, if $A$ is such a set then by Lemma~\ref{lemma_mcs_41}~(4) 
$X = \ssets{A} \approx A \cup \{a\}$ for any element $a \notin A$, and hence $A \approx X \setminus \{c\}$ for any 
$c \in X$. Finally, if $A$ is a finite set with $A \approx X \setminus \{c\}$ then $A \approx A'$ (since 
$A' \approx X \setminus \{c\}$), hence by Lemma~\ref{lemma_mcs_41}~(1) and Proposition~\ref{prop_cfs_21} $A$ is 
$\#$-regular and by Lemma~\ref{lemma_mcs_11}  $\ssets{A} = \ssets{A'} = X$.

This completes the proof of Theorem~\ref{theorem_mcs_21}.
\eop

In what follows we assume that $X$ is finite and let $z$ be the unique element given in Theorem~\ref{theorem_mcs_21}
such that $(X,f,x_0)$ is $z$-minimal; $z$ will be called the \definition{end-point} of $(X,f,x_0)$. Note that
$z \ne x_0$ if $X \ne \{x_0\}$ (since $(X,f,x_0)$ being $x_0$-minimal implies that $X = \{x_0\}$).
Theorem~\ref{theorem_mcs_21} essentially characterises the behaviour of $f$ on the set $X \setminus \{z\}$;
it does not say anything about the value $f(z)$.

\begin{proposition}\label{prop_mcs_11}
(1)\enskip
If $f(z) = x_0$ then $f$ is bijective. In this case $x_0 \in f(X)$. 

(2)\enskip
If $f(z) \ne x_0$ then $f$ is not injective and $x_0 \notin f(X)$. 
\end{proposition}

\proof
This follows immediately from Theorem~\ref{theorem_mcs_21}~(2).
\eop

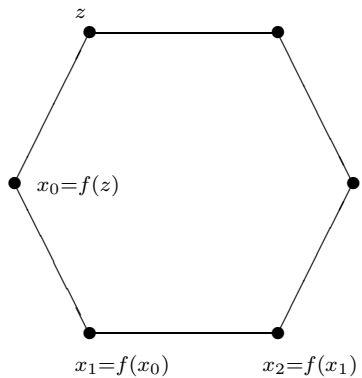
\setlength{\graphicthick}{0.1mm}
\setlength{\graphicmid}{0.1mm}
\setlength{\graphicthin}{0.1mm}

\begin{center}
\setlength{\unitlength}{1.0mm}
\begin{picture}(140,50)

\linethickness{\graphicthick}

\linethickness{\graphicthin}

\linethickness{\graphicmid}
\mythicklines

\put(49,24){$\bullet$}
\sput{53}{24}{x_0 = f(z)}

\sput{58}{0}{x_1 = f(x_0)}
\put(59,4){$\bullet$}

\sput{83}{0}{x_2 = f(x_1)}
\put(84,4){$\bullet$}

\put(59,44){$\bullet$}
\sput{58}{47}{z}

\put(84,44){$\bullet$}
\put(94,24){$\bullet$}

\put(50,25){\line(1,2){10}}
\put(50,25){\line(1,-2){10}}
\put(60,45){\line(1,0){25}}
\put(60,5){\line(1,0){25}}
\put(85,45){\line(1,-2){10}}
\put(85,5){\line(1,2){10}}

\end{picture}

\end{center}

\begin{center}
\setlength{\unitlength}{1.0mm}
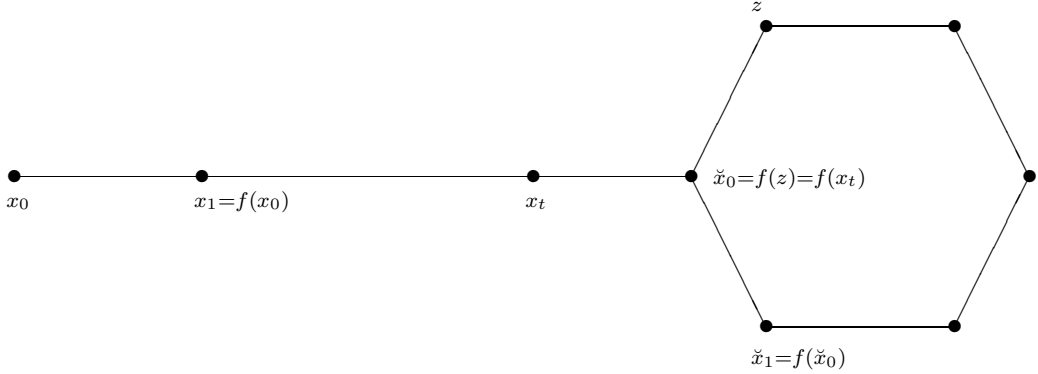
\begin{picture}(140,50)

\linethickness{\graphicthick}

\linethickness{\graphicthin}

\linethickness{\graphicmid}
\mythicklines

\put(5,25){\line(1,0){90}}

\put(95,25){\line(1,2){10}}
\put(95,25){\line(1,-2){10}}
\put(105,45){\line(1,0){25}}
\put(105,5){\line(1,0){25}}
\put(130,45){\line(1,-2){10}}
\put(130,5){\line(1,2){10}}

\sput{4}{21}{x_0}
\put(4,24){$\bullet$}

\sput{29}{21}{x_1 = f(x_0)}
\put(29,24){$\bullet$}

\sput{73}{21}{x_t}
\put(73,24){$\bullet$}

\sput{98}{24}{\breve{x}_0 = f(z) = f(x_t)}
\put(94,24){$\bullet$}

\sput{103}{0}{\breve{x}_1 = f(\breve{x}_0)}
\put(104,4){$\bullet$}

\sput{103}{47}{z}
\put(104,44){$\bullet$}

\put(129,4){$\bullet$}
\put(129,44){$\bullet$}
\put(139,24){$\bullet$}

\end{picture}

\end{center}

\bigskip
Let $w \in X$. Then we can modify $f$ to obtain a new mapping $f_w : X \to X$ by changing the value at the argument 
$z$ from $f(z)$ to $w$ (and leaving all other values unchanged). More precisely, $f_w(x) = f(x)$  for 
$x \in X \setminus \{z\}$ and $f_w(z) = w$. This gives us a new counting system $(X,f_w,x_0)$.

\begin{proposition}\label{prop_mcs_21}
$(X,f_w,x_0)$ is minimal with end-point $z$ for each $w \in X$.
\end{proposition}

\proof 
Let $X' \subset X$ with $x_0 \in X$ and $f_w(X' \setminus \{z\}) \subset X'$; then $f(X' \setminus \{z\}) \subset X'$
also holds and thus $X' = X$, since $(X,f,x_0)$ is $z$-minimal. Hence $(X,f_w,x_0)$ is $z$-minimal, and it follows 
from Lemma~\ref{lemma_mcs_31} that $(X,f_w,x_0)$ is minimal and from Theorem~\ref{theorem_mcs_21} that $z$ is the 
end-point of $(X,f_w,x_0)$.
\eop

\begin{lemma}\label{lemma_mcs_51}
$f(x) \ne x$ for all $x \in X \setminus \{z\}$.
\end{lemma}

\proof 
Let $A$ be a $\#$-regular finite set with $\ssets{A} = X$, and so $\#(A) = z$. If $x \in X \setminus \{z\}$ then 
$x = \#(B)$ for some proper subset $B$ of $A$; let $a \in A \setminus B$. Then $B \cup \{a\} \subset A$ and by 
Proposition~\ref{prop_fs_101} $B \not\approx B \cup \{a\}$. Hence by Lemma~\ref{lemma_mcs_41}~(3)
$f(x) = f(\#(B)) = \#(B \cup \{a\}) \ne \#(B) = x$.
\eop

An element $x \in X$ is a \definition{fixed point} of $f$ if $f(x) = x$.
By Lemma~\ref{lemma_mcs_51} $f$ has a fixed point if and only if $z$ is a fixed point, and in this case 
$z$ is the unique fixed point of $f$. If $z$ is a fixed point then
we say that $(X,f,x_0)$ is a \definition{segment}.

If $(X,f,x_0)$ has end-point $z$ then it follows from Proposition~\ref{prop_mcs_21} that $(X,f_z,x_0)$ is a 
segment. Moreover, given the segment $(X,f_z,x_0)$ and knowing the value $f(z)$ we can recover the original
counting system $(X,f,x_0)$. This means that all the properties of finite minimal counting systems can be deduced
from those of segments.

We typically use $(S,\sigma,s_0)$ (and not $(X,f,x_0)$) to denote a segment. A segment $(S,\sigma,s_0)$ will be 
called \definition{non-trivial} if $S \ne \{s_0\}$. Proposition~\ref{prop_mcs_31} below shows that for each 
non-empty finite set $A$ there is essentially a unique segment $(S,\sigma,s_0)$ with $S \approx A$.

\begin{lemma}\label{lemma_mcs_61}
Let  $(S,\sigma,s_0)$ be a segment with end-point $z$. Then:

(1)\enskip 
Let $S_\diamond = S \cup \{\diamond\}$, with $\diamond$ an element not in $S$, and let 
$\sigma_\diamond : S_\diamond \to S_\diamond$ be the mapping with $\sigma_\diamond(s) = \sigma(s)$ for $s \in S$ and 
$\sigma_\diamond(\diamond) = s_0$. Then $(S_\diamond,\sigma_\diamond,\diamond)$ is also a segment with end-point $z$. 

(2)\enskip
If $(S,\sigma,s_0)$ is non-trivial then the set $S' = S \setminus \{s_0\}$ is $\sigma$-invariant, and if 
$\sigma'$ is the restriction of $\sigma$ to $S'$ then $(S',\sigma',\sigma(s_0))$ is a segment with end-point $z$. 
\end{lemma}

\proof 
(1)\enskip 
Let $S'_\diamond$ be $\sigma_\diamond$-invariant and contain $\diamond$ and let
$S' = S'_\diamond \setminus \{\diamond\}$. Then $S'$ is an $\sigma$-invariant subset of $S$. (If $s \in S'$ then 
$\sigma(s) = \sigma_\diamond(s) \in S'_\diamond$ and so $\sigma(s) \in S'$, since $\sigma(s) \in S$.) Therefore 
$S' = S$, since $s_0 = \sigma_\diamond(\diamond) \in S'$ and $(S,\sigma,s_0)$ is minimal, and thus 
$S'_\diamond = S_\diamond$. Hence $(S_\diamond,\sigma_\diamond,\diamond)$ is minimal. Moreover, by 
Lemma~\ref{lemma_fs_21} $S_\diamond$ is finite, and $z$ is a fixed point of $\sigma_\diamond$. Therefore by 
Lemma~\ref{lemma_mcs_51} $(S_\diamond,\sigma_\diamond,\diamond)$ is a segment with end-point $z$.

(2)\enskip
Proposition~\ref{prop_mcs_11}~(2) implies that $s_0 \notin \sigma(S)$ (since $\sigma(z) = z \ne s_0$), hence
$\sigma(S) \subset S \setminus \{s_0\} = S'$ and in particular $\sigma(S') \subset S'$. Now if $S'_0$ is 
$\sigma'$-invariant and contains $\sigma(s_0)$ then $S'_0 \cup \{s_0\}$ is an $\sigma$-invariant subset of $S$ 
containing $s_0$ and so $S'_0 \cup \{s_0\} = S$, since $(S,\sigma,s_0)$ is minimal. Thus $(S',\sigma',\sigma(s_0))$ 
is minimal. It follows from Lemma~\ref{lemma_mcs_51} that $(S',\sigma',\sigma(s_0))$ is a segment with end-point $z$,
since by Proposition~\ref{prop_fs_21} $S'$ is finite, and $z$ is a fixed point of $\sigma'$.
\eop

\begin{proposition}\label{prop_mcs_31}
(1)\enskip 
For each non-empty finite set $A$ there exists a segment $(S,\sigma,s_0)$ with $S \approx A$.

(2)\enskip
If $(S,\sigma,s_0)$  and $(T,\tau,t_0)$ are segments with $S \approx T$ then there exists a unique mapping 
$p : S \to T$ with $p(s_0) = t_0$ such that $\tau \circ p = p \circ \sigma$.
\end{proposition}

\proof 
(1)\enskip
For each finite set $A$ let $\prop(A)$ be the proposition that if $A$ is non-empty then there exists a segment 
$(S,\sigma,s_0)$ with $S \approx A$.

($\diamond$)\enskip 
$\prop(\varnothing)$ holds trivially.

($\star$)\enskip 
Let $A$ be a finite set for which $\prop(A)$ holds and let $a \notin A$. If $A = \varnothing$ then $(\{a\},\sigma,a)$ 
with $\sigma(a) = a$ is a segment and $\{a\} \approx A \cup \{a\}$, and so in this case $\prop(A \cup \{a\})$ holds.
On the other hand, if $A \ne \varnothing$ then there exists a segment $(S,\sigma,s_0)$ with $S \approx A$, and if
$(S_\diamond,\sigma_\diamond,\diamond)$ is the segment given in Lemma~\ref{lemma_mcs_61}~(1) then 
$S_\diamond = S \cup \{\diamond\} \approx A \cup \{a\}$ and so again $\prop(A \cup \{a\})$ holds.

Therefore by the induction principle for finite sets $\prop(A)$ holds for every finite set $A$, thus for each 
non-empty finite set $A$ there exists a segment $(S,\sigma,s_0)$ with $S \approx A$.

(2)\enskip
For each finite set $A$ let $\prop(A)$ be the proposition that if $A$ is non-empty and $(S,\sigma,s_0)$ and 
$(T,\tau,t_0)$ are segments with $S \approx A \approx T$ then there exists a mapping $p : S \to T$ with 
$p(s_0) = t_0$ such that $\tau \circ p = p \circ \sigma$.

($\diamond$)\enskip 
$\prop(\varnothing)$ holds trivially.

($\star$)\enskip 
Let $A$ be a finite set for which $\prop(A)$ holds and let $a \notin A$. If $A = \varnothing$ then 
$A \cup \{a\} = \{a\}$ and it is clear that $\prop(A \cup \{a\})$ holds here. Assume then that 
$A \ne \varnothing$ and let $(S,\sigma,s_0)$  and $(T,\tau,t_0)$ be segments with $S \approx A \cup \{a\} \approx T$.
 These segments are non-trivial and so by Lemma~\ref{lemma_mcs_61}~(2) we have the segments $(S',\sigma',\sigma(s_0))$ 
and $(T',\tau',\tau(t_0))$, where $S' = S \setminus \{s_0\}$, $T' = T \setminus \{t_0\}$ and $\sigma'$ and $\tau'$ are 
respectively the restrictions of $\sigma$ to $S'$ and $\tau$ to $T'$. Now $S' \approx A \approx T'$ and $\prop(A)$ 
holds and there thus exists a mapping $p' : S' \to T'$ with $p'(\sigma(s_0)) = \tau(t_0)$ such that 
$\tau' \circ p' = p' \circ \sigma'$. Extend $p'$ to a mapping $p : S \to T$ by letting $p(s_0) = t_0$. Then 
$\tau \circ p = p \circ \sigma$, which implies that again $\prop(A \cup \{a\})$ holds.

Therefore by the induction principle for finite sets $\prop(A)$ holds for every finite set $A$, which shows that 
if $(S,\sigma,s_0)$ and $(T,\tau,t_0)$ are segments with $S \approx T$ then there exists a mapping $p : S \to T$ with 
$p(s_0) = t_0$ such that $\tau \circ p = p \circ \sigma$. The uniqueness of $p$ follows as usual from the fact that 
$(S,\sigma,s_0)$ is minimal.
\eop

Let $(S,\sigma,s_0)$  and $(T,\tau,t_0)$ be segments with $S \approx T$. Then by Proposition~\ref{prop_mcs_31}~(2)
there exists a unique mapping $p : S \to T$ with $p(s_0) = t_0$ such that $\tau \circ p = p \circ \sigma$, and the 
mapping $p$ is a bijection: Reversing the roles of $(S,\sigma,s_0)$  and $(T,\tau,t_0)$ gives us a unique mapping 
$q : T \to S$ with $q(t_0) = s_0$ such that $\sigma \circ q = q \circ \tau$, and it easily follows that 
$q \circ p = \id_S$ and $p \circ q = \id_T$. Moreover, if $z$ is the end-point of $(S,\sigma,s_0)$  and $z'$ the 
end-point of $(T,\tau,t_0)$ then $p(z) = z'$, since $\tau(p(z)) = p(\sigma(z)) = p(z)$ and $z'$ is the unique fixed 
point of $\tau$.

Finally, the following result shows how two segments can be joined together to make a `longer' segment.

\begin{proposition}\label{prop_mcs_41}
Let $(S,\sigma,s_0)$ and $(T,\tau,t_0)$ be segments with end-points $z$ and $z'$ respectively, and such that the 
sets $S$ and $T$ are disjoint. Put $R = S \cup T$ and define a mapping $\varrho : R \to R$ by
\[ 
  \varrho(r) = \left\{ \begin{array}{cl}
                    \sigma(r) &\ \mbox{if $r \in S \setminus \{z\}$}\;,\\
                     t_0 &\ \mbox{if $r = z$}\;,\\
                    \tau(r) &\ \mbox{if $r \in T$}\;.
                   \end{array} \right.
\]
Then $(R,\varrho,s_0)$ is a segment with end-point $z'$.
\end{proposition}

\proof 
Let $R_0$ be an $\varrho$-invariant subset of $R$ containing $s_0$, and put $S_0 = R_0 \cap S$ and $T_0 = R_0 \cap T$.
If $s \in S_0 \setminus \{z\}$ then $\sigma(s) = \varrho(s) \in R_0$, and so $\sigma(s) \in S_0$; on the other hand, 
if $z \in S_0$ then $\sigma(z) = z \in S_0$. Thus $S_0$ is a $\sigma$-invariant subset of $S$ containing $s_0$ and 
hence $S_0 = S$, since $(S,\sigma,s_0)$ is minimal. In particular, $z \in S_0 \subset R_0$, which implies that
$t_0 = \varrho(z) \in R_0$, i.e., $t_0 \in T_0$. But $T_0$ is clearly $\tau$-invariant and therefore $T_0 = T$, since 
$(T,\tau,t_0)$ is minimal. It follows that $R_0 = S_0 \cup T_0 = S \cup T = R$, which shows that $(R,\varrho,s_0)$ is 
minimal. But by Proposition~\ref{prop_fs_31} $R$ is finite, and $z'$ is a fixed point of $\varrho$, and so by
Lemma~\ref{lemma_mcs_51} $(R,\varrho,s_0)$ is a segment with end-point $z'$.
\eop


\startsection{Addition and multiplication}

\label{am}

In this section we show how an addition and a multiplication can be defined for any minimal counting system. 
These operations are associative and commutative and can be specified by the rules (a0), (a1), (m0) and (m1) 
below, which are usually employed when defining the operations on $\Nat$ via the Peano axioms.

In the following let $(X,f,x_0)$ be a minimal counting system and let $\#$ be the unique iterator for 
$(X,f,x_0)$. We first state the main results (Theorems \ref{theorem_am_11} and \ref{theorem_am_21}) and then 
develop the machinery required to prove them. In the following section we give alternative proofs for these 
theorems.

\begin{theorem}\label{theorem_am_11}
There exists a unique binary operation $\oplus$ on $X$ such that
\[ 
\#(A) \oplus \#(B) = \#(A \cup B) 
\]
whenever $A$ and $B$ are disjoint finite sets. This operation $\oplus$ is both associative and commutative, 
$x \oplus x_0 = x$ for all $x \in X$ and for all $x_1,\,x_2 \in X$ there is an $x \in X$ such that either 
$x_1 = x_2 \oplus x$ or $x_2 = x_1 \oplus x$. Moreover, $\oplus$ is the unique binary operation $\oplus$ on $X$ 
such that
\begin{evlist}{15pt}{6pt}
\item[(a0)]  $x \oplus x_0 = x$ for all $x \in X$. 

\item[(a1)]  $x \oplus f(x') = f(x \oplus x')$ for all $x,\,x' \in X$.
\end{evlist}
\end{theorem}

\begin{theorem}\label{theorem_am_21}
There exists a unique binary operation $\otimes$ on $X$ such that
\[ 
\#(A) \otimes \#(B) = \#(A\times B) 
\]
for all finite sets $A$ and $B$. This operation $\otimes$ is both associative and commutative, 
$x \otimes x_0 = x_0$ and $x \otimes f(x_0) = x$ for all $x \in X$ (and so $f(x_0)$ is a multiplicative unit) and 
the distributive law holds for $\oplus$ and $\otimes$: For all $x,\,x_1,\,x_2 \in X$
\[     
x \otimes (x_1 \oplus x_2) = (x \otimes x_1) \oplus (x \otimes x_2)\;. 
\]
Moreover, $\otimes$ is the unique binary operation on $X$ such that 
\begin{evlist}{15pt}{6pt}
\item[(m0)]  $x \otimes x_0 = x_0$ for all $x \in X$. 

\item[(m1)]  $x \otimes f(x') = x \oplus (x \otimes x')$ for all $x,\,x' \in X$.
\end{evlist}
\end{theorem}

Before beginning with the proofs of Theorems \ref{theorem_am_11} and \ref{theorem_am_21} we give two results 
about the operation $\oplus$ for special cases of $(X,f,x_0)$.

\begin{proposition}\label{prop_am_11}
If $f$ is injective then the cancellation law holds for $\oplus$ (meaning  that $x_1 = x_2$ whenever 
$x_1 \oplus x = x_2 \oplus x$ for some $x \in X$).
\end{proposition}

\proof 
Let $x_1,\,x_2 \in X$ with $x_1 \ne x_2$ and let $X_0 = \{ x \in X : x_1 \oplus x \ne x_2 \oplus x \}$; then 
$x_0 \in X_0$, since by (a0) $x_1 \oplus x_0 = x_1 \ne x_2 = x_2 \oplus x_0$. Let $x \in X_0$, then by (a1), and 
since $f$ is injective, $x_1 \oplus f(x) = f(x_1 \oplus x) \ne f(x_2 \oplus x) = x_1 \oplus f(x)$, i.e.,
$f(x) \in X_0$. Thus $X_0$ is an $f$-invariant subset of $X$ containing $x_0$ and so $X_0 = X$, since
$(X,f,x_0)$ is minimal. Hence if $x_1 \ne x_2$ then $x_1 \oplus x \ne x_2 \oplus x$ for all $x \in X$, which 
shows that the cancellation law holds for $X$. 
\eop

If $(X,f,x_0)$ is a Dedekind system then $f$ is injective and so by Proposition~\ref{prop_am_11} the cancellation 
law holds for $\oplus$. In particular, $x \ne x \oplus x'$ for all $x,\, x' \in X$ with $x' \ne x_0$ (since 
$x = x \oplus x_0$).

\begin{proposition}\label{prop_am_21}
If $x_0 \in f(X)$ (and so by Theorem~\ref{theorem_mcs_21} and Proposition~\ref{prop_mcs_11} $X$ is finite and 
$f$ is bijective) then $(X,\oplus,x_0)$ is an abelian group: For each $x \in X$ there exists $x' \in X$ such that 
$x \oplus x' = x_0$.
\end{proposition}

\proof 
Let $z$ be the end-point of $(X,f,x_0)$, so by Proposition~\ref{prop_mcs_11}~(1) $f(z) = x_0$ and by 
Theorem~\ref{theorem_mcs_21} $z = \#(A)$ where $A$ is any $\#$-regular finite set with $\ssets{A} = X$. Let 
$x \in X$ and so $x = \#(B)$ for some $B \subset A$. Put $C = (A \setminus B) \cup \{a\}$, where $a$ is some element 
not in $A$ and let $x' = \#(C)$. Then $B$ and $C$ are disjoint and hence 
\[ 
x \oplus x' = \#(B) \oplus \#(C) = \#(B \cup C) = \#(A \cup \{a\}) = f(\#(A)) = f(z) = x_0\;.\ \eop
\]

It is not difficult to show that the group in Proposition~\ref{prop_am_21} is cyclic and generated by the 
element $f(x_0)$. 

We now prepare for the proof of Theorem~\ref{theorem_am_11}. Among other things, the theorem states that there 
exists a binary operation $\oplus$ on $X$ such that
\begin{evlist}{26pt}{6pt}
\item[$\mathrm{(\alpha)}$]
$\#(A) \oplus \#(B) = \#(A \cup B)$ whenever $A$ and $B$ are disjoint finite sets.
\end{evlist}
To understand what this implies it is convenient to introduce some notation. For pairs of finite sets let us 
write $(A,B) \approx (A',B')$ if  $A \approx A'$ and $B \approx B'$, and employ $\#(A,B)$ as shorthand for the 
element $(\#(A),\#(B))$ of $X \times X$. A pair $(A,B)$ will be called \definition{disjoint} if the sets $A$ and 
$B$ are disjoint. Consider disjoint pairs $(A,B)$ and $(A',B')$ with $\#(A,B) = \#(A',B')$. Then 
$\#(A) \oplus \#(B) = \#(A') \oplus \#(B')$ (regardless of how $\oplus$ is defined) which shows that if 
$\,\oplus$ exists satisfying $\mathrm{(\alpha)}$ then the following must hold: 
\begin{evlist}{26pt}{6pt}
\item[$\mathrm{(\beta)}$]
If $(A,B)$ and $(A',B')$ are disjoint pairs with $\#(A,B) = \#(A',B')$ then $\#(A \cup B) = \#(A'\cup B')$.
\end{evlist}
Conversely, $\mathrm{(\beta)}$ is sufficient to ensure the existence of an operation $\oplus$ satisfying 
$\mathrm{(\alpha)}$. To see this is the case, first note the following:

\begin{lemma}\label{lemma_am_11}
For all elements $x,\,x' \in X$ there exists a disjoint pair $(A,B)$ with $(x,x') = \#(A,B)$.
\end{lemma}

\proof 
By Lemma~\ref{lemma_cfs_41} there exists a pair $(C,D)$ with $\#(C,D) = (x,x')$, put 
$A = C \times \{\triangleleft\}$ and $B = D \times \{\triangleright\}$, where $\triangleleft$ and 
$\triangleright$ are distinct elements. Then $(A,B)$ is a disjoint pair with $(A,B) \approx (C,D)$ and so
$\#(A,B) = \#(C,D) = (x,x')$.
\eop

Now suppose that $\mathrm{(\beta)}$ holds and consider disjoint pairs $(A,B)$ and $(A',B')$ with 
$\#(A,B) = (x,x') = \#(A',B')$. Then $\#(A \cup B) = \#(A'\cup B')$ and we can thus define a binary operation 
$\oplus$ on $X$ by letting $x \oplus x' = \#(A \cup B)$, where $(A,B)$ is any disjoint pair with 
$\#(A,B) = (x,x')$. In particular, $\#(A) \oplus \#(B) = \#(A \cup B)$ whenever $A$ and $B$ are disjoint,
i.e., $\mathrm{(\alpha)}$ holds.

The main step in the proof of Theorem~\ref{theorem_am_11} will be to establish that $\mathrm{(\beta)}$ holds. 
For a Dedekind system $(X,f,x_0)$ this is not a problem, since if $(A,B)$ and $(A',B')$ are disjoint pairs with 
$\#(A,B) = \#(A',B')$ then Theorem~\ref{theorem_cfs_21} implies that $(A,B) \approx (A',B')$, from which it 
easily follows that $A \cup B \approx A'\cup B'$ and hence by Proposition~\ref{prop_cfs_21}  that 
$\#(A \cup B) = \#(A'\cup B')$. Establishing $\mathrm{(\beta)}$ for a general minimal counting system involves
more work (which is done in Lemma~\ref{lemma_am_21}).

Once it is known that an operation $\oplus$ exists satisfying $\mathrm{(\alpha)}$ then the remaining properties 
of $\oplus$ listed in Theorem~\ref{theorem_am_11} follow from the corresponding properties of the union operation 
$\cup$ (for example, that it is associative and commutative).

Theorem~\ref{theorem_am_21} will be dealt with in a similar manner. Here we are looking for an operation
$\otimes$ on $X$ such that
\begin{evlist}{26pt}{6pt}
\item[$\mathrm{(\gamma)}$]
$\#(A) \otimes \#(B) = \#(A \times B)$ whenever $A$ and $B$ are finite sets
\end{evlist}
and the condition corresponding to $\mathrm{(\beta)}$ is clearly the following:
\begin{evlist}{26pt}{6pt}
\item[$\mathrm{(\delta)}$]
If $(A,B)$ and $(A',B')$ are (arbitrary) pairs with $\#(A,B) = \#(A',B')$ then
$\#(A \times B) = \#(A'\times B')$.
\end{evlist}
If $\mathrm{(\delta)}$ holds and $(A,B)$ and $(A',B')$ are pairs with $\#(A,B) = (x,x') = \#(A',B')$ then 
$\#(A \times B) = \#(A'\times B')$ and so we can define a binary operation $\otimes$ on $X$ by letting
$x \otimes x' = \#(A \times B)$, where $(A,B)$ is any pair with $\#(A,B) = (x,x')$. In particular, 
$\#(A) \otimes \#(B) = \#(A \times B)$  for all finite sets $A$ and $B$, i.e., $\mathrm{(\gamma)}$ holds. Again, 
there is no problem to show that $\mathrm{(\delta)}$ holds for a Dedekind system, since if 
$(A,B) \approx (A',B')$ then $A \times B \approx A'\times B'$. The general minimal counting system is dealt with 
in Lemma~\ref{lemma_am_31}.

As with the addition $\oplus$, once it is known that an operation $\otimes$ exists satisfying $\mathrm{(\gamma)}$ 
then the remaining properties of $\otimes$ listed in Theorem~\ref{theorem_am_21} follow from the corresponding 
properties of the cartesian product operation $\times$ (for example, that it is associative and commutative) and 
the relationship between $\cup$ and $\times$.

The following shows that condition $\mathrm{(\beta)}$ holds.

\begin{lemma}\label{lemma_am_21}
If $(A,B)$ and $(A',B')$ are disjoint pairs with $\#(A,B) = \#(A',B')$ then $\#(A \cup B) = \#(A'\cup B')$.
\end{lemma}

\proof 
Start by considering a finite sets $A$ and $A'$ with $\#(A) = \#(A')$ and for each finite set $B$ let $\prop(B)$ 
be the proposition that if $B$ is disjoint from $A$ and $A'$ then $\#(A \cup B) = \#(A' \cup B)$.

($\diamond$)\enskip 
$\prop(\varnothing)$ holds, since $\#(A \cup \varnothing) = \#(A) = \#(A') = \#(A' \cup \varnothing)$.

($\star$)\enskip 
Let $B$ be a finite set for which $\prop(B)$ holds and $b \notin B$. If $B \cup \{b\}$ is not disjoint from $A$ 
and $A'$ then $\prop(B \cup \{b\})$ holds trivially, and so we can assume that this is not the case. In 
particular, $B$ is then disjoint from $A$ and $A'$ and so $\#(A \cup B) = \#(A' \cup B)$, and also 
$b \notin A \cup B$ and $b \notin A' \cup B$. Thus
\begin{eqnarray*}
\#(A \cup (B \cup \{b\})) &=& \#((A \cup B) \cup \{b\}) =  f(\#(A \cup B))\\ 
&=& f(\#(A' \cup B)) = \#((A' \cup B) \cup \{b\}) = \#(A' \cup (B \cup \{b\})) 
\end{eqnarray*}
and hence $\prop(B \cup \{b\})$ holds.

Therefore by the induction principle for finite sets $\prop(B)$ holds for every finite set $B$, thus if $A$ and 
$A'$ are finite sets with $\#(A) = \#(A')$ then $\#(A \cup B) = \#(A' \cup B)$ for every finite set $B$ disjoint 
from $A$ and $A'$.

For each set $C$ and each element $d$ put $C_d = C \times \{d\}$ (and so $C_d \approx C$). Now let $(A,B)$ and 
$(A',B')$ be disjoint pairs with $\#(A,B) = \#(A',B')$, and choose distinct elements $\triangleleft$ and 
$\triangleright$; then $\#(A \cup B ) = \#(A_\triangleleft \cup B_\triangleleft)$ (since 
$A \cup B \approx A_\triangleleft \cup B_\triangleleft$), 
$\#(A'_\triangleright \cup B'_\triangleright) = \#(A' \cup B')$ 
(since $A'_\triangleright \cup B'_\triangleright \approx A' \cup B'$),
$\#(A_\triangleleft) = \#(A'_\triangleright)$ (since 
$A_\triangleleft \approx A$ and  $A' \approx A'_\triangleright)$)
and $\#(B_\triangleleft) = \#(B'_\triangleright)$ (since
$B_\triangleleft \approx B$ and  $B' \approx B'_\triangleright$), which gives us the following data: 
\begin{evlist}{16pt}{6pt}
\item[--] $\#(A \cup B) = \#(A_\triangleleft \cup B_\triangleleft)$, 
\item[--] $\#(A_\triangleleft) = \#(A'_\triangleright)$ and $B_\triangleleft$ is disjoint from both $A_\triangleleft$ 
and $A'_\triangleright$,
\item[--] $\#(B_\triangleleft) = \#(B'_\triangleright)$ and $A'_\triangleright$ is disjoint from both $B_\triangleleft$ 
and $B'_\triangleright$,
\item[--] $\#(A'_\triangleright \cup B'_\triangleright) = \#(A' \cup B')$. 
\end{evlist}
Thus by two applications of the first part of the proof
\begin{eqnarray*} 
\#(A \cup B) &=& \#(A_\triangleleft \cup B_\triangleleft) = \#(A'_\triangleright \cup B_\triangleleft)\\
&=& \#(B_\triangleleft \cup A'_\triangleright) = \#(B'_\triangleright \cup A'_\triangleright)
= \#(A'_\triangleright \cup B'_\triangleright) = \#(A' \cup B')\;.\eop
\end{eqnarray*}

\textit{Proof of Theorem~\ref{theorem_am_11}:}\enskip
Let $x,\,x' \in X$. By Lemma~\ref{lemma_am_11} there exists a disjoint pair $(A,B)$ with $(x,x') = \#(A,B)$ and 
if $(A',B')$ is another disjoint pair $(A',B')$ with $(x,x') = \#(A',B')$ then $\#(A,B) = \#(A',B')$ and so 
Lemma~\ref{lemma_am_21} implies that $\#(A \cup B) = \#(A' \cup B')$. We can therefore define $x \oplus x'$ to 
be $\#(A \cup B)$, where $(A,B)$ is any disjoint pair with  $(x,x') = \#(A,B)$. In particular, if $A$ and $B$ 
are disjoint finite sets then $\#(A) \oplus \#(x') = \#(A \cup B)$, since $(A,B)$ is a disjoint pair for
$(\#(A),\#(B))$. Moreover, $\oplus$ is uniquely determined by this requirement: Consider any binary operation 
$\oplus'$ on $X$ for which $\#(A) \oplus' \#(B) = \#(A \cup B)$ whenever $A$ and $B$ are disjoint finite sets. If 
$C$ and $D$ are any finite sets then there exists a disjoint pair $(A,B)$ with $(A,B) \approx (C,D)$ and hence by 
Proposition~\ref{prop_cfs_21} 
\[ 
\#(C) \oplus' \#(D) = \#(A) \oplus' \#(B) = \#(A \cup B) = \#(A) \oplus \#(B) = \#(C) \oplus \#(D) 
\] 
and so by Lemma~\ref{lemma_cfs_41} ${\oplus'} = {\oplus}$.

We show that $\oplus$ is associative and commutative: Let $x_1,\,x_2,\,x_3 \in X$; then by 
Lemma~\ref{lemma_cfs_41} there exists finite sets $A_1,\,A_2$ and $A_3$ with $x_1 = \#(A_1)$, $x_2 = \#(A_2)$ and 
$x_3 = \#(A_3)$. Let $\triangleleft,\,\triangleright$ and $\diamond$ be distinct elements and put
$B_1 = A_1 \cup \{\triangleleft\}$, $B_2 = A_2 \cup \{\triangleright\}$ and $B_3 = A_3 \cup \{\diamond\}$. Then 
$B_1,\,B_2,\,B_3$ are disjoint with $B_j \approx A_j$ and hence with $\#(B_j) = \#(A_j)$ for $j = 1,\,2,\,3$.
Therefore 
\begin{eqnarray*}
(x_1 \oplus x_2) \oplus x_3 &=& (\#(B_1) \oplus \#(B_2)) \oplus \#(B_3) 
 = \#(B_1 \cup B_2) \oplus \#(B_3) \\
 &=& \#((B_1 \cup B_2) \cup B_3) = \#(B_1 \cup (B_2 \cup B_3)) \\
 &=& \#(B_1) \oplus \#(B_2 \cup B_3) = \#(B_1) \oplus (\#(B_2) \oplus \#(B_3))\\ 
&=& x_1 \oplus (x_2 \oplus x_3) \;.
\end{eqnarray*}
In the same way $\oplus$ is commutative. Let $x_1,\,x_2 \in X$; then there exists a disjoint pair $(A_1,A_2)$ with 
$(x_1,x_2) = \#(A_1,A_2)$. Thus
\begin{eqnarray*}
x_1 \oplus x_2 &=& \#(A_1) \oplus \#(A_2)\\ 
&=& \#(A_1 \cup A_2) = \#(A_2 \cup A_1) = \#(A_2) \oplus \#(A_1) = x_2 \oplus x_1 \;.
\end{eqnarray*}
Moreover, if $x \in X$ and $A$ is a finite set with $x = \#(A)$ then
\[ 
x \oplus x_0 = \#(A) \oplus \#(\varnothing)  = \#(A \cup \varnothing) = \#(A) =  x\,,
\]
and so $x \oplus x_0 = x$ for all $x \in X$.

Let $x_1,\,x_2 \in X$, and so by Lemma~\ref{lemma_cfs_41} there exist finite sets $A$ and $B$ such that
$x_1 = \#(A)$ and $x_2 = \#(B)$. By Proposition~\ref{prop_fs_91} there either exists an injective mapping 
$g : A \to B$ or an injective mapping $h : B \to A$. Assume the former holds and put $B' = g(A)$ and 
$C = B \setminus B'$. Then $B'$ and $C$ are disjoint and $B = B' \cup C$; moreover, $A \approx B'$ (since $g$ 
considered as a mapping from $A$ to $B'$ is a bijection) and so by Proposition~\ref{prop_cfs_21} $\#(A) = \#(B')$.
Thus, putting $x = \#(C)$, it follows that
$x_2 = \#(B) = \#(B' \cup C) = \#(B') \oplus \#(C) = \#(A) \oplus \#(C) = x_1 \oplus x$. On the other hand, if 
there exists an injective mapping $h : B \to A$ then the same argument shows that $x_1 = x_2 \oplus x$ for some 
$x \in X$.

Now to (a0) and (a1), and we have seen above that (a0) holds. Let $x,\,x' \in X$, so by Lemma~\ref{lemma_am_11} 
there exists a disjoint pair $(A,B)$ for $(x,x')$.  Let $b \notin A \cup B$; then 
\begin{eqnarray*}
x \oplus f(x') &=& \#(A) \oplus f(\#(B)) = \#(A) \oplus \#(B \cup \{b\}) = \#(A \cup (B \cup \{b\})) \\
 &=& \#((A \cup B) \cup \{b\}) =  f(\#(A \cup B)) = f(\#(A) \oplus \#(B)) = f(x \oplus x') 
\end{eqnarray*}
and hence (a1) holds. If $\oplus'$ is another binary operation on $X$ satisfying (a0) and (a1) then it is easy 
to see that $X_0 = \{ x' \in X : x \oplus' x' = x \oplus x'\ \mbox{for all $x \in X$} \}$ is an $f$-invariant 
subset of $X$ containing $x_0$. Hence $X_0 = X$, since $(X,f,x_0)$  is minimal, which implies that 
${\oplus'} = {\oplus}$.

This completes the proof of Theorem~\ref{theorem_am_11}.
\eop

We prepare for the proof of Theorem~\ref{theorem_am_21} by showing that $\mathrm{(\delta)}$ holds.

\begin{lemma}\label{lemma_am_31}
If $(A,B)$ and $(A',B')$ are any pairs with $\#(A,B) = \#(A',B')$ then $\#(A \times B) = \#(A' \times  B')$.
\end{lemma}

\proof 
Start by considering finite sets $A$ and $A'$ with $\#(A) = \#(A')$ and for each finite set $B$ let $\prop(B)$ 
be the proposition that $\#(A \times B) = \#(A' \times B)$.

($\diamond$)\enskip 
$\prop(\varnothing)$ holds, since $A \times \varnothing = \varnothing = A' \times \varnothing$ and so 
$\#(A \times \varnothing) = \#(A' \times \varnothing)$.

($\star$)\enskip 
Let $B$ be a finite set for which $\prop(B)$ holds and let $b \notin B$. Then the sets $A \times B$ and 
$A \times \{b\}$ are disjoint and $A \times (B \cup \{b\}) = (A \times B) \cup (A \times \{b\})$. It thus follows 
that
\[ 
\#(A \times (B \cup \{b\})) = \#((A \times B) \cup (A \times \{b\}) = \#(A \times B) \oplus \#(A \times \{b\}) 
\]
and in the same way $\#(A' \times (B \cup \{b\})) = \#(A' \times B) \oplus \#(A' \times \{b\})$. Clearly 
$A \times \{b\} \approx A$ and so by Proposition~\ref{prop_cfs_21} $\#(A \times \{b\}) = \#(A)$, and in the same 
way $\#(A' \times \{b\}) = \#(A')$. Therefore
\begin{eqnarray*}
\#(A \times (B \cup \{b\})) &=& \#(A \times B) \oplus \#(A \times \{b\}) = \#(A' \times B) \oplus \#(A)\\ 
&=& \#(A' \times B) \oplus \#(A' \times \{b\}) = \#(A' \times (B \cup \{b\})) 
\end{eqnarray*}
and so $\prop(B \cup \{b\})$ holds.

Therefore by the induction principle for finite sets $\prop(B)$ holds for every finite set $B$, thus if 
$\#(A) = \#(A')$ then $\#(A \times B) = \#(A' \times B)$ for every finite set $B$.

Now let $(A,B)$ and $(A',B')$ be any pairs with $\#(A,B) = \#(A',B')$. Then clearly we have
$A' \times B \approx B \times A'$ and $A' \times B' \approx B' \times A'$ and therefore by 
Proposition~\ref{prop_cfs_21} $\#(A' \times B) = \#(B \times A')$ and $\#(A' \times B') = \#(B' \times A')$.
Hence by the first part 
\[ 
\#(A \times B) = \#(A' \times B) = \#(B \times A') = \#(B' \times A') = \#(A' \times B')\;.\eop
\]

\textit{Proof of Theorem~\ref{theorem_am_21}:}\enskip
Let $x,\,x' \in X$. By Lemma~\ref{lemma_cfs_41} there exist finite sets $A$ and $B$ with $x = \#(A)$ and 
$x' = \#(B)$ and if $x = \#(A')$ and $x' = \#(B')$ for some other finite sets $A'$ and $B'$ then $\#(A') = \#(A)$
and $\#(B') = \#(B)$ and hence by Lemma~\ref{lemma_am_31} $\#(A \times B) = \#(A' \times B')$. We can thus define 
$x \otimes x'$ to be $\#(A \times B)$, where $A$ and $B$ are any finite sets with $x = \#(A)$ and $x' = \#(B)$. 
Then $\#(A) \otimes \#(B) = \#(A \times B)$ for all finite sets $A$ and $B$, a requirement which clearly 
determines $\otimes$ uniquely.

We show that $\otimes$ is associative and commutative: Let $x_1,\,x_2,\,x_3 \in X$; then by 
Lemma~\ref{lemma_cfs_41} there exists finite sets $A_1,\,A_2,\,A_3$ with $x_1 = \#(A_1)$, $x_2 = \#(A_2)$ and 
$x_3 = \#(A_3)$. Now it is easy to check that $(A_1 \times A_2) \times A_3 \approx A_1 \times (A_2 \times A_3)$
and so by Proposition~\ref{prop_cfs_21} $\#((A_1 \times A_2) \times A_3) = \#(A_1 \times (A_2 \times A_3))$.
Therefore 
\begin{eqnarray*}
(x_1 \otimes x_2) \otimes x_3 &=& (\#(A_1) \otimes \#(A_2)) \otimes \#(A_3) = \#(A_1 \times A_2) \otimes \#(A_3) \\
 &=& \#((A_1 \times A_2) \times A_3) = \#(A_1 \times (A_2 \times A_3)) \\
 &=& \#(A_1) \otimes \#(A_2 \times A_3) = \#(A_1) \otimes (\#(A_2) \otimes \#(A_3))\\ 
&=& x_1 \otimes (x_2 \otimes x_3) 
\end{eqnarray*}
which shows $\otimes$ is associative. Let $x_1,\,x_2 \in X$; by Lemma~\ref{lemma_cfs_41} there exist finite sets 
$A_1$ and $A_2$ with $x_1 = \#(A_1)$ and $x_2 = \#(A_2)$. Then by Proposition~\ref{prop_cfs_21} 
$\#(A_1 \times A_2) = \#(A_2 \times A_1)$, since clearly $A_1 \times A_2 \approx A_2 \times A_1$. Thus
\begin{eqnarray*}
x_1 \otimes x_2 &=& \#(A_1) \otimes \#(A_2)\\ 
&=& \#(A_1 \times A_2) = \#(A_2 \times A_1) = \#(A_2) \otimes \#(A_1) = x_2 \otimes x_1 
\end{eqnarray*}
which shows that $\otimes$ is also commutative. 

Let $x \in X$, so by Lemma~\ref{lemma_cfs_41} there exists a finite set $A$ with $x = \#(A)$. Then
\[
x \otimes x_0 = \#(A) \otimes \#(\varnothing) = \#(A \times \varnothing) = \#(\varnothing) = x_0\;.
\]
Moreover, if $a$ is any element then by Proposition~\ref{prop_cfs_21} $\#(A \times \{a\}) = \#(A)$, since 
$A \times \{a\} \approx A$, and hence
\begin{eqnarray*}
 x \otimes f(x_0) &=& \#(A) \otimes f(\#(\varnothing)) = \#(A) \otimes \#(\varnothing \cup \{a\})\\
&=& \#(A) \otimes \#(\{a\}) = \#(A \times \{a\}) = \#(A) = x \;.
\end{eqnarray*}
Thus $x \otimes x_0 = x_0$ and $x \otimes f(x_0) = x$ for each $x \in X$ (and note that the first statement is 
(m0)).

Now for the distributive law. Let $x,\,x_1,\,x_2 \in X$. There exists a finite set $A$ with $x = \#(A)$ and 
a disjoint pair $(B,C)$ with $(x_1,x_2) = \#(B,C)$. Then $A \times (B \cup C)$ is the disjoint union of 
$A \times B$ and $A \times C$ and thus
\begin{eqnarray*}
(x \otimes x_1) \oplus (x \otimes x_2) &=&  (\#(A) \otimes \#(B)) \oplus (\#(A) \otimes \#(C))\\ 
&=&  \#(A \times B) \oplus \#(A \times C) =  \#((A \times B) \cup (A \times C))\\ 
&=&  \#(A \times (B \cup C)) = \#(A) \otimes \#(B \cup C)\\ 
&=& \#(A) \otimes (\#(B) \oplus \#(C)) =  x \otimes (x_1 \oplus x_2)\;.
\end{eqnarray*}

We have already seen that (m0) holds and, since $f(x_0)$ is a unit, (m1) is a special case of the distributive 
law: Let $x,\,x' \in X$; then by (a0) and (a1) and since $\oplus$ is commutative it follows that
$f(x') = f(x' \oplus x_0) = x' \oplus f(x_0) = f(x_0) \oplus x'$, and hence 
$x \otimes f(x') = x \otimes (f(x_0) \oplus x')
= (x \otimes f(x_0)) \oplus (x \otimes x') = x \oplus (x \otimes x')$, which is (m1). Finally, if $\otimes'$ is 
another binary operation satisfying (m0) and (m1) then it is easy to see that 
$X_0 = \{ x' \in X : x \otimes' x' = x \otimes x'\ \mbox{for all $x \in X$} \}$
is an $f$-invariant subset of $X$ containing $x_0$. Hence $X_0 = X$, since $(X,f,x_0)$  is minimal, which implies 
that ${\otimes'} = {\otimes}$.
\eop

We end the section by looking at the operation of exponentiation. Here we have to be more careful: For example, 
$2 \cdot 2 \cdot 2 = 2$ in $\Int_3$ and so $2^3$ is not well-defined if the exponent $3$ is considered as an 
element of $\Int_3$ (since we would also have to have $2^0 = 1$). However, $2^3$ does make sense if $2$ is 
considered as an element of $\Int_3$ and the exponent $3$ as an element of $\Nat$.

In general we will see that if $(Y,g,y_0)$ is a Dedekind system then we can define an element of $X$ which is 
`$x$ to the power of $y$' for each $x \in X$ and each $y \in Y$ and this operation has the properties which might 
be expected. 

In what follows let $(Y,g,y_0)$ be a Dedekind system and let $\#'$ be the unique iterator for $(Y,g,y_0)$. (As 
before $(X,f,x_0)$ is assumed to be minimal with unique iterator $\#$.)

\begin{theorem}\label{theorem_am_31}
There exists a unique operation ${\uparrow} : X \times Y \to X$ such that  
\[ 
\#(A) \uparrow \#'(B) = \#(A^B) 
\]
for all finite sets $A$ and $B$. This operation ${\uparrow}$ satisfies
\[     
x \uparrow (y_1 \oplus y_2) = (x \uparrow y_1) \otimes (x \uparrow y_2) 
\]
for all $x \in X$ and all $y_1,\,y_2 \in Y$ and
\[     
(x_1 \otimes x_2) \uparrow y = (x_1 \uparrow y) \otimes (x_2 \uparrow y) 
\]
for all $x_1,\,x_2 \in X$ and $y \in Y$. Moreover, ${\uparrow}$ is the unique operation such that
\begin{evlist}{15pt}{6pt}
\item[(e0)]  $x \uparrow y_0 = f(x_0)$ for all $x \in X$. 

\item[(e1)]  $x \uparrow g(y) = x \otimes (x \uparrow y)$ for all $x \in X$, $y \in Y$.
\end{evlist}
\end{theorem}

\begin{lemma}\label{lemma_am_41}
If $B,\,C$ are finite sets with $\#(B) = \#(C)$ then for all finite sets $A$ we have $\#(B^A) = \#(C^A)$. 
\end{lemma}

\proof
Let $B,\, C$ be finite sets with $\#(B) = \#(C)$ and for each finite set $A$ let $\prop(A)$ be the proposition 
that $\#(B^A) = \#(C^A)$.

($\diamond$)\enskip 
$\prop(\varnothing)$ holds since $\#(B^\varnothing) =\#(C^\varnothing)$. (For any set $X$ the set $X^\varnothing$ 
consists of the single element $\{\varnothing\}$.)

($\star$)\enskip Let $A$ be a finite set for which $\prop(A)$ holds and $a \notin A$. Now $\#(B^A) = \#(C^A)$ 
(since $\prop(A)$ holds) and $\#(B) = \#(C)$; hence by Proposition~\ref{prop_cfs_21} and Lemma~\ref{lemma_am_31}
\[ 
\#(B^{A \cup \{a\}}) = \#(B^A \times B)  = \#(C^A \times C) = \#(C^{A \cup \{a\}}) 
\]
(since $E^{A\cup \{a\}} \approx E^A \times E$ for each set $E$), and so $\prop(A \cup \{a\})$ holds.

Therefore by the induction principle for finite sets $\prop(A)$ holds for every finite set $A$, and so 
$\#(B^A) = \#(C^A)$ holds for all finite sets $A$.
\eop

\textit{Remark:}\enskip
If $B,\, C$ are finite sets with $\#(B) = \#(C)$ then $\#(A^B) = \#(A^C)$ does not hold in general for a finite 
set $A$.

\medskip
\textit{Proof of Theorem~\ref{theorem_am_31}:}\enskip
Let $A_1,\,A_2,\,B_1,\,B_2$ be finite sets with $\#(A_1) = \#(A_2)$ and $\#'(B_1) = \#'(B_2)$; then by 
Lemma~\ref{lemma_am_41} $\#(A_1^{B_1}) = \#(A_2^{B_1})$ and by Theorem~\ref{theorem_cfs_21} $B_1 \approx B_2$. 
Since $B_1 \approx B_2$ it follows that $A_2^{B_1} \approx A_2^{B_2}$ and then by Proposition~\ref{prop_cfs_21}
$\#(A_2^{B_1}) = \#(A_2^{B_2})$. This shows that $\#(A_1^{B_1}) = \#(A_2^{B_2})$. Therefore by 
Lemma~\ref{lemma_cfs_41} we can define $x \uparrow y$ to be $\#(A^B)$, where $A$ and $B$ are any finite sets 
with $x = \#(A)$ and $y = \#'(B)$. Then
\[ 
\#(A) \uparrow \#'(B) = \#(A^B) 
\]
for all finite sets $A$ and $B$ and this requirement clearly determines ${\uparrow}$ uniquely.
 
Let $x \in X$ and $y_1,\,y_2 \in Y$; then by Lemma~\ref{lemma_am_11} there exists a disjoint pair $(B_1,B_2)$ 
with  $(y_1,y_2) = \#'(B_1,B_2)$ and by Lemma~\ref{lemma_cfs_41} there exists a finite set $A$ with $x = \#(A)$.
Moreover, it is easily checked that $A^{B_1 \cup B_2} \approx A^{B_1} \times A^{B_2}$ and thus by 
Proposition~\ref{prop_cfs_21}
\begin{eqnarray*}
x \uparrow (y_1 \oplus y_2)  
  &=& \#(A) \uparrow (\#'(B_1) \oplus \#(B_2))\\ 
&=& \#(A) \uparrow \#'(B_1 \cup B_2) = \#(A^{B_1 \cup B_2}) = \#(A^{B_1} \times A^{B_2})\\
&=& \#(A^{B_1}) \otimes \#(A^{B_2}) = (x \uparrow y_1) \otimes (x \uparrow y_2) \;.
\end{eqnarray*}
Now let $x_1,\,x_2 \in X$ and $y \in Y$. By Lemma~\ref{lemma_cfs_41} there exist finite sets $A_1,\,A_2$ and $B$ 
such that $x_1 = \#(A_1)$, $x_2 = \#(A_2)$ and $y = \#'(B)$ and $(A_1\times A_2)^B  \approx A_1^B \times A_2^B$. 
Thus by Proposition~\ref{prop_cfs_21}
\begin{eqnarray*}
(x_1 \otimes x_2) \uparrow y  &=& (\#(A_1) \otimes \#(A_2)) \uparrow \#'(B)\\ 
&=& \#(A_1\times A_2) \uparrow \#'(B) = \#((A_1\times A_2)^B) = \#(A_1^B \times A_2^B)\\
&=& \#(A_1^B) \otimes \#(A_2^B) = (x_1 \uparrow y) \otimes (x_2 \uparrow y) \;.
\end{eqnarray*}
It remains to consider the properties (e0) and (e1). Now for each finite set $A$ we have
$\#(A) \uparrow \#'(\varnothing) = \#(A^\varnothing)= \#(\{\varnothing\}) = f(x_0)$ and hence 
$x \uparrow y_0 = f(x_0)$ for each $x \in X$, i.e., (e0) holds. Let $A$ and $B$ be finite sets and let 
$b \notin B$. Then, since $A^{B\cup \{b\}} \approx A \times A^B$, it follows from Proposition~\ref{prop_cfs_21} that
\begin{eqnarray*}
\#(A) \uparrow g(\#'(B)) &=& \#(A) \uparrow \#'(B\cup \{b\}) = \#(A^{B\cup \{b\}}) = \#(A \times A^B)\\ 
&=& \#(A) \otimes \#(A^B) = \#(A) \otimes (\#(A) \uparrow \#'(B)) 
\end{eqnarray*}
and this shows 
$x \uparrow g(y) = x \otimes (x \uparrow y)$ for all $x \in X$, $y \in Y$, i.e., (e1) holds. Finally, if 
${\uparrow}'$ is another operation satisfying (e0) and (e1) then 
\[
Y_0 = \{ y \in Y : x \uparrow' y = x \uparrow y\ \mbox{for all $x \in X$} \}
\]
is a $g$-invariant subset of $Y$ containing $y_0$. Therefore $Y_0 = Y$, since $(Y,g,y_0)$  is minimal, which 
implies that ${\uparrow'} = {\uparrow}$.
\eop


\startsection{Another take on addition and multiplication}

\label{mam}

In the following again let $(X,f,x_0)$ be a minimal counting system and let $\#$ be the unique iterator for 
$(X,f,x_0)$. In this section we give alternative proofs for Theorems \ref{theorem_am_11} and \ref{theorem_am_21}. 

In Section~\ref{am} only the single iterator $\#$ was used. Here we make use of a family of iterators 
$\{ \#_x : x \in X \}$, which arise as follows: For each $x \in X$ there is the counting system $(X,f,x)$ (which 
will usually not be minimal) and by Theorem~\ref{theorem_cfs_11} there is then the unique iterator for $(X,f,x)$, 
which we denote by $\#_x$. Thus $\#_x(\varnothing) = x$ and $\#_x(A \cup \{a\}) = f(\#_x(A))$ whenever $A$ is a 
finite set and $a \notin A$. In particular we have $\# = \#_{x_0}$. For each finite set $A$ the element $\#_x(A)$ 
can be thought of as the `number' obtained by counting the elements in $A$, but starting with $x$ instead of $x_0$. 
Now it is more convenient to repackage the information given by the iterators $\#_A$, $x \in X$, by introducing 
for each finite set $A$ the mapping $f_A : X \to X$ with $f_A(x) = \#_x(A)$ for all $x \in X$, and so 
$\#(A) = \#_{x_0}(A) = f_A(x_0)$.

Consider disjoint finite sets $A$ and $B$; then $\#(A \cup B)$ gives the `number' of elements in $A \cup B$. But 
this `number' can also be determined by first counting the elements in $B$, giving the result $\#(B)$, and then 
counting the elements in $A$, but starting the counting with $\#(B)$ and not with $x_0$. The result is thus 
$\#_z(A)$, where $z = \#(B)$, and $\#_z(A) = f_A(z) = f_A(\#(B)) = f_A(f_B(x_0)) = (f_A \circ f_B)(x_0)$, and so 
we would expect that $\#(A \cup B) =  (f_A \circ f_B)(x_0)$. But if $\oplus$ is the operation given by 
Theorem~\ref{theorem_am_11} then $\#(A \cup B) = \#(A) \oplus \#(B) = f_A(x_0) \oplus f_B(x_0)$, which suggests 
that the following should hold:
\begin{evlist}{26pt}{6pt}
\item[$\mathrm{(\mu)}$]
$f_A(x_0) \oplus f_B(x_0) = (f_A \circ f_B)(x_0)$ whenever $A$ and $B$ are disjoint finite sets.
\end{evlist}
It will be seen later that $\mathrm{(\mu)}$ does hold. What is perhaps more important, though, is that
$\mathrm{(\mu)}$ can actually be used to define $\oplus$, as we now explain.

Denote by $\,\Self{X}$ the set of all mappings from $X$ to itself. Then $(\Self{X},\circ,\id_X)$, where $\circ$ 
is functional composition and $\id_X : X \to X$ is the identity mapping, is a monoid. (A \definition{monoid} is 
any triple $(M,\bullet,e)$ consisting of a set $M$, an associative operation $\bullet$ on $X$ and a unit element 
$e \in M$ satisfying $a \bullet e = e \bullet a = a$ for all $a \in M$.) Lemma~\ref{lemma_mam_61} shows that
\[ 
M_f = \{ u \in \Self{X} : \mbox{$u = f_A$ for some finite set $A$} \}
\]
is a submonoid of $(\Self{X},\circ, \id_X)$, meaning that $\id_X \in M_f$ and $u_1 \circ u_2 \in M_f$ for all 
$u_1,\,u_2 \in M_f$, and that this submonoid is commutative, i.e., $u_1 \circ u_2 = u_2 \circ u_1$ for all 
$u_1,\,u_2 \in M_f$. (The monoid $(\Self{X},\circ,\id_X)$ itself is not commutative except when $X = \{x_0\}$.)

Let $\Phi_{x_0} : M_f \to X$ be the mapping with $\Phi_{x_0}(u) = u(x_0)$ for each $u \in M_f$, and so in 
particular $\Phi_{x_0}(f_A) = f_A(x_0) = \#(A)$ for each finite set $A$. Lemma~\ref{lemma_mam_71} will show that 
$\Phi_{x_0}$ is a bijection, and therefore there exists a unique operation $\oplus$ on $X$ such that
\begin{evlist}{26pt}{6pt}
\item[$\mathrm{(\nu)}$]
$\Phi_{x_0}(u) \oplus \Phi_{x_0}(v) = \Phi_{x_0}(u \circ v)$ for all $u,\,v \in M_f$.
\end{evlist}
This is how $\oplus$ will be defined below. Note that if $A$ and $B$ are (not necessarily disjoint) finite sets 
then by $\mathrm{(\nu)}$
\[ 
f_A(x_0) \oplus f_B(x_0) = \Phi_{x_0}(f_A) \oplus \Phi_{x_0}(f_B) = \Phi_{x_0}(f_A \circ f_B) = (f_A \circ f_B)(x_0)
\]
and so in particular $\mathrm{(\mu)}$ holds.

We now give the details of the approach outlined above.

\begin{lemma}\label{lemma_mam_11}
The assignment $A \mapsto f_A$ satisfies $f_\varnothing = \id_X$ and $f_{A \cup \{a\}} = f \circ f_A$ whenever 
$A$ is a finite set and $a \notin A$. Moreover, it is uniquely determined by these requirements.
\end{lemma}

\proof 
We have $f_\varnothing(x) = \#_x(\varnothing) = x = \id_X(x)$ for all $x \in X$, and so $f_\varnothing = \id_X$.
Moreover, if $A$ is a finite set and $a \notin A$ then
\[ 
f_{A \cup \{a\}}(x) = \#_x(A \cup \{a\}) = f(\#_x(A)) = f(f_A(x)) = (f \circ f_A)(x)
\]
for all $x \in X$ and hence $f_{A \cup \{a\}} = f \circ f_A$. Finally, consider a further assignment 
$A \mapsto f'_A$ with $f'_\varnothing = \id_X$ and such that $f'_{A \cup \{a\}} = f \circ f'_A$ whenever $A$ is a 
finite set and $a \notin A$. For each finite set $A$ let $\prop(A)$ be the proposition that $f'_A = f_A$.

($\diamond$)\enskip 
$\prop(\varnothing)$ holds since $f'_\varnothing = \id_X = f_\varnothing$.

($\star$)\enskip 
Let $A$ be a finite set for which $\prop(A)$ holds (and so $f'_A = f_A$) and let $a \notin A$. Then
$f'_{A \cup \{a\}} = f \circ f'_A = f \circ f_A = f_{A \cup \{a\}}$ and so $\prop(A \cup \{a\})$ holds.

Therefore by the induction principle for finite sets $\prop(A)$ holds for every finite set $A$, which means that 
$f'_A = f_A$ for each finite set $A$.
\eop

The assignment $A \mapsto f_A$ will be called the \definition{$f$-iterator}. Note that in particular 
$f_{\{a\}} = f$ for each element $a$, since 
$f_{\{a\}} = f_{\varnothing \cup \{a\}} = f \circ f_\varnothing = f \circ \id_X = f$.

\begin{lemma}\label{lemma_mam_21}
$f \circ f_A = f_A \circ f$  for each finite set $A$.
\end{lemma}

\proof 
For each finite set $A$ let $\prop(A)$ be the proposition that $f \circ f_A = f_A \circ f$.

($\diamond$)\enskip 
$\prop(\varnothing)$ holds since 
$f \circ f_\varnothing = f \circ \id_X = f = \id_X \circ f = f_\varnothing \circ f$.

($\star$)\enskip 
Let $A$ be a finite set for which $\prop(A)$ holds and let $a \notin A$. Then 
\[ 
f \circ f_{A \cup \{a\}} = f \circ f \circ f_A = f \circ f_A \circ f = f_{A \cup \{a\}} \circ f
\]
and so $\prop(A \cup \{a\})$ holds.

Therefore by the induction principle for finite sets $\prop(A)$ holds for every finite set $A$ and hence
$f \circ f_A = f_A \circ f$  for each finite set $A$.
\eop

The next result establishes an important relationship between the $f$-iterator and the iterator $\#$ for 
$(X,f,x_0)$.

\begin{proposition}\label{prop_mam_11}
If $A$ and $B$ are finite sets then $f_A = f_B$ holds if and only if $\#(A) = \#(B)$.
\end{proposition}

\proof 
By definition $\#(C) = f_C(x_0)$ for each finite set $C$, and so $\#(A) = \#(B)$ whenever $f_A = f_B$. Suppose 
conversely that $\#(A) = \#(B)$ and consider the set $X_0 = \{ x \in X : f_A(x) = f_B(x) \}$. Then $X_0$ is 
$f$-invariant, since if $x \in X_0$ then by Lemma~\ref{lemma_mam_21} 
$f_A(f(x)) = f(f_A(x)) = f(f_B(x)) = f_B(f(x))$, i.e., $f(x) \in X_0$. Also $x_0 \in X_0$, since 
$f_A(x_0) = \#(A) = \#(B) = f_B(x_0)$. Hence $X_0 = X$, since $(X,f,x_0)$ is minimal. This shows that 
$f_A(x) = f_B(x)$ for all $x \in X$, i.e., $f_A = f_B$. 
\eop

There is another way of obtaining the $f$-iterator: Consider the counting system $(\Self{X},f_*,\id_X)$, where 
$f_* : \Self{X} \to \Self{X}$ is defined by $f_*(h) = f \circ h$ for all $h \in \Self{X}$.

\begin{lemma}\label{lemma_mam_31}
If $\#_*$ is the unique iterator for $(\Self{X},f_*,\id_X)$ then $\#_*(A) = f_A$ for each finite set $A$.
\end{lemma}

\proof 
By definition $\#_*(\varnothing) = \id_X$ and $\#_*(A \cup \{a\}) = f_*(\#_*(A)) = f \circ \#_*(A)$ whenever $A$ 
is a finite set and $a \notin A$. Thus by the uniqueness in Lemma~\ref{lemma_mam_11} $\#_*(A) = f_A$ for each 
finite set $A$.
\eop

\begin{proposition}\label{prop_mam_21}
If $A$ and $B$ are finite sets with $A \approx B$ then $f_A = f_B$.
\end{proposition}

\proof 
This follows from Proposition~\ref{prop_cfs_21} (applied to $\#_*$) and Lemma~\ref{lemma_mam_31}. 
\eop

\begin{lemma}\label{lemma_mam_41}
If $A$ and $B$ are disjoint finite sets then $f_{A \cup B} = f_A \circ f_B$.
\end{lemma}

\proof 
For each finite set $A$ let $\prop(A)$ be the proposition that $f_{A \cup B} = f_A \circ f_B$ for each finite 
set $B$ disjoint from $A$.

($\diamond$)\enskip 
$\prop(\varnothing)$ holds since $f_{\varnothing \cup B} = f_B = \id_X \circ f_B = f_\varnothing \circ f_B$ for each 
finite set $B$.

($\star$)\enskip 
Let $A$ be a finite set for which $\prop(A)$ holds and let $a \notin A$. Consider a finite set $B$ disjoint from 
$A \cup \{a\}$; then $B$ is disjoint from $A$ and so $f_{A \cup B} = f_A \circ f_B$. Moreover $a \notin A \cup B$ 
and hence
\[
f_{(A \cup \{a\}) \cup B} = f_{(A \cup B) \cup \{a\}} = f \circ f_{A \cup B} =  f \circ f_A \circ f_B  
= f_{A \cup \{a\}} \circ f_B \;.
\]
This shows that $\prop(A \cup \{a\})$ holds.

Therefore by the induction principle for finite sets $\prop(A)$ holds for every finite set $A$, which means 
$f_{A \cup B} = f_A \circ f_B$ whenever $A$ and $B$ are disjoint finite sets.
\eop

\begin{lemma}\label{lemma_mam_51}
(1)\enskip
If $f$ is bijective then $f_A$ is bijective for each finite set $A$. 

(2)\enskip
If $f$ is injective then $f_A$ is also injective for each finite set $A$. 
\end{lemma}

\proof
(1)\enskip
For each finite set $A$ let $\prop(A)$ be the proposition that $f_A$ is bijective.

($\diamond$)\enskip 
$\prop(\varnothing)$ holds since $f_{\varnothing} = \id_X$ is bijective.

($\star$)\enskip 
Let $A$ be a finite set for which $\prop(A)$ holds and let $a \notin A$. Then $f_{A \cup \{a\}} = f \circ f_A$, as 
the composition of two bijective mappings, is itself bijective, and so $\prop(A \cup \{a\})$ holds.

Therefore by the induction principle for finite sets $\prop(A)$ holds for every finite set $A$, which means that 
$f_A$ is bijective for each finite set $A$.

(2)\enskip 
Just replace `bijective' by `injective' in (1).
\eop

As above let $M_f = \{ u \in \Self{X} : \mbox{$u = f_A$ for some finite set $A$} \}$. Thus in particular 
$\id_X \in M_f$, since $\id_X = f_\varnothing$, and $f \in M_f$, since $f = f_{\{a\}}$ for each element $a$. 
Moreover, if $f$ is injective (resp.\ bijective) then by Lemma~\ref{lemma_mam_51} each element in $M_f$ is 
injective (resp.\ bijective).

\begin{lemma}\label{lemma_mam_61}
For all $u_1,\,u_2 \in M_f$ we have $u_1 \circ u_2 \in M_f$ and $u_1 \circ u_2 = u_2 \circ u_1$. (Since also
$\id_X \in M_f$ this means that $M_f$ is a commutative submonoid of the monoid $(\Self{X},\circ,\id_X)$.)
\end{lemma}

\proof 
Let $u_1,\,u_2 \in M_f$ and so there exist finite sets $A$ and $B$ with $u_1 = f_A$ and $u_2 = f_B$. There then 
exists a disjoint pair $(A',B')$ with $(A',B') \approx (A,B)$ and hence by Proposition~\ref{prop_mam_21} and 
Lemma~\ref{lemma_mam_41}
\[ 
u_1 \circ u_2 = f_A \circ f_B = f_{A'} \circ f_{B'} = f_{A' \cup B'} = f_{B' \cup A'}
= f_{B'} \circ f_{A'} = f_B \circ f_A = u_2 \circ u_1\;, 
\]
i.e., $u_1 \circ u_2 = u_2 \circ u_1$. Moreover, since $u_1 \circ u_2 = f_{A' \cup B'}$ and 
$f_{A' \cup B'} \in M_f$, this also shows that $u_1 \circ u_2 \in M_f$.
\eop

As above let $\Phi_{x_0} : M_f \to X$ be the mapping with $\Phi_{x_0}(u) = u(x_0)$ for all $u \in M_f$. Then 
$\Phi_{x_0}(\id_X) = x_0$ and $\Phi_{x_0}(f_A) = f_A(x_0) = \#(A)$ for each finite set $A$. An important property 
of $\Phi_{x_0}$ is that 
\begin{evlist}{26pt}{6pt}
\item[$\mathrm{(\sharp)}$]
$\ u(\Phi_{x_0}(v)) = \Phi_{x_0}(u \circ v)\,$ for all $u,\,v \in M_f$, 
\end{evlist}
which holds since $u(\Phi_{x_0}(v)) = u(v(x_0)) = (u \circ v)(x_0) = \Phi_{x_0}(u \circ v)$. The special case
of this with $u = f$ gives us $f(\Phi_{x_0}(v)) = \Phi_{x_0}(f \circ v)$ for all $v \in M_f$.

\begin{lemma}\label{lemma_mam_71}
The mapping $\Phi_{x_0}$ is a bijection.
\end{lemma}

\proof 
If $x \in X$ then by Lemma~\ref{lemma_cfs_41} there exists a finite set $A$ with $x = \#(A)$ and it follows that
$\Phi_{x_0}(f_A) = f_A(x_0) = \#(A) = x$. Thus $\Phi_{x_0}$ is surjective. Now let $u_1,\,u_2 \in M_f$ with 
$\Phi_{x_0}(u_1) = \Phi_{x_0}(u_2)$. By the definition of $M_f$ there exist finite sets $A$ and $B$ with 
$u_1 = f_A$ and $u_2 = f_B$, and hence 
\[ 
\#(A) = \Phi_{x_0}(f_A) = \Phi_{x_0}(u_1) = \Phi_{x_0}(u_2)= \Phi_{x_0}(f_B) = \#(B)\;.
\]
Therefore by Proposition~\ref{prop_mam_11} $f_A = f_B$, i.e., $u_1 = u_2$, which shows that $\Phi_{x_0}$ is also 
injective. 
\eop

\textit{Proof of Theorem~\ref{theorem_am_11}:}\enskip
Since $\Phi_{x_0} : M_f \to X$ is a bijection there clearly exists a unique binary relation $\oplus$ on $X$ such 
that
\[ 
\Phi_{x_0}(u_1) \oplus \Phi_{x_0}(u_2) = \Phi_{x_0}(u_1 \circ u_2)
\]
for all $u_1,\,u_2 \in M_f$. The operation $\oplus$ is associative since $\circ$ has this property: 
If $x_1,\,x_2,\,x_3 \in X$ and $u_1,\,u_2,\,u_3 \in M_f$ are such that $x_j = \Phi_{x_0}(u_j)$ for each $j$ then 
\begin{eqnarray*}
(x_1 \oplus x_2) \oplus x_3 &=& (\Phi_{x_0}(u_1) \oplus \Phi_{x_0}(u_2)) \oplus \Phi_{x_0}(u_3) \\
&=& \Phi_{x_0}(u_1 \circ u_2) \oplus \Phi_{x_0}(u_3) = \Phi_{x_0}( (u_1 \circ u_2) \circ u_3)\\ 
&=& \Phi_{x_0}( u_1 \circ (u_2 \circ u_3)) =  \Phi_{x_0}(u_1) \oplus \Phi_{x_0}(u_2 \circ u_3)\\
&=& \Phi_{x_0}(u_1) \oplus (\Phi_{x_0}(u_2) \oplus \Phi_{x_0}(u_3))
=  x_1 \oplus (x_2 \oplus x_3) \;.
\end{eqnarray*}
In the same way $\oplus$ is commutative, since by Lemma~\ref{lemma_mam_61} the restriction of $\circ$ to $M_f$ 
has this property: If $x_1,\,x_2 \in X$ and $u_1,\,u_2 \in M_f$ are such that $x_1 = \Phi_{x_0}(u_1)$ and 
$x_2 = \Phi_{x_0}(u_2)$ then $u_1 \circ u_2 = u_2 \circ u_1$ and so 
\begin{eqnarray*}
x_1 \oplus x_2  &=& \Phi_{x_0}(u_1) \oplus \Phi_{x_0}(u_2)\\ 
&=& \Phi_{x_0}(u_1 \circ u_2) = \Phi_{x_0}(u_2 \circ u_1) = \Phi_{x_0}(u_2) \oplus \Phi_{x_0}(u_1) = x_2 \oplus x_1\;.
\end{eqnarray*}
Moreover, if $x \in X$ and $u \in M_f$ is such that $x = \Phi_{x_0}(u)$ then
\[ 
x \oplus x_0 = \Phi_{x_0}(u) \oplus \Phi_{x_0}(\id_X) = \Phi_{x_0}(u \circ \id_X) = \Phi_{x_0}(u) =  x\,,
\]
and so $x \oplus x_0 = x$ for all $x \in X$.

Let $x_1,\,x_2 \in X$; we next show that for some $x \in X$ either $x_1 = x_2 \oplus x$ or $x_2 = x_1 \oplus x$. 
Let $u_1,\,u_2 \in M_f$ be such that $x_1 = \Phi_{x_0}(u_1)$ and $x_2 = \Phi_{x_0}(u_2)$ and let $A$ and $B$ be 
finite sets with $u_1 = f_A$ and $u_2 = f_B$. By Proposition~\ref{prop_fs_91} there either exists an injective 
mapping $g : A \to B$ or an injective mapping $h : B \to A$. Assume the former holds and put $B' = g(A)$ and 
$C = B \setminus B'$. Then $B'$ and $C$ are disjoint and $B = B' \cup C$; moreover, $A \approx B'$ (since $g$ 
considered as a mapping from $A$ to $B'$ is a bijection) and so by Proposition~\ref{prop_mam_21} $f_A = f_{B'}$. 
Thus, putting $x = \Phi_{x_0}(f_C)$, it follows that
\begin{eqnarray*}
x_2 &=& \Phi_{x_0}(u_2) = \Phi_{x_0}(f_B) = \Phi_{x_0}(f_{B' \cup C}) = \Phi_{x_0}(f_{B'}\circ f_C)\\
&=& \Phi_{x_0}(f_{B'}) \oplus \Phi_{x_0}(f_C) = \Phi_{x_0}(f_A) \oplus \Phi_{x_0}(f_C)
= \Phi_{x_0}(u_1) \oplus x = x_1 \oplus x\:.
\end{eqnarray*}
On the other hand, if there exists an injective mapping $h : B \to A$ then the same argument shows there exists 
$x \in X$ with $x_1 = x_2 \oplus x$.

Now to (a0) and (a1), and we have seen above that (a0) holds. Let $x,\,x' \in X$ and let $u,\,u' \in M_f$ with 
$x = \Phi_{x_0}(u)$ and $x' = \Phi_{x_0}(u')$. Then by $\mathrm{(\sharp)}$
\begin{eqnarray*}
x \oplus f(x') &=& \Phi_{x_0}(u) \oplus f(\Phi_{x_0}(u')) = \Phi_{x_0}(u) \oplus \Phi_{x_0}(f \circ u'))\\
&=& \Phi_{x_0}(u \circ f \circ u') = \Phi_{x_0}(f \circ u \circ u') = f(\Phi_{x_0}(u \circ u'))\\
&=& f(\Phi_{x_0}(u) \oplus \Phi_{x_0}(u')) = f(x \oplus x')
\end{eqnarray*}
and so (a1) holds. If $\oplus'$ is another binary operation on $X$ satisfying (a0) and (a1) then it is easy to 
see that $X_0 = \{ x' \in X : x \oplus' x' = x \oplus x'\ \mbox{for all $x \in X$} \}$ is an $f$-invariant subset 
of $X$ containing $x_0$. Hence $X_0 = X$, since $(X,f,x_0)$  is minimal, which implies that ${\oplus'} = {\oplus}$.

Finally, if $A$ and $B$ are disjoint finite sets then
\[ 
\#(A) \oplus \#(B) = \Phi_{x_0}(f_A) \oplus \Phi_{x_0}(f_B) 
= \Phi_{x_0}(f_A \circ f_B) = \Phi_{x_0}(f_{A\cup B}) = \#(A \cup B)\;.
\] 
Moreover, $\oplus$ is uniquely determined by this requirement: Consider any binary operation $\oplus'$ on $X$ 
for which $\#(A) \oplus' \#(B) = \#(A \cup B)$ whenever $A$ and $B$ are disjoint finite sets. If $C$ and $D$ are 
any finite sets then there exists a disjoint pair $(A,B)$ with $(A,B) \approx (C,D)$ and hence by 
Proposition~\ref{prop_cfs_21} 
\[ 
\#(C) \oplus' \#(D) = \#(A) \oplus' \#(B) = \#(A \cup B) = \#(A) \oplus \#(B) = \#(C) \oplus \#(D) 
\] 
and so by Lemma~\ref{lemma_cfs_41} ${\oplus'} = {\oplus}$. This completes the proof of 
Theorem~\ref{theorem_am_11}. 
\eop

Let $\oplus$ be the operation given in Theorem~\ref{theorem_am_11}. The theorem shows in particular that  
$(X,\oplus,x_0)$ is a commutative monoid. The next result, which generalises Propositions \ref{prop_am_11} and 
\ref{prop_am_21}, shows how properties of the mapping $f$ correspond to properties of the monoid $(X,\oplus,x_0)$.

\begin{proposition}\label{prop_mam_31}
(1)\enskip
$(X,\oplus,x_0)$ is a group if and only if $f$ is a bijection.

(2)\enskip
The cancellation law holds in $(X,\oplus,x_0)$ if and only if $f$ is injective.
\end{proposition}

\proof
We have the commutative monoid $(X,\oplus,x_0)$, and also the commutative monoid $(M_f,\circ,\id_X)$. Now the 
operation $\oplus$ was defined so that
\[ 
\Phi_{x_0}(u_1) \oplus \Phi_{x_0}(u_2) = \Phi_{x_0}(u_1 \circ u_2)
\]
for all $u_1,\,u_2 \in M_f$ and, since  $\Phi_{x_0}(\id_X) = x_0$ and $\Phi_{x_0}$ is a bijection, this means 
$\oplus$ was defined to make $\Phi_{x_0} : (M_f,\circ,\id_X) \to (X,\oplus,x_0)$ a monoid isomorphism. It follows 
that the cancellation law holds in $(X,\oplus,x_0)$ if and only if it holds in $(M_f,\circ,\id_X)$ and that 
$(X,\oplus,x_0)$ will be a group if an only if $(M_f,\circ,\id_X)$ is. It is thus enough to prove the statements 
in the proposition with $(X,\oplus,x_0)$ replaced by $(M_f,\circ,x_0)$.

(1)\enskip
We first show that $u^{-1} \in M_f$ whenever $u \in M_f$ is a bijection. This follows from the fact that 
$u^{-1}(x_0) \in X$ and $\Phi_{x_0}$ is surjective and so there exists $v \in M_f$ with 
$\Phi_{x_0}(v) = u^{-1}(x_0)$; thus by $\mathrm{(\sharp)}$
\[  
\Phi_{x_0}(u \circ v) = u(\Phi_{x_0}(v)) = u(u^{-1}(x_0)) = x_0 = \Phi_{x_0}(\id_X) 
\]
and therefore $u \circ v = \id_X$, since $\Phi_{x_0}$ is injective. Hence $u^{-1} = v \in M_f$. Now clearly $M_f$ 
is a group if and only if  each mapping $u \in M_f$ is a bijection and $u^{-1} \in M_f$, and we have just seen that
$u^{-1} \in M_f$ holds automatically whenever $u \in M_f$ is a bijection. Moreover, by 
Lemma~\ref{lemma_mam_51}~(1) each element of $M_f$ is a bijection if and only if $f$ is a bijection.

(2)\enskip
Suppose  the cancellation law holds in $(M_f,\circ,\id_X)$, and let $x_1,\,x_2 \in X$ with $f(x_1) = f(x_2)$. Then 
there exist $u_1,\,u_2 \in M_f$ with $\Phi_{x_0}(u_1) = x_1$ and $\Phi_{x_0}(u_2) = x_2$ (since $\Phi_{x_0}$ is 
surjective), and hence by $\mathrm{(\sharp)}$
\[   
\Phi_{x_0}(f\circ u_1) = f(\Phi_{x_0}(u_1)) = f(x_1) = f(x_2) = f(\Phi_{x_0}(u_2)) = \Phi_{x_0}(f\circ u_2)\;.
\] 
It follows that $f \circ u_1 = f \circ u_2$ (since $\Phi_{x_0}$ is injective) and so $u_1 = u_2$. In particular 
$x_1 = x_2$, which shows that $f$ is injective. The converse is immediate, since if $f$ is injective then by 
Lemma~\ref{lemma_mam_51}~(2) so is each $u \in M_f$ and hence $u_1 = u_2$ whenever $u \circ u_1 = u \circ u_2$.
\eop

We now begin the preparations for the proof of Theorem~\ref{theorem_am_21}.

For each finite set $B$ there is the mapping $f_B : X \to X$ and thus there exists a unique $f_B$-iterator. This 
is the unique assignment $A \mapsto (f_B)_A$ with $(f_B)_\varnothing = \id_X$ such that
$(f_B)_{A \cup \{a\}} = f_B \circ (f_B)_A$ for each finite set $A$ and each element $a \notin A$.

\begin{lemma}\label{lemma_mam_81}
(1)\enskip $(f_B)_A = f_{B\times A}$ for all finite sets $A$ and $B$.

(2)\enskip $(f_B)_A = (f_A)_B$ for all finite sets $A$ and $B$.
\end{lemma}

\proof
(1)\enskip
Consider the finite set $B$ to be fixed and for each finite set $A$ let $\prop(A)$ be the proposition that 
$(f_B)_A = f_{B\times A}$.

($\diamond$)\enskip 
$\prop(\varnothing)$ holds since $(f_B)_{\varnothing} = \id_X = f_\varnothing = f_{B \times \varnothing}$.

($\star$)\enskip 
Let $A$ be a finite set for which $\prop(A)$ holds (and so $(f_B)_A = f_{B\times A}$) and let $a \notin A$. 
Then $B \times (A \cup \{a\})$ is the disjoint union of the sets $B \times \{a\}$ and $B \times A$ and 
$B \times \{a\} \approx B$; thus by Proposition~\ref{prop_mam_21} and Lemma~\ref{lemma_mam_41}
\[  
(f_B)_{A \cup \{a\}} = f_B \circ (f_B)_A = f_B \circ f_{B\times A} = f_{B\times \{a\}} \circ f_{B\times A} 
= f_{(B\times \{a\}) \cup (B\times A)} = f_{B \times (A\cup \{a\})}
\]
and so $\prop(A \cup \{a\})$ holds.

Therefore by the induction principle for finite sets $\prop(A)$ holds for every finite set $A$, which means  
$(f_B)_A = f_{B\times A}$ for all finite sets $A$ and $B$.

(2)\enskip
By Proposition~\ref{prop_mam_21} $f_{B\times A} = f_{A\times B}$, since clearly $B \times A \approx A \times B$, 
and therefore by (1) $(f_B)_A = f_{B\times A} = f_{A \times B} = (f_A)_B$. 
\eop

The next result will not be needed in what follows, but it shows that $\mathrm{(\delta)}$ in Section~\ref{am}
holds, and could be used instead of Lemma~\ref{lemma_am_31} in the previous proof of Theorem~\ref{theorem_am_21}.

\begin{lemma}\label{lemma_mam_91}
Let $(A,B)$ and $(A',B')$ be pairs of finite sets.

(1)\enskip
If $f_A = f_{A'}$ and $f_B = f_{B'}$ then $f_{A \times B} = f_{A' \times B'}$.

(2)\enskip
If $\#(A,B) = \#(A',B')$ then $\#(A \times B) = \#(A' \times B')$.
\end{lemma}

\proof 
(1)\enskip
By several applications of Lemma~\ref{lemma_mam_81} (1) and (2) we have
\[
f_{A \times B} = (f_B)_A = (f_{B'})_A = (f_A)_{B'} = (f_{A'})_{B'} = (f_{B'})_{A'} = f_{A' \times B'}\;.
\]

(2)\enskip
This follows immediately from (1) and Proposition~\ref{prop_mam_11}. \eop

\begin{lemma}\label{lemma_mam_101}
Let $v \in M_f$ and $A \to v_A$ be the unique $v$-iterator (and so $v_\varnothing = \id_X$ and
$v_{A \cup \{a\}} = v \circ v_A$ for each finite set $A$ and each element $a \notin A$). Then:

(1)\enskip
$v_A \in M_f$ for every finite set $A$.

(2)\enskip
If $A$ and $B$ are finite sets with $f_A = f_B$ then $v_A = v_B$.
\end{lemma}

\proof
(1)\enskip
For each finite set $A$ let $\prop(A)$ be the proposition that $v_A \in M_f$.

($\diamond$)\enskip 
$\prop(\varnothing)$ holds since $v_{\varnothing} = \id_X \in M_f$.

($\star$)\enskip Let $A$ be a finite set for which $\prop(A)$ holds (and so $v_A \in M_f$) and let $a \notin A$. 
Then $v_{A \cup \{a\}} = v \circ v_A \in M_f$, since $M_f$ is a submonoid of $(\Self{X},\circ,\id_X)$, and so 
$\prop(A \cup \{a\})$ holds.

Therefore by the induction principle for finite sets $\prop(A)$ holds for every finite set $A$, which means  
$v_A \in M_f$ for all finite sets $A$.

(2)\enskip
There exists a finite set $C$ with $v = f_C$ and thus by Lemma~\ref{lemma_mam_81}~(2)
\[ 
v_A = (f_C)_A = (f_A)_C = (f_B)_C = (f_C)_B = v_B\;.\ \eop
\]

Let $v \in M_f$; then by Lemma~\ref{lemma_mam_101} there exists a unique mapping $\psi_v : M_f \to M_f$ such 
that $\psi_v(f_A) = v_A$ for each finite set $A$, and in particular $\psi_v(f) = v$ (since if $a$ is any element 
then $f = f_{\{a\}}$ and $v = v_{\{a\}}$). Moreover, if $v = f_B$ then by Lemma~\ref{lemma_mam_81}~(1) 
$\psi_v(f_A) = f_{A \times B}$. 

A mapping $\psi : M_f \to M_f$ is an \definition{endomorphism} (of the monoid $(M_f,\circ,\id_X)$) if 
$\psi(\id_X) = \id_X$ and $\psi(u_1 \circ u_2) = \psi(u_1) \circ \psi(u_2)$ for all $u_1,\,u_2 \in M_f$.

\begin{lemma}\label{lemma_mam_111}
(1)\enskip
$\psi_v$ is an endomorphism for each $v \in M_f$. 

(2)\enskip
$\psi_v(u) = \psi_u(v)$ for all $u,\,v \in M_f$.
\end{lemma}

\proof
(1)\enskip
If $u_1,\,u_2 \in M_f$ then there exist disjoint finite sets $A$ and $B$ with $u_1 = f_A$ and $u_2 = f_B$ and 
hence by Lemma~\ref{lemma_mam_41}
\[ 
\psi_v(u_1 \circ u_2) = \psi_v(f_A \circ f_B) = \psi_v(f_{A \cup B}) = v_{A \cup B} = 
  v_A \circ v_B = \psi_v(f_A) \circ \psi_v(f_B) 
\] 
(noting that Lemma~\ref{lemma_mam_41} can also be applied to the $v$-iterator). Moreover, we have
$\psi_v(\id_X) = \psi_v(f_\varnothing) = v_\varnothing = \id_X$ and hence $\psi_v$ is an endomorphism.

(2)\enskip
Let $A,\,B$ be finite sets with $u = f_A$ and $v = f_B$. Then by Lemma~\ref{lemma_mam_81}~(2)
\[ 
\psi_v(u) = \psi_v(f_A) = v_A = (f_B)_A = (f_A)_B = u_B = \psi_u(f_B) = \psi_u(v)\;.\ \eop 
\]

\textit{Proof of Theorem~\ref{theorem_am_21}:}\enskip
Define a binary operation $\diamond : M_f \times M_f \to M_f$ by letting
\[ 
u \diamond v = \psi_u(v) 
\] 
for all $u,\,v \in M_f$. In particular, if $A$ and $B$ are finite sets then by Lemma~\ref{lemma_mam_81}~(1) 
$f_A \diamond f_B = \psi_{f_A}(f_B) = (f_A)_B = f_{A \times B}$ and therefore
\[ 
f_A \diamond f_B = f_{A \times B} 
\]
for all finite sets $A$ and $B$. Let $u,\,v,\,w \in M_f$ and $A$, $B$ and $C$ be finite sets with $u = f_A$, 
$v = f_B$ and $w = f_C$. Then clearly $A \times (B \times C) \approx (A \times B) \times C$, so by 
Proposition~\ref{prop_mam_21} $f_{A \times (B \times C)} = f_{(A \times B) \times C}$ and thus
\begin{eqnarray*}
u \diamond (v \diamond w) &=& f_A \diamond (f_B \diamond f_C) = f_A \diamond f_{B\times C}\\ 
&=& f_{A\times (B\times C)} = f_{(A\times B)\times C} = f_{A\times B} \diamond f_C 
= (f_A \diamond f_B) \diamond f_C = (u \diamond v) \diamond w\;.
\end{eqnarray*}
Hence $\diamond$ is associative. Moreover, Lemma~\ref{lemma_mam_111}~(2) shows that $\diamond$ is commutative, 
since $u \diamond v = \psi_u(v) = \psi_v(u) = v \diamond u$ for all $u,\,v \in M_f$. Also (with $a$ any element) 
$u \diamond f = u \diamond f_{\{a\}} = u_{\{a\}} = u$, i.e., $u \diamond f = u$ for all $u \in M_f$, and by 
Lemma~\ref{lemma_mam_111}~(1) $u \diamond \id_X = \id_X$ and 
$u \diamond (v_1 \circ v_2) = (u \diamond v_1) \circ (u \diamond v_2)$ for all $u,\,v_1,\,v_2 \in M_f$.

Now since $\Phi_{x_0} : M_f \to X$ is a bijection there clearly exists a unique binary relation $\otimes$ on 
$X$ such that
\[ 
\Phi_{x_0}(u_1) \otimes \Phi_{x_0}(u_2) = \Phi_{x_0}(u_1 \diamond u_2)
\]
for all $u_1,\,u_2 \in M_f$, and exactly as in the proof of Theorem~\ref{theorem_am_11} the operation $\otimes$ 
is associative and commutative since $\diamond$ has these properties. The same holds true of the distributive 
law: Let $x,\,x_1,\,x_2 \in X$, and $u,\,v_1,\,v_2 \in M_f$ be such that $x = \Phi_{x_0}(u)$, 
$x_1 = \Phi_{x_0}(v_1)$ and $x_2 = \Phi_{x_0}(v_2)$. Then
\begin{eqnarray*}
x \otimes (x_1 \oplus x_2) &=& \Phi_{x_0}(u) \otimes (\Phi_{x_0}(v_1) \oplus \Phi_{x_0}(v_2)) \\
&=& \Phi_{x_0}(u) \otimes \Phi_{x_0}(v_1 \circ v_2) = \Phi_{x_0}(u \diamond (v_1 \circ v_2)) \\
&=& \Phi_{x_0}((u \diamond v_1) \circ (u \diamond v_2))
= \Phi_{x_0}(u \diamond v_1) \oplus \Phi_{x_0}(u \diamond v_2)\\ 
&=& (\Phi_{x_0}(u) \otimes \Phi_{x_0}(v_1)) \oplus (\Phi_{x_0}(u) \otimes \Phi_{x_0}(v_2)) 
=  (x \otimes x_1) \oplus (x \otimes x_2)
\end{eqnarray*}

Next, if $x \in X$ and $u \in M_f$ is such that $x = \Phi_{x_0}(u)$ then
\begin{eqnarray*} 
&&x \otimes x_0 = \Phi_{x_0}(u) \otimes \Phi_{x_0}(\id_X) = \Phi_{x_0}(u \diamond \id_X) = \Phi_{x_0}(\id_X) 
=  x_0\,,\\
&&x \otimes f(x_0) = \Phi_{x_0}(u) \otimes \Phi_{x_0}(f) = \Phi_{x_0}(u \diamond f) = \Phi_{x_0}(u) =  x
\end{eqnarray*}
and so $x \otimes x_0 = x_0$ and $x \otimes f(x_0) = x$ for all $x \in X$.

We have already seen  (m0) holds and, since $f(x_0)$ is a unit, (m1) is a special case of the distributive law: 
Let $x,\,x' \in X$; then by (a0) and (a1) and since $\oplus$ is commutative it follows that
$f(x') = f(x' \oplus x_0) = x' \oplus f(x_0) = f(x_0) \oplus x'$, and hence
$x \otimes f(x') = x \otimes (f(x_0) \oplus x')
= (x \otimes f(x_0)) \oplus (x \otimes x') = x \oplus (x \otimes x')$, which is (m1). Finally, if $\otimes'$ is 
another binary operation satisfying (m0) and (m1) then it is easy to see that 
$X_0 = \{ x' \in X : x \otimes' x' = x \otimes x'\ \mbox{for all $x \in X$} \}$
is an $f$-invariant subset of $X$ containing $x_0$. Hence $X_0 = X$, since $(X,f,x_0)$ is minimal, which implies 
that ${\otimes'} = {\otimes}$.
\eop


\sbox{\ttt}{\textsc{References}}
\thispagestyle{plain}

\bigskip
\bigskip


{\sc Fakult\"at f\"ur Mathematik, Universit\"at Bielefeld}\\
{\sc Postfach 100131, 33501 Bielefeld, Germany}\\
\textit{E-mail address:} \texttt{preston@math.uni-bielefeld.de}\\
\textit{URL:} \texttt{http://www.math.uni-bielefeld.de/\symbol{126}preston}


\end{document}